\documentclass[3p,times]{elsarticle}

\usepackage{amsmath,amssymb}

\biboptions{sort&compress}

\usepackage{url}

\newcommand{\ii}{\mathrm{i}}

\newcommand{\dt}{\ensuremath{{\delta t}}}

\journal{Journal of Computational Physics}

\begin{document}

\begin{frontmatter}

\title{High-order commutator-free exponential time-propagation \\ of driven quantum systems}

\author[aa]{A. Alvermann\corref{cor1}}
\ead{aha26@cam.ac.uk}
\author[hf]{H. Fehske}

\cortext[cor1]{Corresponding author.}
\address[aa]{Theory of Condensed Matter, Cavendish Laboratory, Cambridge CB3 0HE, United Kingdom}
\address[hf]{Institut f\"ur Physik, Ernst-Moritz-Arndt-Universit\"at, 17487 Greifswald, Germany}

\begin{abstract}
We discuss the numerical solution of the Schr{\"o}dinger equation with a time-dependent Hamilton operator using commutator-free time-propagators.
These propagators are constructed as products of exponentials of simple weighted sums of the Hamilton operator. Owing to their exponential form they strictly preserve the unitarity of time-propagation. 
The absence of commutators or other computationally involved operations allows for straightforward implementation and application also to large-scale and sparse matrix problems.
We explain the derivation of commutator-free exponential time-propagators in the context of the Magnus expansion, and provide optimized propagators up to order eight. 
An extensive theoretical error analysis is presented together with practical efficiency tests for different problems. 
Issues of practical implementation, in particular the use of the Krylov technique for the calculation of exponentials, are discussed. 
We demonstrate for two advanced examples, the hydrogen atom in an electric field and pumped systems of multiple interacting two-level systems or spins that this approach enables fast and accurate computations.
\end{abstract}

\begin{keyword}

 time-dependent Schr{\"o}dinger equation \sep  Magnus expansion \sep driven quantum systems \sep Lie group integrators

\end{keyword}

\end{frontmatter}


\section{Introduction}
\label{sec:Intro}

The time-evolution of a driven quantum system is determined by the Schr{\"o}dinger equation
\begin{equation}
 \ii \partial_t \psi(t) = H(t) \psi(t)  
\end{equation}
with a time-dependent Hamilton operator $H(t)$, which one tries to solve
for a given initial wave function $\psi(t_0)$ and times $t\ge t_0$. 
Prominent examples are atoms in laser fields, spins in magnetic fields or quantum dots contacted to AC voltage sources (see e.g. Ref.~\cite{Hae97} for an introductory discussion). 
The Schr{\"o}dinger equation is a special case of a general linear differential equation 
\begin{equation}\label{DGL}
 \partial_t x(t) = A(t) x(t)
\end{equation}
with time-dependent coefficients, where $A(t) = - \ii H(t)$.
Other examples from quantum mechanics are the Liouville-von-Neumann equation for the density operator $\rho(t)$ or master equations for dissipative systems~\cite{BP02}.
Analytical solutions of such equations can be found only in a very limited number of cases. 
In most situations one must resort to numerical computations.
In the present paper we study an efficient numerical solution technique,
which is related to the Magnus expansion but avoids the use of commutators.

The propagator $U(t_1,t_2)$ of Eq.~\eqref{DGL} satisfies the initial value problem
\begin{equation}\label{DefU}
  \partial_{t_1} U(t_1,t_2) = A(t_1) U(t_1,t_2) \;, \quad U(t_0,t_0)= I \;,
\end{equation}
with the identity operator, or matrix, $I$. 
The solutions $x(t)$ of Eq.~\eqref{DGL} fulfill $x(t_1)=U(t_1,t_2) x(t_2)$.
We note the group property $U(t_1,t_2) U(t_2,t_3) = U(t_1,t_3)$.

For time-independent $A \equiv A(t)$ the propagator is given by a (matrix) exponential
\begin{equation}\label{UasExp}
 U(t_1-t_2) \equiv U(t_1,t_2) = \exp \, [ (t_1-t_2) A ] \;.
\end{equation}	
The generalization of this expression for time-dependent $A(t)$ is due to W.~Magnus~\cite{Ma54}.
The Magnus expansion (we refer the reader to the recent review~\cite{BCOR09}) expresses the propagator in the form
\begin{equation}\label{Magn1}
 U(t) \equiv U(t,0) = \exp [\Omega(t)]	\;.
\end{equation}	
Notice that we often set the initial time $t_0=0$. Expressions for arbitrary initial time $t_0$ are obtained by the variable substitution $t \mapsto t+t_0$.
The operator $\Omega(t)$ is given as an infinite series
\begin{equation}\label{Magn2}
 \Omega(t) = \int_0^t dt_1 A(t_1) + \frac{1}{2}  \int_0^t dt_1 \! \int_0^{t_1} dt_2 \, [A(t_1),A(t_2)] + \dots
\end{equation}	
involving nested commutators of $A(t)$ at different times.
Only if $[A(t_1),A(t_2)]=0$ for all $t_1$, $t_2$,
Eqs.~\eqref{Magn1},~\eqref{Magn2} reduce to the simpler expression Eq.~\eqref{UasExp}.
Otherwise the nested commutators provide the necessary correction terms.

The Magnus expansion is important from a theoretical and practical point of view.
In many cases the differential equation Eq.~\eqref{DGL} has an underlying Lie group structure,
where the propagator $U(t_1,t_2)$ is element of a Lie group and $A(t)$ of the associated Lie algebra. For the Schr{\"o}dinger equation, skew-hermiticity of $A(t) = -\ii H(t)$ implies unitarity of $U(t_1,t_2)$.
Violating unitarity results leads to artificial decay or growth of relevant components of the wave function, which spoils the stability of numerical time propagation.
In particular, only unitary propagators preserve the normalization of the wave function.
The Magnus expansion respects the Lie group structure,
since the exponential function $\Omega \mapsto \exp[\Omega]$ maps $\Omega(t)$, which as a sum of commutators of $A(t)$ is itself a Lie algebra element, onto a Lie group element $U(t)$.

From the practical point of view, a truncation of the infinite Magnus expansion provides an approximate propagator $\tilde{U}(t+\dt,t)$, which can be used to propagate a solution $x(t)$ over a small time-step $\dt$.
For an $N$th-order approximation, where the approximation error scales as $\dt^{N+1}$, all terms with $N$ or less commutators in Eq.~\eqref{Magn2} must be kept.
It is the virtue of the Magnus expansion that for every truncation $\tilde{U}(t+\dt,t)$ is a Lie group element (whenever a Lie group structure is present).
In this way the Magnus expansion allows for the systematic construction of geometric integrators~\cite{IN99,IMNZ00,HLW06}, which preserve Lie group structures.

The practical evaluation of the Magnus expansion is however rather involved. 
The number of terms in $\Omega(t)$ is large already for moderate approximation order, and their calculation is complicated because of the nested commutators.
Our starting point for better numerical algorithms are approximations of the form
\begin{equation}\label{UasCFree}
  \tilde{U}(t+\dt,t) = e^{A_1} e^{A_2} \cdot \dots \cdot e^{A_s} \;,
\end{equation}
where each $A_i = \sum_n g_{i,n} A(t_n)$ is a (finite) linear combination of $A(t)$ at different times $t_n \in [t,t+\dt]$ (which will later be chosen as Gauss-Legendre quadrature points).
Such commutator-free exponential time-propagators (CFETs) preserve Lie group structures through the exponential form of the approximation but avoid the use of commutators.
Their application is thus straightforward and requires only slight adjustments of existing programs for the calculation of matrix exponentials.
No complicated scheme for the computation of nested commutators or the storage of intermediate results is needed. CFETs are examples for Crouch-Grossman methods~\cite{HLW06},
and have been studied with a focus on linear differential equations in Refs.~\cite{BM06,Thal06}.
In particular the work of Blanes and Moan~\cite{BM06}, together with the review~\cite{BCOR09}, provided the initial motivation for the work reported here.


In the present paper we discuss CFETs from a practitioner's point of view.
Our intention is to provide a comprehensive account of the theoretical background
and a demonstration of the practical usefulness of this approach.
A specific goal is the construction of optimized high-order CFETs, which can be applied to the Schr{\"o}dinger equation in general situations where the resource consumption of naive computational approaches, e.g. a second-order approximation, would be intolerably large.
To pursue these goals we first revisit the derivation of the Magnus expansion (Secs.~\ref{sec:Magnus},~\ref{sec:Prop}) and of the order conditions for the CFET coefficients (Sec.~\ref{sec:CF}).
A notable deviation from the literature is the replacement of a power series expansion of $A(t)$ with an expansion in Legendre polynomials.
Their orthogonality properties allow to simplify the presentation in two important aspects.
First, the rather non-obvious fact that, effectively, only terms of order $\dt^{N/2}$ of $A(t)$ must be taken into account for the construction of $N$th-order approximations is evident from the structure of the order conditions.
Second, the application of Gauss-Legendre quadrature (Sec.~\ref{sec:GaussLeg}) is straightforward, and the corresponding coefficients are obtained without additional work.
We believe that our presentation is not only simpler than others in the literature,
but allows the reader to understand the derivation without taking unexplained aspects for granted.

Extending previous results we construct CFETs up to order $8$.
Their error is analyzed theoretically in Sec.~\ref{sec:ErrorAnalysis},
complemented by a practical error analysis in Sec.~\ref{sec:Error}.
Minimization of the CFET error requires inclusion of higher order terms from the Magnus expansion, specifically of the $N/2+1$-order term of $A(t)$ for an $N$th-order approximation.
This in contrast to the error analysis for split-operator techniques found in the literature.
Our improved analysis leads to optimized 4th- and 6th-order CFETs.
Again, the use of Legendre polynomials is vital for the analysis.

In practical applications with large Hamiltonian matrices the evaluation of the exponentials in Eq.~\eqref{UasCFree} is the determining factor for the actual efficiency.
We discuss the combination of CFETs with the Krylov technique in Secs.~\ref{sec:Implement},~\ref{sec:ParaOsci}. In Sec.~\ref{sec:Floquet} we compare CFETs with the $(t,t')$-method, a Floquet-based approach. 
Finally, we demonstrate in Sec.~\ref{sec:App} the application of CFETs in two situations
where precise results are hard to obtain otherwise, e.g. with the original Magnus expansion,
before we conclude in Sec.~\ref{sec:Conc}.
The appendices give the recursion for the Magnus expansion in a form suitable for computer algebra computations, a short discussion of free Lie algebras and Hall bases, and the explicit solution of the order conditions for 6th-order CFETs.


\section{The Magnus expansion}
\label{sec:Magnus}

The Magnus expansion provides $\Omega(t)$ in Eq.~\eqref{Magn1} as a series
\begin{equation}\label{Omsum}
 \Omega(t) = \sum_{n=1}^\infty \Omega_n (t) \;,
 \end{equation}	
where $\Omega_n(t)$ is the $n$-fold integral of a sum of $n-1$-fold nested commutators of $A(t)$.
We say that a function $f(t)$ is of order $t^N$ if $\lim_{t \to 0} f(t)/t^{N-1} = 0$, i.e. its power series in $t$ starts with $t^N$.
Since each integration over $t$ increases the order by one, the term $\Omega_n(t)$ is of order $t^n$.
Derivations of the Magnus expansion can be found at many places in the literature (cf. Ref.~\cite{BCOR09}). For our presentation, we follow Ref.~\cite{PL97}.
The principal idea is to find an implicit equation relating $A(t)$ with $\Omega(t)$,
which is solved order by order for the $\Omega_n(t)$.
Notice that we always assume that $A(t)$, and the solutions of Eq.~\eqref{DGL}, are sufficiently regular to permit a local power series expansion.

\subsection{Derivation}
By definition (Eqs.~\eqref{DefU},~\eqref{Magn1}), $\Omega(t)$ is the solution of the implicit differential equation
\begin{equation}\label{DGLOm}
\partial_t e^{\Omega(t)} = A(t) e^{\Omega(t)} \;, \quad \Omega(0)=0  \;.
\end{equation}
To evaluate the derivative of the matrix exponential on the left hand side,
consider the function $f(s,t) =  \partial_t e^{ s \Omega(t) }$.
It fulfills the differential equation
\begin{equation}
 \partial_s f(s,t) = \partial_t \partial_s e^{s \Omega(t)} = 
  \partial_t   e^{s \Omega(t)} \Omega(t) =
   f(s,t) \Omega(t) +  e^{s \Omega(t)} \dot{\Omega}(t)
\end{equation}
with initial condition $f(0,t)=0$,
whose solution is given by 
$f(s,t)=\int_0^s e^{r \Omega(t)} \, \dot{\Omega}(t) \,e^{(s-r) \Omega(t)} dr$.
For $s=1$, we obtain 
\begin{equation}\label{Dexp}
 \partial_t e^{\Omega(t)} = \left(\int_0^1 e^{r \Omega(t)} \dot{\Omega}(t) e^{-r \Omega(t)} dr \right) \, e^{\Omega(t)} = \sum_{m=0}^\infty \left( \int_0^1 \frac{r^m}{m!} [\Omega(t),\dot{\Omega}(t)]_m dr \right) \, e^{\Omega(t)} = \sum_{m=0}^\infty \frac{1}{(m+1)!} [\Omega(t),\dot{\Omega}(t)]_m \, e^{\Omega(t)} \;,
\end{equation}
where we used the identity $e^X Y e^{-X} = \sum_{m=0}^\infty (1/m!) [X,Y]_m$ with
the iterated commutators  
\begin{equation}\label{itComm}
  [X,Y]_0 = Y \,, \; [X,Y]_1 = [X,Y] \,, \; [X,Y]_{m+1} = [X,[X,Y]_m] \;,
\end{equation}
which follows, e.g., from comparison of the derivatives of
$s \mapsto e^{sX} Y e^{-sX}$
and
$s \mapsto \sum_{m=0}^\infty (s^m/m!) [X,Y]_m$.
If the $\Omega(t)$ at different $t$ commute, only the first term $m=0$ in the sum contributes.
Using Eq.~\eqref{Dexp} in Eq.~\eqref{DGLOm} gives
\begin{equation}\label{AandOm}
  A(t) =  (\partial_t e^{\Omega(t)}) \, e^{-\Omega(t)}  = \sum_{m=0}^\infty \frac{1}{(m+1)!} [\Omega(t),\dot{\Omega}(t)]_m \;.
\end{equation}

We now insert the ansatz for the Magnus series Eq.~\eqref{Omsum} into Eq.~\eqref{AandOm}.
This gives
(we drop the argument $t$ in $\Omega(t)$)
\begin{equation}
\begin{split}
 A(t) &= \sum_{m=0}^\infty \frac{1}{(m+1)!}  [\sum_{n=1}^\infty \Omega_n, \sum_{k=1}^\infty \dot{\Omega}_k ]_m \\
  &= \sum_{m=0}^\infty \frac{1}{(m+1)!} \sum_{n_1,\dots,n_m=1}^\infty \sum_{k=1}^\infty [\Omega_{n_1},[\Omega_{n_2},\dots,[\Omega_{n_m},\dot{\Omega}_k] \dots ] \;.
 \end{split}
\end{equation}

To solve for $\dot{\Omega}_n$ we collect all terms of order $t^{n-1}$.
A nested commutator $[\Omega_{n_1},[\Omega_{n_2},\dots,[\Omega_{n_m},\dot{\Omega}_k] \dots ]$ is of order ${n_1+\dots+n_m+k-1}$ in $t$.
The only $(n-1)$th order term that contains $\dot{\Omega}_n$ is the term with $m=0$.
Thus, 
\begin{equation} 
\dot{\Omega}_1 = A \;, \quad \dot{\Omega}_n = - \sum_{m=1}^{n-1} \frac{1}{(m+1)!} \sum_{\substack{n_1, \dots, n_m,k \ge 1 \\ n_1 + \dots + n_m + k = n}} [\Omega_{n_1},[\Omega_{n_2},\dots,[\Omega_{n_m},\dot{\Omega}_k] \dots ] \;.
\end{equation}
Notice that all sums are finite (the last term in the sum over $m$, for $m=n-1$, is $(1/n!) [\Omega_1,\dot{\Omega}_1]_{n-1}$).
A final integration gives the explicit expressions
\begin{equation}\label{OmRec}
\Omega_1(t) = \int\limits_0^t A(t') dt' \;,\quad
\Omega_n (t)  =  - \sum_{m=1}^{n-1} \frac{1}{(m+1)!} \sum_{\substack{n_1, \dots, n_m,k \ge 1 \\ n_1 + \dots + n_m + k = n}} \int\limits_0^t dt' [\Omega_{n_1}(t'),[\Omega_{n_2}(t'),\dots,[\Omega_{n_m}(t'),\dot{\Omega}_k(t')] \dots ] \;,
\end{equation}
which allow for the recursive calculation of the $\Omega_n(t)$.
As stated before, every term in $\Omega_n(t)$ involves an $n$-fold integral of an $n-1$-fold nested commutator of $A(t)$.
We obtain explicitly, up to order $t^3$,
\begin{equation}\label{MagnusOm}
\begin{split}
 \Omega(t) &= \Omega_1(t) + \Omega_2(t) + \Omega_3(t) + O(t^4) \;, \quad\text{with} \\ 
 \Omega_1(t) &= \int_0^t dt_1 A(t_1) \;, \quad
 \Omega_2(t) = \frac{1}{2} \int_0^t dt_1 \int_0^{t_1} dt_2 [A(t_1),A(t_2)] \;, \\
 \Omega_3(t) &= \frac{1}{6} \int_0^t dt_1 \int_0^{t_1} dt_2 \int_0^{t_2} dt_3 [A(t_1),[A(t_2),A(t_3)]] + [[A(t_1),A(t_2)],A(t_3)] \;.
 \end{split}
\end{equation}

It is convenient to write the $\Omega_n(t)$ as time-ordered integrals~\cite{PL97}.
This requires additional manipulation of the integration domains of the terms found by straightforward integration in Eq.~\eqref{OmRec}. It is possible to derive a systematic recursion (see~\ref{app:RecMag}), which is very useful for symbolic calculations on a computer.

An alternative route to solve Eq.~\eqref{AandOm} is to note that the commutator expression on the right hand side involves the Taylor expansion of the function $(e^x-1)/x = \sum_{m=0}^\infty x^m/(m+1)!$. Solving for $\Omega(t)$ is thus possible using the inverse function, where the Bernoulli numbers $B_n$ appear as the Taylor coefficients in $x/(e^x-1) = \sum_{n=0}^\infty B_n \frac{x^n}{n!}$.
After a few additional manipulations one obtains again a recursive definition of the $\Omega_n(t)$  (see e.g. Ref~\cite{BCOR09}).
Our experience is that the present approach is better suited for an algorithmic implementation.
Interestingly enough, it avoids the use of Bernoulli numbers.

\subsection{The Baker-Campbell-Hausdorff formula}\label{sec:BCH}

A special case of the Magnus expansion is the Baker-Campbell-Hausdorff (BCH) formula 
\begin{equation}\label{BCH1}
\begin{split}
e^X e^Y = \exp \Big[ 
X  + Y +\frac{1}{2} [X,Y] +\frac{1}{12} [X,[X,Y]] -\frac{1}{12} [Y,[X,Y]] -\frac{1}{24} [Y,[X,[X,Y]]] 
+ \dots \Big] \;.
\end{split}
\end{equation}
We note that the left hand side of this equation is the exact propagator for a stepwise constant $A(t)$,
with $A(t) =Y$ for $0 \le t < 1$, $A(t) =X$ for $1 \le t \le 2$. 
Inserting this $A(t)$ into the recursion Eq.~\eqref{OmRec} for $\Omega(t)$ provides the exponential on the right hand side.
In a similar spirit, 
we can obtain the BCH formula for several exponentials
\begin{equation}\label{BCH2}
\begin{split}
e^{X_1} \cdots e^{X_s} = \exp \Big[ 
 \sum_{i=1}^s X_i + \frac{1}{2} \sum_{1 \le i<j \le s} [X_i,X_j] + \dots
 \Big] \;.
\end{split}
\end{equation}

\section{Approximate Magnus propagators}

\label{sec:Prop}

By construction, the Magnus expansion is an expansion in orders of $t$.
It thus provides a systematic way to obtain $N$th-order approximations $\tilde{U}^{(N)}(t) =\exp[ \sum_{n=1}^N \Omega_n (t) ]$, which coincide with the exact propagator $U(t)$ for all terms of order $t^N$ or less, from direct truncation of the infinite series Eq.~\eqref{Omsum}.
Notice that we call a function $f(t)$ an $N$th-order approximation of another function $g(t)$ if the difference $f(t)-g(t)$ is of order $t^{N+1}$.
 
The expression for $\Omega_n(t)$, given through Eq.~\eqref{OmRec}, involves $n$-fold integrals. 
These can be simplified since each integral needs to be evaluated only up to order $t^{N+1}$
for $N$th-order approximations.
Starting from an expansion of $A(t)$ in orders of $t$, 
all multi-dimensional integrals in Eq.~\eqref{MagnusOm} can be replaced by one-dimensional integrals. 
In the literature, it is common to expand $A(t)$ in powers of $t$ (or centered powers $(t-\delta t/2)^n$ for a given time-step \dt).
Contrary to these treatments, we use an expansion in Legendre polynomials. Although both expansion are principally equivalent, the choice of Legendre polynomials proves itself useful because of their orthogonality properties.

\subsection{Legendre expansion of $A(t)$}

The (shifted) Legendre polynomials $P_n(x)$ are
defined for $n=0,1,2,\dots$ through the recurrence
\begin{equation}
  P_0(x) = 1 \;, P_1(x) = 2x-1 \;, P_{n+1}(x) = \frac{2n+1}{n+1} (2x-1) P_n(x) - \frac{n}{n+1} P_{n-1}(x)  \;.
\end{equation}
By definition, $P_n(x)$ is a polynomial of degree $n$.
Explicitly,
\begin{equation}
  P_2(x) = 6x^2-6x+1 \;, P_3(x)=20x^3-30x^2+12x-1 \:,
  P_4(x) = 70 x^4 - 140 x^3 + 90 x^2 - 20 x + 1 \;.
\end{equation}
The polynomials $P_n(x)$ are symmetric with respect to $x=1/2$, i.e.
\begin{equation}\label{PSymm}
  P_n(1-x)= (-1)^n P_n (x) \;.
\end{equation}
Furthermore, they form a complete set of orthogonal functions on the interval $[0,1]$, with scalar product
\begin{equation}
 \int\limits_0^1 P_m(x) P_n(x) \,dx = \frac{1}{2n+1} \delta_{mn} \;.
\end{equation}
In particular, $\int_0^1 p(x) P_n(x) \,dx =0$ for every polynomial $p(x)$ of degree less than $n$.

We now fix a time-step $\dt$, for which an approximate $N$th-order propagator $\tilde{U}^{(N)}(\dt) \equiv \tilde{U}^{(N)}(\dt,0)$ should be constructed.
The function $A(t)$ is expanded on the interval $[0,\dt]$ in a series of Legendre polynomials
\begin{equation}\label{Aexpansion}
 A(t) = \frac{1}{\dt} \sum_{n=1}^{N} A_n P_{n-1}\Big(\frac{t}{\dt}\Big) + O(\dt^{N+1}) \qquad (0 \le t \le \dt) \;.
\end{equation}
Notice the index shift of $P_{n-1}$ versus $A_n$.
The (matrix-valued) coefficients are obtained as
\begin{equation}\label{An}
  A_n = (2n-1) \int\limits_0^\dt A(t) P_{n-1}\Big(\frac{t}{\dt}\Big) \, dt = (2n-1) \dt \int\limits_0^1 A(x\dt) P_{n-1}(x) dx \;.
\end{equation}
To see that $A_n$ is a term of order $\dt^n$, compare this expansion with an expansion $A(t) = \sum_{m \ge 1} a_m t^{m-1}$ in powers of $t$.
Since $P_n(x)$ is orthogonal to all $x^m$ with $m<n$,
we see from Eq.~\eqref{An} that
$A_n$ starts with the term
$\dt \int_0^1 a_n (x \dt)^{n-1} P_{n-1}(x) \, dx$ of order $\dt^n$.
In particular, it is a linear combination only of $a_m$ with $m \ge n$.

\subsection{Legendre expansion of $\Omega(t)$}\label{sec:LegOm}

If we insert the expansion Eq.~\eqref{Aexpansion} of $A(t)$ into the recursion Eq.~\eqref{OmRec} we obtain $\Omega(\dt)$ as a sum of nested commutators of the expansion coefficients $A_n$.
A nested commutator $[A_{n_1}, \dots, A_{n_m}]$ is of order $n_1 + \dots +n_m$ in $\dt$. The prefactor of this term is obtained as the $n$-fold integral $\xi(n_1,\dots,n_m) = \int_0^1 dx_1 \dots \int_0^{x_{m-1}} dx_m P_{n_1-1}(x_1)  \cdots P_{n_m-1}(x_m)$, which is a rational number independent of $A(t)$. For example,
\begin{equation}
\begin{split}
 \Omega_2(\dt) & = \frac{1}{2 \dt^2} \int_0^t \! dt_1 \int_0^{t_1} \! dt_2 \, [\sum_{n_1\ge 1} A_{n_1} P_{n_1-1}\Big(\frac{t_1}{\dt}\Big), \sum_{n_2 \ge 1} A_{n_2} P_{n_2-1}\Big(\frac{t_2}{\dt}\Big) ] \\
 & = \frac{1}{2} \sum_{n_1,n_2 \ge 1} \left( \int_0^1 \! dx_1 \int_0^{x_1} \! dx_2 \,  P_{n_1-1} (x_1) P_{n_2-1}(x_2) \right) \, [A_{n_1},A_{n_2}] \\  
& = \frac{1}{2} \left( \int_0^1 \! dx_1 \int_0^{x_1} \! dx_2 \,  P_0 (x_1) P_1(x_2) - P_1(x_1) P_0(x_2) \right) \, [A_1,A_2] = - \frac{1}{6} [A_1,A_2] \;.
\end{split}
\end{equation}
Notice that the only non-zero contributions in the second line come from $n_1=1, n_2=0$ and $n_1=0,n_2=1$, since the integral of $P_{n_1}(x_1) P_{n_2}(x_2)$ vanishes in all other cases.
This hints at a general pattern to be discussed below.

Collecting all terms up to order $\dt^9$, we find
\begin{equation}\label{OmTo8}
\begin{split}
\Omega(\dt) &=  A_1 - \frac{1}{6} [A_1,A_2]   \\
& +\frac{1}{60} [A_1,[A_1,A_3]] -\frac{1}{60} [A_2,[A_1,A_2]] 
 +\frac{1}{360} [A_1,[A_1,[A_1,A_2]]] -\frac{1}{30} [A_2,A_3] \\
 &-\frac{1}{70} [A_3,A_4] + \frac{1}{140} [A_2,[A_1,A_4]] -\frac{1}{210} [A_2,[A_2,A_3]] -\frac{1}{420} [A_3,[A_1,A_3]] -\frac{1}{210}[A_4,[A_1,A_2]]\\
 & -\frac{1}{840} [A_1,[A_1,[A_1,A_4]]] -\frac{1}{504} [[A_1,A_2],[A_1,A_3]] +\frac{1}{504} [A_2,[A_1,[A_1,A_3]]] -\frac{1}{840} [A_2,[A_2,[A_1,A_2]]] \\ & +\frac{1}{2520} [A_3,[A_1,[A_1,A_2]]]
 -\frac{1}{2520} [A_1,[A_1,[A_1,[A_1,A_3]]]]  -\frac{1}{7560} [[A_1,A_2],[A_1,[A_1,A_2]]] \\
 & +\frac{1}{2520} [A_2,[A_1,[A_1,[A_1,A_2]]]] -\frac{1}{15120} [A_1,[A_1,[A_1,[A_1,[A_1,A_2]]]]] + O(\dt^{10}) \;.
\end{split}
\end{equation}
The first line contains the 4th-order terms, the second line the 6th-order terms,
and the remaining lines the 8th-order terms.
This expression for $\Omega(\dt)$ avoids multi-dimensional integrals.

\subsection{Properties of the expansion}\label{sec:PropExp}

As seen above for $\Omega_2(\dt)$, only few out of the many possible commutators contribute to $\Omega(\dt)$.
In Eq.~\eqref{OmTo8} several nested commutators of order $\dt^9$ or less are missing, 
e.g. the terms $A_2, \dots, A_9$ or $[A_1,A_3], \dots, [A_1,A_8]$. This is a consequence of two general properties of the expansion that result in a zero prefactor $\xi(n_1,\dots,n_m)$ of $[A_{n_1},\dots,A_{n_m}]$.\\[1ex]
\emph{(P1)}  Time-reversal symmetry $U(\dt,0) = U(0,\dt)^{-1}$ of the propagator 
implies that $\Omega(t)$ changes sign if $A(t)$ is replaced with $-A(\dt-t)$.
According to the parity Eq.~\eqref{PSymm} of the Legendre polynomials it follows that even order terms in the expansion, i.e.  terms $[A_{n_1},\dots,A_{n_m}]$ with even $n_1 + \dots + n_m$, vanish.
This follows also from the calculation of the prefactor $\xi(n_1,\dots,n_m)$ as an $n$-fold integral: Each of the inner integrations over $x_2,\dots,x_n$ changes the parity of the integrand. The parity also changes by multiplication with a polynomial $P_n$ for odd $n$. Hence, the integrand in the final integration over $x_1$ has odd parity for odd $(n_1-1)+\dots+(n_m-1) + (-1)^{m-1}$, i.e. if $n_1+\dots+n_m$ is even.
Then, the integration gives zero and the respective term vanishes in Eq.~\eqref{OmTo8}.
\\[0.5ex]
\emph{(P2)} As a consequence of the orthogonality of the Legendre polynomials a term $[A_{n_1},\dots,A_{n_m}]$ vanishes if some index $n_k$ exceeds the sum of the others by two, i.e. $n_k > 1+ \sum_{i \ne k} n_i $ for a $k=1,\dots,m$. To see this change the integration order in the integral for the prefactor $\xi(n_1,\dots,n_m)$ such that the outermost integration is over $x_k$. This final integration is of the form $\int_0^1 d x_k P_{n_k-1}(x_k) p(x_k)$, where the polynomial $p(x_k)$ results from the previous $m-1$ integrations of the other polynomials $P_{n_i-1}(x_i)$.
The degree of $p(x_k)$ is at most $\sum_{i \ne k} n_i$. If, by assumption, this sum is smaller than $n_k-1$ the final integral is zero since $P_{n_k-1}$ is orthogonal to polynomials with smaller degree.
Now suppose $[A_{n_1},\dots,A_{n_m}]$ is a term of order $N$,
and one index $n_k > N/2$.
Then, $N \ge \sum_i n_i = n_k + \sum_{i \ne k} n_i > N/2 + \sum_{i \ne k} n_i$,
or $\sum_{i \ne k} n_i < N/2$.
The above condition applies, and it follows that this term gives no contribution.\\[1ex]
Both properties considerably simplify the derivation of approximate propagators since they reduce the number of terms in the expansion Eq.~\eqref{OmTo8}. 
According to (P1), approximately only half of the commutators contribute.
In particular, an expansion including all terms up to some odd order $\dt^N$ is automatically correct up to order $\dt^{N+1}$.
According to (P2) there are no contributions from terms such as $A_m$ for $m \ge 2$ or $[A_m,A_n]$ for $|m-n| \ne 1$.
It also explains why in the $N$th-order expansion of $\Omega(t)$ only commutators of terms $A_1,\dots,A_{N/2}$ occur (i.e. $A_1,\dots,A_4$ in Eq.~\eqref{OmTo8}).
This has the remarkable consequence that for the construction of $N$th-order propagators terms $A_n$ for $n \ge N/2+1$ can be neglected even if $n < N$.
Notice that property  (P1) is shared by an expansion in centered powers $(t-\dt/2)^n$, while (P2) requires orthogonality of the Legendre polynomials.
This fact motivated our use of a Legendre expansion of $A(t)$ 
instead of the apparently simpler Taylor expansion.

\subsection{Uniqueness of the expansion: Hall basis}\label{sec:Hall}

The expression for $\Omega(\dt)$ in Eq.~\eqref{OmTo8} is not unique.
Non-trivial identities between nested commutators, e.g. the Jacobi identity $[A,[B,C]]+[B,[C,A]]+[C,[A,B]]=0$, allow to replace one commutator by others.
To compare nested commutator expressions by equating the coefficients we must therefore first eliminate the ensuing linear dependencies.
Technically, this amounts to calculations using a vector space basis of the free Lie algebra generated by the $A_n$.
Since every nested commutator is a unique linear combination of the basis elements, uniqueness of the entire Magnus expansion is achieved. 

A systematic construction of free Lie algebra bases is provided by a Hall basis~\cite{MO99}. 
Algorithms exist for the rewriting of nested commutators in terms of the Hall basis elements,
and for their enumeration.
The number of Hall basis elements grows rapidly with the maximal order considered.
As listed in the following table,
\begin{center}
\begin{tabular}{lrrrrr}
 order N & 2 & 4 & 6 & 8  & 10\\ \hline
 full set of elements & 2  & 7 & 22 & 70 & 225\\
 relevant according to (P1), (P2) & 1 & 2 & 7 & 22 & 73
\end{tabular}
\end{center}
there are $70$ elements up to order $8$ in the Hall basis.
As a consequence of the two properties (P1), (P2) from Sec.~\ref{sec:PropExp}
only the $22$ elements in Table~\ref{tab:Hall} are relevant for our purposes.
Notice that in Eq.~\eqref{OmTo8} the elements $A_3$ and $[A_1,A_4]$ from the Hall basis are missing according to (P2), but yield order conditions for the CFETs as discussed in Sec.~\ref{sec:OrdCond}.
We do all calculations using the 22 Hall basis elements, rewriting commutators as necessary, e.g.
$[A_1,[A_2,[A_1,[A_1,A_2]]]] = [A_2,[A_1,[A_1,[A_1,A_2]]]] + [[A_1,A_2],[A_1,[A_1,A_2]]]$.

\begin{table}
\begin{equation*}
\begin{gathered}
A_1 \,, \; A_3 \,, \;
[A_1,A_2]  \,, \; 
[A_1,A_4]  \,, \; 
[A_2,A_3]  \,, \; 
[A_3,A_4] \\
[A_1,[A_1,A_3]]  \,, \; 
[A_2,[A_1,A_2]]  \,, \; 
[A_2,[A_1,A_4]]  \,, \; 
[A_2,[A_2,A_3]]  \,, \; 
[A_3,[A_1,A_3]]  \,, \; 
[A_4,[A_1,A_2]] \\
[A_1,[A_1,[A_1,A_2]]]  \,, \; 
[A_1,[A_1,[A_1,A_4]]]  \,, \; 
[A_2,[A_1,[A_1,A_3]]]   \\
[A_2,[A_2,[A_1,A_2]]]  \,, \; 
[A_3,[A_1,[A_1,A_2]]]  \,, \; 
[[A_1,A_2],[A_1,A_3]] \\
[A_1,[A_1,[A_1,[A_1,A_3]]]]  \,, \; 
[A_2,[A_1,[A_1,[A_1,A_2]]]]  \,, \; 
[[A_1,A_2],[A_1,[A_1,A_2]]] \\
[A_1,[A_1,[A_1,[A_1,[A_1,A_2]]]]]
\end{gathered}
\end{equation*}
\caption{The 22 odd elements up to order $8$ of the Hall basis with generators $A_1,\dots,A_4$.}
\label{tab:Hall}
\end{table}


\section{Commutator-free exponential time-propagators}
\label{sec:CF}

The expansion Eq.~\eqref{OmTo8} of $\Omega(\dt)$ still contains nested commutators.
An $N$th-order commutator-free exponential time-propagator (CFET) is based on the ansatz 
\begin{equation}\label{Ucf}
\tilde{U}^{(N)}_\mathrm{CF}(\dt) = e^{\Omega_1} e^{\Omega_2} \cdots e^{\Omega_s}  \;, 
\end{equation}
where each of the $s$ exponentials $\Omega_i$ is a linear combination
\begin{equation}\label{Omi}
  \Omega_i = \sum_{n=1}^N f_{i,n}  A_n
\end{equation}
of the $A_1,\dots,A_N$ from the Legendre expansion Eq.~\eqref{Aexpansion} of $A(t)$.
The CFET is completely determined through the coefficients $f_{i,n}$, which are fixed once and independently of the concrete $A(t)$ used in a calculation.
The practical evaluation of Eq.~\eqref{Ucf}, avoiding commutators and multi-dimensional integrals, is considerably simpler than for the original Magnus expansion. It will be discussed in more detail in Sec.~\ref{sec:Implement}.

Effectively, Eq.~\eqref{Ucf} is the exact propagator for an auxiliary problem with a fictitious, stepwise constant $\tilde{A}(t)$. The CFET coefficients $f_{i,n}$ must be determined in such a way that the replacement of the complicated time-dependent problem by the simpler auxiliary problem introduces only an error $\propto \dt^{N+1}$, independently of $A(t)$.
Now consider an $A(t) = x(t) X + y(t) Y$ which is the sum of two contributions $X,Y$. This situation arises, e.g., for a particle moving in a time-dependent field.
By construction, each $\Omega_i = x_i X + y_i Y$ itself is a sum of $X,Y$, with constant $x_i, y_i$ replacing $x(t), y(t)$.
Therefore, the CFET describes again a particle moving in a field, and thus preserves the principal physical situation. Notice, however, that fictitious negative time-steps can occur. The analogous statement does not hold for the original Magnus expansion involving commutators of $X,Y$.

The simplest example of a CFET is the 2nd-order midpoint rule
\begin{equation}\label{UCF2}
 \tilde{U}^{(2)}_\mathrm{CF2:1}(\dt) = \exp [ A_1 ]  = \exp \Big[ \int_0^\dt \! dt \, A(t) \Big ] 
  \simeq  \exp[ \dt \, A(\dt/2)] \;,
 \end{equation}	
corresponding to $s=1$ and $f_{1,1}=1$.
The second exponential is identical to the first  
according to the definition Eq.~\eqref{An} of $A_1$.
The last exponential is obtained by approximation of the integral through Gauss-Legendre quadrature (addressed later in Sec.~\ref{sec:GaussLeg}), which here reduces to evaluation of $A(t)$ at the midpoint $t=\dt/2$.

\subsection{Derivation of order conditions}\label{sec:OrdCond}

The construction of higher-order CFETs is substantially more difficult, and a systematic procedure is missing.
We adopt the following strategy:
Starting from the CFET ansatz Eq.~\eqref{Ucf}, the BCH formula~\eqref{BCH2} allows us to combine the $s$ exponentials until
we obtain $\tilde{U}^{(N)}_\mathrm{CF}(\dt)=e^{\tilde{\Omega}}$ with
\begin{equation}\label{tildeOm}
 \tilde{\Omega} = \sum_{i=1}^s f_{i,1} A_1 + \sum_{1 \le i < j \le s} \frac{f_{i,1} f_{j,2} -f_{j,1} f_{i,2}}{2} [A_1,A_2] + \dots 
\end{equation}
The $\tilde{\Omega}$ has to be compared with $\Omega(\dt)$ from the Magnus expansion Eq.~\eqref{OmTo8}, demanding equality of terms of order $\dt^{N}$ or less.
Working in a Hall basis, this implies equality of their prefactors which results
in equations for the coefficients $f_{i,n}$, the so-called order conditions.
Specifically, we find 
\begin{equation}\label{Order1}
  \sum_{i=1}^s f_{i,n} = \delta_{n,1} 
\end{equation}
arising from the terms $A_n$, and from $[A_1,A_2]$ 
\begin{equation}\label{Order2}
 \sum_{1 \le i < j \le s} f_{i,1} f_{j,2} -f_{j,1} f_{i,2} = -\frac{1}{3} \;. 
\end{equation}
For higher-order commutators,
the derivations become increasingly cumbersome, and calculations are best delegated to a computer. Since standard computer algebra systems are less useful for calculations in non-commutative algebras we used self-written programs that perform the Lie algebra manipulations, based on algorithms from Ref.~\cite{DeG00}.

Counting all $N$th-order elements in the Hall basis (Sec.~\ref{sec:Hall}), we see that the number of order conditions is $22$ ($70$) for order $6$ (order $8$), and thus appears to be too large for a practical solution of the multivariate polynomial equations that arise.
As we found in Sec.~\ref{sec:LegOm} several commutators do not appear in $\Omega(\dt)$ as a consequence of the two properties (P1), (P2). The key observation is that the corresponding order conditions can be satisfied by a suitably restricted choice of the $f_{i,n}$ according to the following two rules. \\[1ex]
\emph{(R1)} 
Since passing from $\tilde{U}_{CF}^{(N)}(\dt,0)$ to $\tilde{U}_{CF}^{(N)}(0,\dt)^{-1}$
changes the sign of $A_n$ by $(-1)^n$, a CFET complies with time-reversal symmetry if the coefficients obey
\begin{equation}\label{FSymm}
 f_{s-i+1,n} = (-1)^{n+1} f_{i,n} \;.
\end{equation}
For a time-symmetric CFET it thus suffices to specify the $f_{i,n}$ for $i \le s/2$, i.e. for the first half of the exponentials $e^{\Omega_i}$, and choose the remaining coefficients according to Eq.~\eqref{FSymm}.
For odd $s$, the coefficients $f_{(s+1)/2,n}$ of the
central exponential $i=(s+1)/2$ must be specified for odd $n$ only, while they are zero for even $n$. 
With this constraint the order conditions for even order terms, which do not contribute to $\Omega(\dt)$ according to (P1), are automatically satisfied. \\[0.5ex]
\emph{(R2)} Property (P2) states that up to order $\dt^{N}$ only terms $A_n$ with $n \le N/2$ contribute to $\Omega(\dt)$.
The order conditions involving higher-order $A_n$  can be satisfied simply by setting $f_{i,n}=0$ for $n > N/2$:
Since all coefficients are zero the corresponding commutators drop out entirely.
The remarkable implication is that an $N$th-order CFET can be built already from the terms $A_1,\dots,A_{N/2}$.
We note that this property is intrinsically connected with Gaussian quadrature using orthogonal polynomials (cf. Sec.~\ref{sec:GaussLeg}). It becomes obvious working with Legendre polynomials, while it requires sophisticated additional arguments in general~\cite{IN99}. \\[1ex]
By rule (R1) the number of relevant coefficients and order conditions is reduced approximately by one half. For this reason we consider only time-symmetric CFETs.
Notice that a symmetric $N$th-order CFET is automatically of order $N+1$, if $N$ is odd.
Rule (R2) implies that the summation index $n$ in Eq.~\eqref{Omi} only has to run from $1$ to $N/2$. We will later relax this rule to allow for minimization of the error, which requires inclusion of the term $A_{N/2+1}$.

With both rules, the number of order conditions is significantly reduced, 
to $2, 7, 22$ for $N=4, 6, 8$ CFETs (cf. the Table in Sec.~\ref{sec:Hall}).
On the other hand, a symmetric $N$th-order CFET with $s$ exponentials has 
$\lfloor{s N/4}\rfloor$ coefficients (rounding down to an integer).
The counting shows that 5 exponentials (11 exponentials) are needed for a 6th-order (8th-order) CFET.
Only  in exceptional cases solutions with less exponentials exist, e.g. CF6:4 in Table~\ref{6thCFET}.

\subsection{Fourth-order CFETs}\label{sec:4th}

\begin{table}
\centering
  \begin{tabular}{lll|lll}
  \hline
  \multicolumn{6}{l}{4th-order} \\\hline\hline
    \multicolumn{3}{l|}{2 exponentials} & \multicolumn{3}{l}{3 exponentials} \\\hline\hline
     \multicolumn{3}{l|}{CF4:2} &  \multicolumn{3}{l}{CF4:3} \\\hline
$f_{1,1}  = 1/2 $ & $f_{1,2} = 1/3$ & $f_{2,1}=0$  & $f_{1,1} = 11/40$ & $f_{1,2} = 20/87 $ & $f_{2,1} = 9/20$   \\\hline
 \end{tabular}
 \caption{Coefficients for 4th-order CFETs with $2$ and $3$ exponentials.
 Notice that CF4:3 is not recommended for use (cf. Sec.~\ref{sec:4therr}).
 }
   \label{4thCFET}
\end{table}

We consider 4th-order propagators ($N=4$) with three exponentials ($s=3$), of the form
\begin{equation}\label{UCF4:3}
  \tilde{U}^{(4)}_\mathrm{CF}(\dt) = \exp [f_{1,1} A_1 + f_{1,2} A_2 ] \exp [f_{2,1} A_1 ] \exp [f_{1,1} A_1 - f_{1,2} A_2 ] \;.
\end{equation}
As explained before (cf. Eqs.~\eqref{Order1},~\eqref{Order2}),
we get the two order conditions
 \begin{equation}
 \begin{split}
  1 & =  2 f_{1,1} + f_{2,1}  \;,  \\
  -\frac{1}{6}  & =  - (f_{1,1}+f_{2,1}) f_{1,2} \;.
\end{split}
 \end{equation}
The first arises from the term $A_1$, and the second from the term $[A_1,A_2]$.
In accordance with the above counting of terms, we have $3$ coefficients and $2$ order conditions.
Using $f_{2,1}$ as the free parameter, we find
\begin{equation}\label{Coeff4th}
  f_{1,1} = \frac{1-f_{2,1}}{2} \;, \quad f_{1,2} = \frac{1}{3(1+f_{2,1})} \;.
\end{equation}
Corresponding coefficients are listed in Table~\ref{4thCFET}.
The parameter $f_{2,1}$ will later allow for optimization of the propagator (see Sec.~\ref{sec:4thOpt}).
Setting $f_{2,1}=0$, we obtain the unique 4th-order CFET with $s=2$ exponentials
\begin{equation}\label{UCF4:2}
  \tilde{U}^{(4)}_\mathrm{CF4:2}(\dt) = \exp \Big[\frac{1}{2} A_1 + \frac{1}{3} A_2 \Big] \, \exp \Big[\frac{1}{2} A_1 - \frac{1}{3} A_2 \Big] \;.
\end{equation}
The notation used here and in the following is CF\textit{N}:\textit{s} for an $N$th-order CFET with $s$ exponentials.

\subsection{Sixth-order CFETs}\label{sec:6th}

For 6th-order ($N=6$), we consider propagators with $s=6$ exponentials.
The $9$ coefficients $f_{i,n}$, for $1 \le i,n \le 3$,
must satisfy $7$ order conditions corresponding to the 7 Hall basis elements 
\begin{equation}
\begin{gathered}
A_1 \,, \; A_3 \,, \;
[A_1,A_2]  \,, \; 
[A_2,A_3]  \,, \; 
[A_1,[A_1,A_3]]  \,, \; 
[A_2,[A_1,A_2]]  \,, \; 
[A_1,[A_1,[A_1,A_2]]] 
\end{gathered}
\end{equation}
from Table~\ref{tab:Hall}.
We note that two coefficients can be chosen as a free parameter.

An explicit solution of the order conditions is possible to a large degree,
and simple explicit expressions for the coefficients can be obtained in some cases (cf.~\ref{app:6th}).
Setting $f_{3,2}=0$, the two central exponentials can be combined, resulting in propagators with $s=5$ exponentials and a single free parameter.
Surprisingly, there is also a solution with $f_{3,1}=f_{3,2}=f_{3,3}=0$, giving a 6th-order CFET with only 4 exponentials (see CF6:4 in Table~\ref{6thCFET}), although there are less coefficients than order conditions.
We do not know whether the existence of this solution is accidental, or hints at a general redundancy pattern of the equations.
For practical purposes, the CFET CF6:5 from Table~\ref{6thCFET} is most relevant,
since it has small approximation error. 
Further optimized 6th-order CFETs will be obtained in Sec.~\ref{sec:6thOpt}.

\begin{table}
  \begin{tabular}{lll}
  \hline
  \multicolumn{3}{l}{6th-order} \\\hline\hline
    \multicolumn{3}{l}{4 exponentials} \\\hline\hline
     \multicolumn{3}{l}{CF6:4} \\\hline\noalign{\smallskip}
     \multicolumn{3}{l}{
$f_{1,1}  = \dfrac{1}{2} + \dfrac{(5400 - 600 \sqrt{6})^{1/3}}{60} + \dfrac{\Big(\frac{1}{5} (9 + \sqrt{6})\Big)^{1/3}}{ 2 \cdot 3^{2/3}}$ } \\[0.3ex]
 & $ f_{1,2} = f_{1,1} - \frac{2}{3} f_{1,1}^2 $ & $f_{1,3} = \dfrac{1}{10-10 f_{1,1}}$ \\[0.75ex]
$f_{2,1} = \frac{1}{2} - f_{1,1}$ & $f_{2,2} = \frac{1}{3}(1-4 f_{1,1} + 2 f_{1,1}^2)$ & $f_{2,3} = - f_{1,3}$ \\[0.75ex]
$f_{3,1}=0$ & $f_{3,2}=0$ & $f_{3,3}=0$ \\\hline\hline
    \multicolumn{3}{l}{5 exponentials} \\\hline\hline
  \multicolumn{3}{l}{CF6:5} \\\hline\noalign{\smallskip}
    $f_{1,1}=\phantom{-} 0.16 $ & $f_{1,2}=\phantom{-}0.14587456942714338561$ & $f_{1,3}=\phantom{-}0.11762370828143015682$ \\
    $f_{2,1}= \phantom{-}0.38752405202531186588 $ & $f_{2,2}=\phantom{-}0.15089113704380764664$ & $f_{2,3}=-0.12805075909013044594$ \\
    $f_{3,1}=\phantom{-}1-2f_{2,1}-2f_{1,1} $ & $f_{3,2}=\phantom{-}0$ & $f_{3,3}=-2f_{2,3}-2f_{1,3}$ \\\hline
     \multicolumn{3}{l}{CF6:5b (cf. Ref.~\cite{BM06})} \\\hline\noalign{\smallskip}
    $f_{1,1}=\phantom{-} 0.2 $ & $f_{1,2}=\phantom{-} 0.1746879190177786220$ & $f_{1,3}=\phantom{-}0.12406375705333586606$ \\
    $f_{2,1}= \phantom{-}0.34815492558797391479 $ & $f_{2,2}=\phantom{-}0.1068765450953683$ & $f_{2,3}=-0.139021313323765096675$ \\
    $f_{3,1}=\phantom{-}1-2f_{2,1}-2f_{1,1} $ & $f_{32}=\phantom{-}0$ & $f_{33}=-2f_{2,3}-2f_{1,3}$ \\\hline\hline
    \multicolumn{3}{l}{6 exponentials} \\\hline\hline
     \multicolumn{3}{l}{CF6:6} \\\hline\noalign{\smallskip}    
$f_{1,1}  = \phantom{-}0.16$ & $ f_{1,2} =\phantom{-}0.15101538937746543493$ & $f_{1,3} =\phantom{-}0.13304616813239630479$ \\[0.75ex]
$f_{2,1} = -0.22738164742696330169$ & $f_{2,2} =-0.087654259755115431662$ & $f_{2,3} =\phantom{-}0.069919836812656575583$ \\[0.75ex]
$f_{3,1}=\phantom{-}1/2-f_{1,1}-f_{2,1}$ & $f_{3,2}=\phantom{-}0.21035154512209824847$ & $f_{3,3}=-f_{1,3}-f_{2,3}$ \\\hline
 \end{tabular}
 \caption{Coefficients for unoptimized 6th-order CFETs with $s=4,5,6$ exponentials (coefficients for optimized 6th-order CFETs are given in Table~\ref{6thCFETOpt}).
 The CFET CF6:5b corresponds to the coefficients of the propagator $\psi_5^{[6]}$ from Ref.~\cite{BM06}.}
   \label{6thCFET}
\end{table}

\subsection{Eighth-order CFETs}\label{sec:8th}

The $22$ order conditions of 8th-order CFETs correspond to the entire set of commutators from Table~\ref{tab:Hall}. Exactly $22$ coefficients exist for $s=11$ exponentials.
Due to their complexity, the order conditions can only be solved numerically.
Several solutions were computed using a root finder based on the Newton iteration~\cite{PFTV86}. 
Severe ill-conditioning of the equations required the use of high-precision arithmetics, based on the MPFUN package~\cite{MPF90}, and repeated restarting of the Newton iteration.
The coefficients of an 8th-order CFET with small approximation error, selected from about $50$ computed solutions of the order conditions, are given in Table~\ref{8thCoeff}.
A systematic search of the coefficient space was not possible.

\begin{table}
\begin{center}
  \begin{tabular}{ll}
\hline
    \multicolumn{2}{l}{8th-order: 11 exponentials} \\\hline\hline
     \multicolumn{2}{l}{CF8:11} \\\hline\noalign{\smallskip}     
$f_{1,1}=\phantom{-}0.169715531043933180094151$ &
$f_{1,2}=\phantom{-} 0.152866146944615909929839$ \\
$f_{1,3}=\phantom{-} 0.119167378745981369601216$ &
$f_{1,4}=\phantom{-} 0.068619226448029559107538$ \\[0.5ex]
$f_{2,1}=\phantom{-} 0.379420807516005431504230$ &
$f_{2,2}=\phantom{-} 0.148839980923180990943008$ \\
$f_{2,3}=-0.115880829186628075021088$ &
$f_{2,4}=-0.188555246668412628269760$ \\[0.5ex]
$f_{3,1}=\phantom{-} 0.469459306644050573017994$ &
$f_{3,2}=-0.379844237839363505173921$ \\
$f_{3,3}=\phantom{-} 0.022898814729462898505141$ &
$f_{3,4}=\phantom{-} 0.571855043580130805495594$ \\[0.5ex]
$f_{4,1}=-0.448225927391070886302766$ &
$f_{4,2}=\phantom{-} 0.362889857410989942809900$ \\
$f_{4,3}=-0.022565582830528472333301$ &
$f_{4,4}=-0.544507517141613383517695$ \\[0.5ex]
$f_{5,1}=-0.293924473106317605373923$ &
$f_{5,2}=-0.026255628265819381983204$ \\
$f_{5,3}=\phantom{-} 0.096761509131620390100068$ &
$f_{5,4}=\phantom{-} 0.000018330145571671744069$ \\[0.5ex]
$f_{6,1}=\phantom{-} 0.447109510586798614120629$ &
$f_{6,3}=-0.200762581179816221704073$ \\\hline
\end{tabular}
\end{center}
\caption{Coefficients for an 8th-order CFET with $11$ exponentials.
}
\label{8thCoeff}
\end{table}

\section{Theoretical error analysis}\label{sec:ErrorAnalysis}

The CFET error is determined by the difference $\chi = \tilde{\Omega}-\Omega$ between the exact $\Omega(\dt)$ from the Magnus expansion Eq.~\eqref{OmTo8}  and the approximate $\tilde{\Omega}$ from Eq.~\eqref{tildeOm}.
The theoretical error analysis aims at minimization of the error term in the general situation,
where no specific information about $A(t)$ is available.

\subsection{General considerations}

By construction the error term is of the form $\chi = \sum_k (p_k - c_k) C_k$, where the $C_k$ are the $N+1$-order commutators from the Hall basis, the $p_k$ are polynomials in the coefficients $f_{i,n}$ such as in Eq.~\eqref{tildeOm}, and the $c_k$ the constant prefactors from Eq.~\eqref{OmTo8}.
The size of $\chi$ can be measured with a matrix norm $\| \cdot \|$. It is
\begin{equation}\label{Delta}
 \|\chi \| = \| \sum_k (p_k-c_k) C_k \| \le \sum_k |p_k-c_k|  \, \cdot \, \| C_k \|  \;.
 \end{equation}

In concrete situations, $\| \chi \|$ depends not only on the size $\| C_k \|$ but also on the amount of dependency between different $C_k$, which is responsible for the difference between left hand and right hand side of the above inequality.
 In the general case we may not assume that the difference is small. Accidental cancellations, i.e. $\| (p_k-c_k) C_k + (p_l-c_l) C_l\| \approx 0$ for a $k \ne l$, are typical.
Optimization of a CFET, that is minimization of $\|\chi\|$ through variation of the coefficients $f_{i,n}$, thus requires that all $|p_k-c_k|$ become simultaneously small.
Only then, the error can be expected to be small in the general case.
Such optimized CFETs are universally applicable and perform equally well in different situations.
Optimization will be achieved for 4th- and 6th-order CFETs, listed below in Tables~\ref{4thCFETOpt},~\ref{6thCFETOpt}. 
For 8th-order CFETs, optimization is not practicable due to the complexity of the order conditions.

An important point, which seems to have been missed in the literature, complicates the error analysis in comparison to split-operator techniques.
While rule (R2) in Sec.~\ref{sec:OrdCond} states that the terms $A_n$ for $n>N/2$ can be
disregarded in the construction of an $N$th-order CFET,
the error term $\chi$ contains a contribution from $A_{N/2+1}$ since 
the prefactors $c_k$ of the corresponding commutators are non-zero in Eq.~\eqref{OmTo8}.
For 4th-order, this applies to the terms $[A_1,[A_1,A_3]]$ and $[A_2,A_3]$ involving $A_3$. 
In contrast to the basic construction of higher-order CFETs with (R2), CFET optimization requires explicit inclusion of $A_{N/2+1}$.
Therefore, the optimized 4th- and 6th-order CFETs include non-zero coefficients for the $A_3$ or $A_4$ term, respectively.
Additional order conditions, e.g. $\sum_{i=1}^s f_{i,N/2+1} =0$ arising from the $A_{N/2+1}$ term itself, must be accounted for.
We note that for split-operator techniques~\cite{MQ02},  where essentially the full Magnus propagator $e^{\Omega(t)}$ is replaced by the term $e^{(X+Y)t}$ for a time-independent $A(t) \equiv X+Y$, the equivalent coefficients $c_k=0$, and no additional provisions are necessary.

Since the error term $\chi$ is of order $\dt^{N+1}$ we must ask whether also terms $A_n$ for $N/2+1<n \le N+1$ need be considered. Property (P1) states that nested commutators involving these terms do not occur in Eq.~\eqref{OmTo8} up to order $N+1$, i.e. the corresponding prefactor $c_k=0$ in Eq.~\eqref{Delta}.
Since we have set the coefficients of these $A_n$ to zero by rule (R2), they do not contribute to $\chi$ and need not be considered.
Notice again that the use of Legendre polynomials in~\eqref{Aexpansion} simplifies the derivation:
With a power series expansion all terms up to order $\dt^{N+1}$ would explicitly contribute to $\chi$, and minimization of $|\chi|$ would result in a number of additional though redundant equations.

\subsection{Optimized fourth-order CFETs}\label{sec:4thOpt}

For a 4th-order CFET including the $A_3$ term, 
we make the ansatz 
\begin{equation}\label{4thOptAnsatz}
  \tilde{U}^{(4)}_\mathrm{CF4}(\dt) = \exp[ f_{1,1} A_1 + f_{1,2} A_2 +f_{1,3} A_3] \exp [f_{2,1} A_1 +f_{2,3} A_3] \exp [f_{1,1} A_1 - f_{1,2} A_2 + f_{1,3} A_3] \;,
\end{equation}	
in extension of Eq.~\eqref{UCF4:3}.
The previous order conditions still apply, and $f_{1,1}$, $f_{1,2}$ are given by Eq.~\eqref{Coeff4th}.
The new order condition arising from the $A_3$ term is $2 f_{1,3} + f_{2,3} =0$,
which gives one additional free parameter $f_{2,3}$ with $f_{1,3}=-f_{2,3}/2$.
With these choices, we obtain for the error term
\begin{equation}\label{Delta4thOpt}
\begin{split}
 \chi =  \tilde{\Omega}-\Omega =&  \Big( \frac{1}{60} - \frac{1+2 f_{2,1}}{54(1+f_{2,1})^2} \Big)[A_2,[A_1,A_2]] + \Big( \frac{1}{1440} - \frac{f_{2,1}^2}{288} \Big)  [A_1,[A_1,[A_1,A_2]]] \\
&-\Big( \frac{1}{24}(1+f_{2,1})f_{2,3} +\frac{1}{60} \Big) [A_1,[A_1,A_3]] + \Big(\frac{f_{2,3}}{6(1+f_{2,1})} +  \frac{1}{30} \Big) [A_2,A_3] + O(\dt^7) \;.
\end{split}
\end{equation}
It has four contributions corresponding to the 5th-order exponentials in the second line of Eq.~\eqref{OmTo8}.

\begin{figure}
  \begin{center}
    \includegraphics[width=0.45\textwidth]{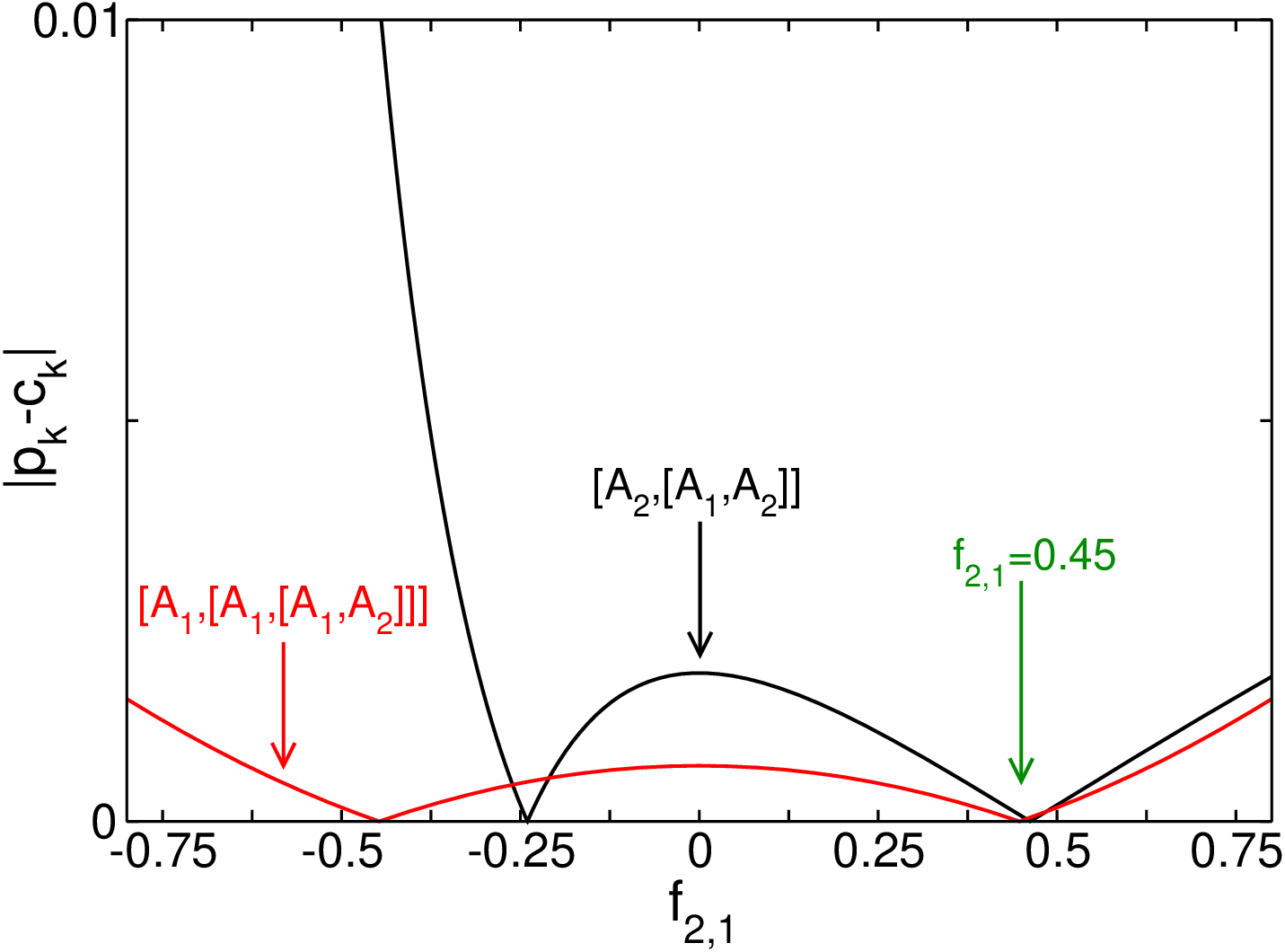}
    \hfill
    \includegraphics[width=0.45\textwidth]{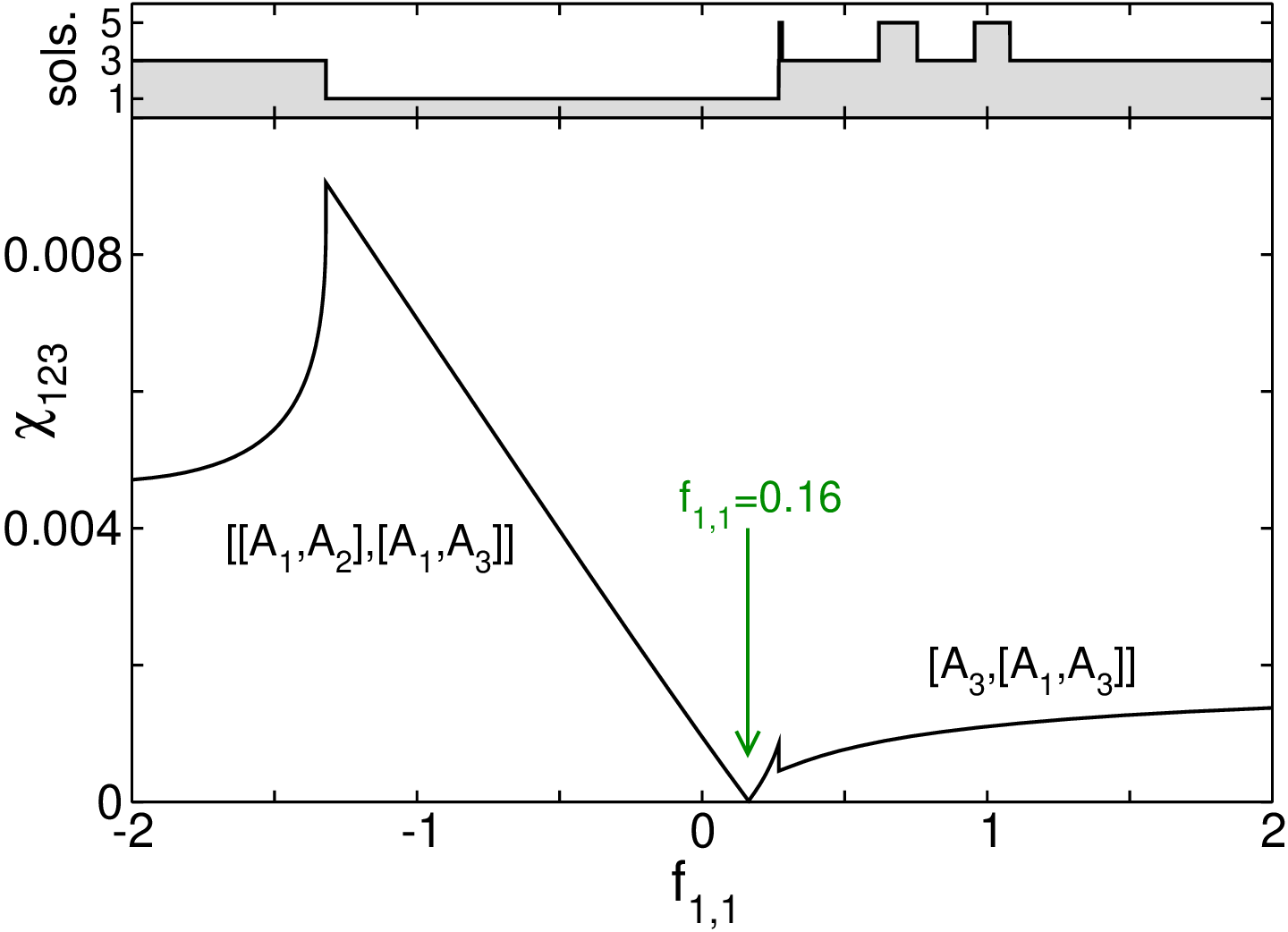}
   \end{center}
\caption{
Left panel: Contributions $|p_k-c_k|$ to the 5th-order error term $\chi$ of 4th-order CFETs (Eq.~\eqref{Delta4thOpt}) as a function of the free parameter $f_{2,1}$,
for commutators $[A_2,[A_1,A_2]]$ (black) and $[A_1,[A_1,[A_1,A_2]]]$ (red).
Both contributions become small in the vicinity of $f_{2,1}=0.45$.
Right panel: Maximum $\chi_{123}$ of the contributions $|p_k-c_k|$ to the 7th-order error term $\chi$ for 6th-order CFETs, excluding contributions from the $A_4$ term.
As indicated in the figure, the maximal contribution comes from either the $[[A_1,A_2],[A_1,A_3]]$
(for $f_{1,1} \lesssim 0.16$) or the $[A_3,[A_1,A_3]]$ term ($f_{1,1} \gtrsim 0.16)$
The upper panel indicates the number of solutions of the $s=5$, $N=6$ order conditions.
}
\label{fig:Chi}
\end{figure}

In Fig.~\ref{fig:Chi} (left panel) we show exemplarily $|p_k-c_k|$ for $k=[A_2,[A_1,A_2]]$ and $k=[A_1,[A_1,[A_1,A_2]]]$ as a function of $f_{2,1}$. The optimal choice is close to $f_{2,1} \approx 0.45=9/20$.
If we only try to minimize the contribution from these two commutators, neglecting the $A_3$ terms, we thus obtain the CFET CF4:3 from Table~\ref{4thCFET}. We will see below in Sec.~\ref{sec:4therr} that this CFET is far from being optimal.
Full minimization of $\chi$ through variation of $f_{2,1}$ and $f_{2,3}$, including the $A_3$ terms, results in $f_{2,1} \approx 0.45$, $f_{2,3} \approx -0.28$.
For an optimized 4th-order CFET we thus propose the choice $f_{2,1} = 9/20$, $f_{2,3} = -7/25$,
 which results in (cf. Table~\ref{4thCFETOpt})
\begin{equation}\label{CF4Optimal}
 U_\mathrm{CF4:3Opt}(\dt) = \exp \Big[\frac{11}{40} A_1 + \frac{20}{87} A_2 + \frac{7}{50} A_3 \Big]  \,
         \exp\Big[\frac{9}{20} A_1 - \frac{7}{25} A_3\Big] \,
         \exp\Big[\frac{11}{40} A_1 - \frac{20}{87} A_2 + \frac{7}{50} A_3\Big] \;. 
\end{equation}

\begin{table}
\centering
  \begin{tabular}{lllll}
  \hline
  \multicolumn{5}{l}{4th-order (optimized)} \\\hline\hline
    \multicolumn{5}{l}{3 exponentials} \\\hline\hline
  \multicolumn{5}{l}{CF4:3Opt} \\\hline\noalign{\smallskip}
  $f_{1,1} = 11/40$ & $f_{1,2} =20/87$ & $f_{1,3} = 7/50$ & $f_{2,1}=9/20$ & $f_{2,3}=-7/25$  \\\hline
 \end{tabular}
 \caption{Coefficients for the optimized 4th-order CFET CF4:3Opt with $3$ exponentials (Eq.~\eqref{CF4Optimal}).
The unoptimized CFET CF4:3 from Table~\ref{4thCFET} is obtained by dropping the $A_3$ terms from CF4:3Opt.
 }
   \label{4thCFETOpt}
\end{table}

\subsection{Optimized sixth-order CFETs}\label{sec:6thOpt}

Extending 6th-order CFETs with $6$ exponentials (Sec.~\ref{sec:6th}) by inclusion of the term $A_4$
provides us with three additional coefficients $f_{1,4}$, $f_{2,4},f_{3,4}$.
The new order condition 
\begin{equation}
 0 = f_{1,1} f_{1,4} + 2 f_{2,1}f_{1,4} + f_{2,1} f_{2,4} + 2 f_{3,1} f_{1,4} +
  2 f_{3,1} f_{2,4} + f_{3,1} f_{3,4}
\end{equation}
arising from the commutator $[A_1,A_4]$ fixes the value of $f_{3,4}$.
Notice that the term $A_4$ itself does not lead to a new order condition, since it is even and the associated $c_k = p_k=0$ by rule (R1).
Using the explicit solution of the order conditions (cf.~\ref{app:6th}),
we can minimize the error term $\chi$ through variation of the four free parameters $f_{1,1}$, $f_{3,2}$, $f_{1,4}$, $f_{2,4}$. 
For $s=5$ exponentials, we set $f_{3,2}=f_{3,4}=0$.
To illustrate the typical behaviour, we show in Fig.~\ref{fig:Chi} (right panel) the partial error $\chi_{123}=\max {|p_k-c_k|}$  including only the contributions from commutators $C_k$ without the $A_4$ term. It depends on the single parameter $f_{1,1}$.
Optimal choices occur around $f_{1,1} \approx 0.16$, corresponding to the CFET CF6:5 from Table~\ref{6thCFET}.
This also provides partial justification for the CFET CF6:5b with $f_{1,1}=0.2$ from Ref.~\cite{BM06}.
Inclusion of the $A_4$ term and subsequent minimization of the associated error contribution, keeping $f_{1,1}$ fixed, results in the improved CFET CF6:5Imp.
The full minimization of $|\chi|$ with free variation of all parameters results in the optimized 6th-order CFETs CF6:5Opt and CF6:6Opt listed in Table~\ref{6thCFETOpt}.

\begin{table}
  \begin{tabular}{lll}
  \hline
  \multicolumn{3}{l}{6th-order (optimized)} \\\hline\hline
    \multicolumn{3}{l}{5 exponentials} \\\hline\hline
  \multicolumn{3}{l}{CF6:5Imp} \\\hline\noalign{\smallskip}
      $f_{1,1}=\phantom{-} 0.16 $ & $f_{1,2}=\phantom{-}0.14587456942714338561$ & $f_{1,3}=\phantom{-}0.11762370828143015682$ \\
    $f_{2,1}=\phantom{-} 0.38752405202531186588 $ & $f_{2,2}=\phantom{-}0.15089113704380764664$ & $f_{2,3}=-0.12805075909013044594$ \\
    $f_{3,1}=\phantom{-}1-2f_{2,1}-2f_{1,1} $ & $f_{3,2}=\phantom{-}0$ & $f_{3,3}=-2f_{2,3}-2f_{1,3}$ \\
    $f_{1,4}=\phantom{-}0.074$ & $f_{2,4}=-0.212530296697694739551$ & $f_{3,4}=\phantom{-}0$ \\\hline
   \multicolumn{3}{l}{CF6:5Opt} \\\hline\noalign{\smallskip}
    $f_{1,1}=\phantom{-} 0.1714 $ & $f_{1,2}=\phantom{-}0.15409059414309687213$ & $f_{1,3}=\phantom{-}0.11947178242929061641$ \\
    $f_{2,1}=\phantom{-}0.37496374319946236513$ & $f_{2,2}=\phantom{-}0.13813675394387646682$ & $f_{2,3}=-0.13090674649282935743$ \\
    $f_{3,1}=\phantom{-}1-2f_{2,1}-2f_{1,1} $ & $f_{3,2}=\phantom{-}0$ & $f_{3,3}=-2f_{2,3}-2f_{1,3}$ \\
    $f_{1,4}=\phantom{-}0.07195$ & $f_{2,4} =-0.21123356253315514306$ & $f_{3,4} =\phantom{-}0$  \\\hline\hline
    \multicolumn{3}{l}{6 exponentials} \\\hline\hline
     \multicolumn{3}{l}{CF6:6Opt} \\\hline\noalign{\smallskip}    
$f_{1,1}  = \phantom{-}0.3952$ & $ f_{1,2} =\phantom{-}0.35629343479227292880$ & $f_{1,3} =\phantom{-}0.27848030437681878641$ \\[0.75ex]
$f_{2,1} = -0.22432144875476807927$ & $f_{2,2} =-0.19935407393749030416$ & $f_{2,3} =-0.15625650102884866893$ \\[0.75ex]
$f_{3,1}=\phantom{-}1/2-f_{1,1}-f_{2,1}$ & $f_{3,2}=\phantom{-}0.1145$ & $f_{3,3}=-f_{1,3}-f_{2,3}$ \\
$f_{1,4}=\phantom{-}0.1579$ & $f_{2,4} =-0.09512$ & $f_{3,4} = -0.16475168057141371958$ \\\hline
 \end{tabular}
 \caption{Coefficients for optimized 6th-order CFETs with $s=5,6$ exponentials.
 In each case, the last row gives the coefficients for the $A_4$ term.
 The CFET CF6:5Imp is obtained from CF6:5 (Table~\ref{6thCFET}) through separate minimization of the $A_4$ error contributions.}
   \label{6thCFETOpt}
\end{table}


\section{Practical error analysis}
\label{sec:Error}

The theoretical error analysis results in optimized CFETs, whose error is expected to be small in the general case.
In a concrete situation dependencies between the nested commutators in the error term may lead to different results. 
To confirm the validity of the theoretical error analysis we study the CFET error for a driven two-level system.
Further issues of practical relevance concern the choice between CFETs of different order,
and the time-step selection.

\subsection{Time-stepping and effective error}\label{sec:EffErr}

In the standard time-stepping approach, the approximate propagator $\tilde{U}(t)$ over longer propagation times is constructed as a product of short-time CFETs $\tilde{U}_\mathrm{CF}^{(N)}(t+\dt,t)$. Equivalently, a concrete solution $x(t)$ is repeatedly propagated over a small time step $\dt$.
The propagator $\tilde{U}(T)$ for the maximal propagation time $T$ is a product of $N_s = T/\dt$ CFETs. Intermediate results are obtained 
at multiples of $\dt$. 

The accuracy of time-stepping is controlled through the size of $\dt$. 
For $N$th-order CFETs $\tilde{U}_\mathrm{CF}^{(N)}(t+\dt,t)$,
the error contributed by each one scales as $\dt^{N+1}$.
Due to accumulation of errors, the propagation error after $N_s$ steps is $\epsilon = c N_s \dt^{N+1} = c T \dt^N$ with an error constant $c$ which depends on the concrete situation.
To achieve a given accuracy requires a time step $\dt \le(\epsilon/cT)^{1/N}$ for a maximal acceptable error $\epsilon$. Usually, $\dt \ll T$.
The computational effort is dominated by the $N_s$-fold evaluation of the $s$ exponentials in Eq.~\eqref{Ucf}.
It is thus proportional to 
$s N_s = sT/\dt \ge \bar{c} T^{1+1/N} \epsilon^{-1/N}$
with the effective error constant
\begin{equation}\label{BarC}
 \bar{c} = s c^{1/N} \;.
\end{equation}
This quantity determines the efficiency of time-propagation with an $N$th-order CFET with $s$ stages. As a rule of thumb we note the relation
\begin{equation}\label{thumb}
\mathsf{effort} \propto \mathsf{error}^{-1/N} \times \mathsf{time}  \;.
\end{equation}

\subsection{Driven two-level system}\label{sec:2level}

Our test problem is a driven two-level system, realized, e.g., by a spin $1/2$ in a magnetic field $\vec{B}(t)=(B_x(t),B_y(t),B_z(t))$.
In the eigenbasis of the z-component of angular momentum, the Hamilton operator is given by the matrix 
\begin{equation}
  H(t) = \frac{1}{2} \sum_{k=x,y,z} B_k(t) \sigma_k = \frac{1}{2} \begin{pmatrix} B_z(t) & B_x(t) - \ii B_y(t) \\ B_x(t) + \ii B_y(t) & -B_z(t)  \end{pmatrix} \;,
\end{equation}
with the standard Pauli matrices~\cite{Mes61}
\begin{equation}
  \sigma_x  = \begin{pmatrix}  0 & 1 \\ 1 & 0\end{pmatrix} \;, \;
  \sigma_y  = \begin{pmatrix}  0 & -\ii \\ \ii & 0 \end{pmatrix} \;, \;
  \sigma_z  = \begin{pmatrix}  1 & 0 \\ 0 & -1 \end{pmatrix} \;.
\end{equation}

For particular choices of $\vec{B}(t)$ the propagator can be expressed in simple, closed form.
One example is the periodically driven two-level system with $\vec{B}(t) =  (2V \cos 2\omega t, 2V \sin 2\omega t,  2\Delta)$, or
\begin{equation}\label{TwoLevelH}
 H(t) = \begin{pmatrix} \Delta & V e^{-2 \ii \omega t}  \\ V e^{2 \ii \omega t}& -\Delta \end{pmatrix} \;,
\end{equation}
where $\Delta$, $V$, $\omega \in \mathbb{R}$.
The exact propagator is given by
\begin{equation}\label{TwoLevelU}
  U(t,0)=
  \begin{pmatrix}
    e^{- \ii\omega t}(\cos \Omega t - \ii\dfrac{\Delta-\omega}{\Omega} \sin \Omega t) & - \ii \dfrac{V}{\Omega} e^{-\ii\omega t} \sin \Omega t \\
 - \ii \dfrac{V}{\Omega} e^{\ii\omega t} \sin \Omega t & e^{\ii\omega t}(\cos \Omega t + \ii\dfrac{\Delta-\omega}{\Omega} \sin \Omega t)
  \end{pmatrix} 
\end{equation}
with $\Omega=\sqrt{(\Delta-\omega)^2+V^2}$.
We note that, in accordance with Floquet theory for periodically driven systems,
$U(\pi n /\omega,0) = U(\pi/\omega,0)^n$ for integer $n$.
The transition probability $\text{spin up} \leftrightarrow \text{spin down}$
\begin{equation}
  P(t) = |U_{21}(t,0)|^2 = \left(\frac{V}{\Omega}\right)^2 \; \sin^2 \Omega t
\end{equation}
is typical for a Breit-Wigner resonance.

In the case of two-level systems, application of a CFET requires evaluation of matrix exponentials $e^{\Omega_i}$, which correspond to propagation with fictitious constant magnetic fields.
Each exponential can be evaluated in closed form with the relation 
\begin{equation}
 \exp( \ii \frac{\phi}{2}  \vec{n} \cdot \vec{\sigma}    )  = \cos \frac{\phi}{2}  + \ii \sin\frac{\phi}{2} \vec{n} \cdot \vec{\sigma}   \qquad ({|\vec{n}}| = 1)
 \end{equation}
for the spin $1/2$ rotation operator.

To quantify the CFET error we calculate the deviation $\epsilon(t)=\|\tilde{U}(t)-U(t)\|$ of the approximate propagator $\tilde{U}(t)$ from the exact $U(t)$.
We use the Frobenius norm for a $L \times L$ square matrix
\begin{equation}\label{Frob}
  \|U-\tilde{U}\|^2 =\frac{1}{L} \mathrm{tr}[(U-\tilde{U})^\dagger (U-\tilde{U})] =\frac{1}{L} \sum_{ij} |U_{ij}-\tilde{U}_{ij}|^2 \;,
\end{equation}
where $\mathrm{tr}[ \cdot ]$ denotes the trace. This choice is particularly convenient for the Schr{\"o}dinger equation, where the propagators are unitary such that $\|U\|=1$ and 
$\|U-\tilde{U}\|^2 = 2 \|1-U^\dagger \tilde{U} \|$. 
Notice that this definition accounts for phase slips of the propagators.
With the BCH formula we find  
that the CFET error
\begin{equation}
\epsilon(\dt) = \|U(\dt)-\tilde{U}(\dt)\| = 2 \|1-e^{-\Omega(\dt)} e^{\tilde{\Omega}(\dt)} \| = 2 \|1 - e^{\chi + O(\dt^{N+2})} \| = 2 \| \chi\| + O(\dt^{N+2}) 
\end{equation}
is indeed determined by the error term $\chi = \tilde{\Omega}(\dt)-\Omega(\dt)$.
The commutator relations of the spin algebra
imply that the nested commutators $C_k$ in Eq.~\eqref{Delta} are not independent.
This allows us to check the theoretical error analysis from Sec.~\ref{sec:ErrorAnalysis} in a situation where cancellation between different $C_k$ plays a role.
We note that cancellation is not a peculiar consequence of the small Hilbert space of the present example, but of commutator relations dictated by physics. Similar relations hold in any relevant situation.

\subsection{Fourth-order CFETs}\label{sec:4therr}

To determine the effective error constant $\bar{c}$, 
we propagate the driven two-level system Eq.~\eqref{TwoLevelH} over $20$ periods of the driving field, i.e. up to a time $T=20\pi/\omega$.
From the maximal propagation error $\epsilon=\max\{\epsilon(t)|0 \le t \le T\}$ we get the effective error constant as $\bar{c}=(s/\dt) (\epsilon/T)^{1/N}$ in the limit $\dt \ll T$.

\begin{figure}
\begin{center}
\includegraphics[width=0.45\textwidth]{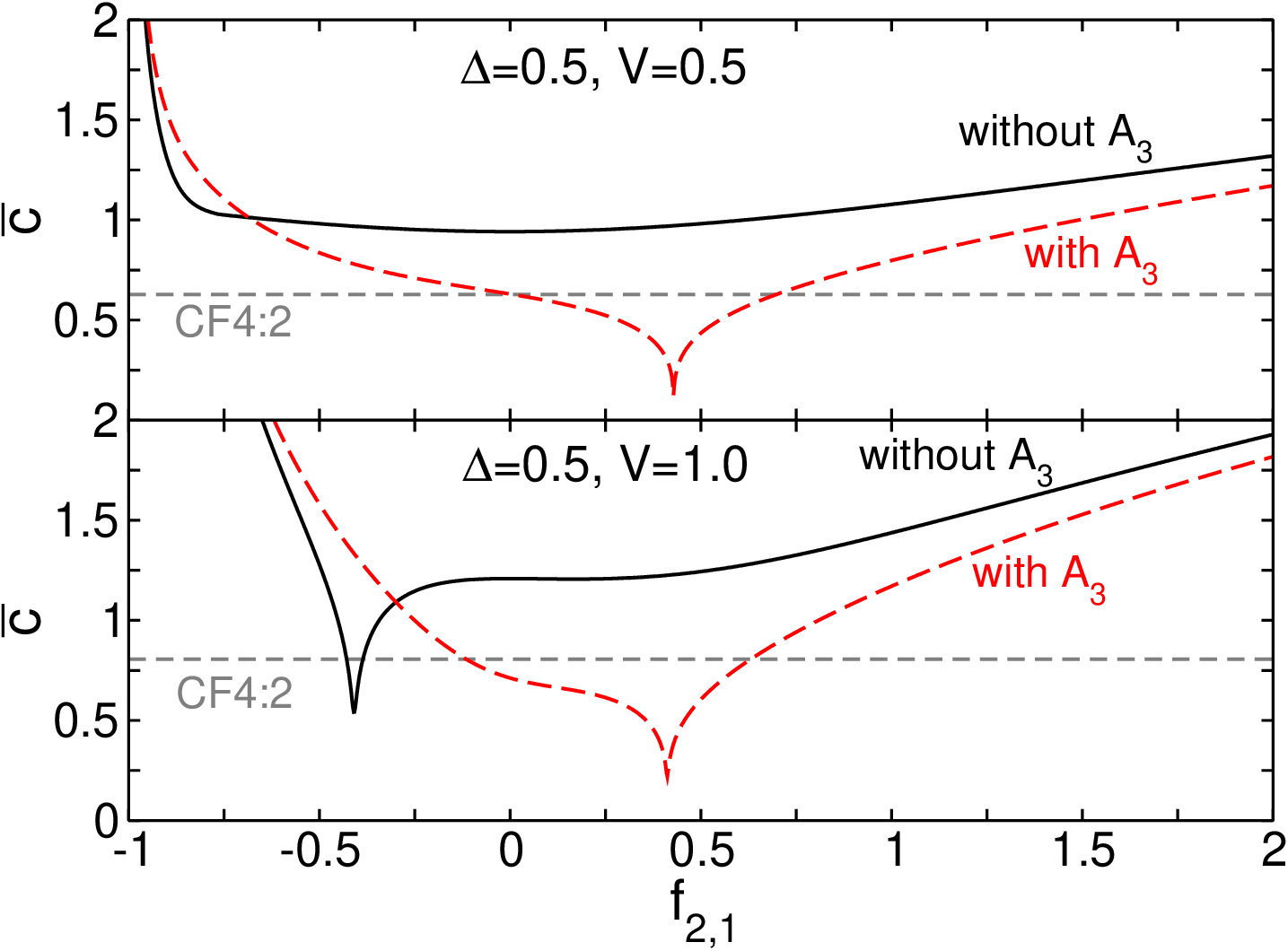}
\end{center}
\caption{Effective error constant $\bar{c}$ of 4th-order CFETs with 3 exponentials, given as a function of the free parameter $f_{2,1}$ in Eq.~\eqref{4thOptAnsatz}.
As explained in the text, the driven two-level system (Eq.~\eqref{TwoLevelH}) is propagated over the time $0 \le t \le T=20\pi/\omega$,
for parameters $\omega=1$, $\Delta=0.5$ and $V=0.5, 1.0$.
Shown are results for CFETs including the $A_3$ term for the optimal choice $f_{2,3}=-0.28$ (dashed red curve),
and without the $A_3$ term ($f_{2,3}=0$, solid black curve).
The horizontal dashed gray line gives $\bar{c}$ for the CFET CF4:2 with 2 exponentials (Eq.~\eqref{UCF4:2}), corresponding to $f_{2,1}=f_{2,3}=0$.
}
\label{fig:Err4thOnF21}
\end{figure}

In Fig.~\ref{fig:Err4thOnF21} we show $\bar{c}$ as a function of the free parameter $f_{2,1}$ used in Sec.~\ref{sec:4thOpt} for optimization of 4th-order CFETs with 3 exponentials.
In both cases (upper and lower panel) $\bar{c}$ is minimal for $f_{2,1} \simeq 0.45=9/20$, which confirms the previous theoretical analysis based on Fig.~\ref{fig:Chi}.
In comparison to CF4:2 with 2 exponentials, which has larger $\bar{c}$ than CF4:3Opt, we see that the error reduction is sufficiently large to outweigh the increased effort arising with an additional third exponential. We conclude that the optimization is successful and results in more efficient CFETs.

The importance of including the $A_3$ term becomes evident when dropping it, i.e. setting $f_{2,3}=0$ (solid black curve). 
Generally, $\bar{c}$ for such CFETs is large because of significant contributions from the $A_3$ terms in Eq.~\eqref{Delta4thOpt}, and the `optimal' value $f_{2,1}=9/20$ does not reduce the error.
Accidental cancellation of different terms occurs for certain parameter combinations and 
leads to the `dip' in $\bar{c}$ for $f_{2,1} \approx -0.4$ (lower panel).
Notice that in contrast to such artificial minima the true optimized value $f_{2,1} \approx 0.45$ gives a stable minimum of $\bar{c}$.

\begin{figure}
\begin{center}
\includegraphics[width=0.45\textwidth]{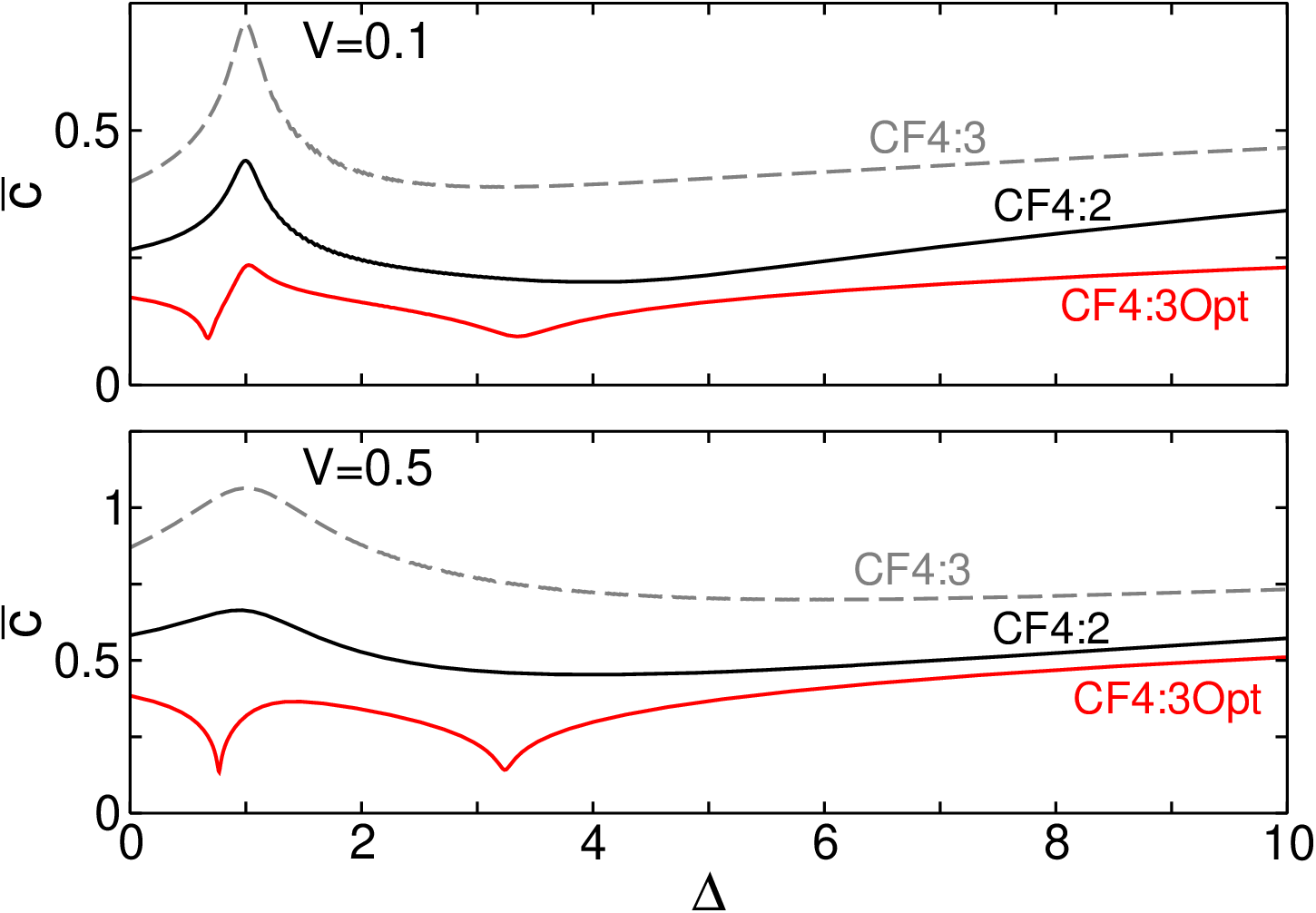} \hfill
\includegraphics[width=0.45\textwidth]{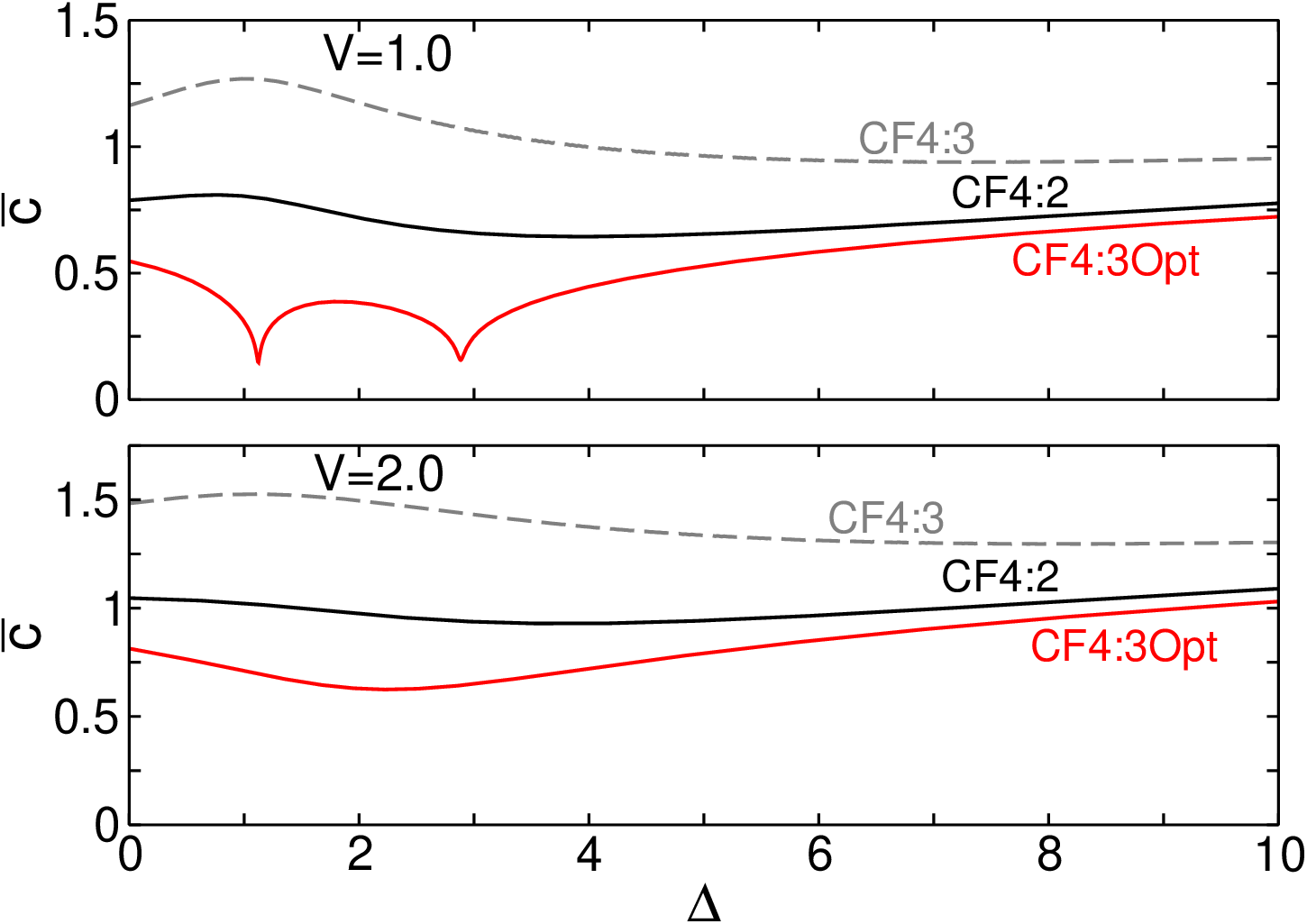}
\end{center}
\caption{Effective error constant $\bar{c}$ for unoptimized (CF4:2, CF4:3 from Table~\ref{4thCFET})
and optimized (CF4:3Opt from Table~\ref{4thCFETOpt}) 4th-order CFETs, as indicated.
As for the previous figure, the driven two-level system is propagated over $T=20\pi/\omega$.
Results are given as a function of $\Delta$, with $\omega=1$ and $V=0.1, 0.5, 1.0, 2.0$ as indicated in the panels.
}
\label{fig:Err4th}
\end{figure}

In Fig.~\ref{fig:Err4th} we show $\bar{c}$ for a range of parameter combinations of the driven two-level system.
Again we see that the optimization of CF4:3Opt is successful and results in smaller values of $\bar{c}$.
As an estimate, CF4:3Opt is about $10\%$ to $50\%$ more efficient than CF4:2.
Optimization attempts without the $A_3$ terms (CF4:3) result in reduced efficiency.

\subsection{Sixth-order CFETs}\label{sec:6thErr}

\begin{figure}
\begin{center}
\includegraphics[width=0.45\textwidth]{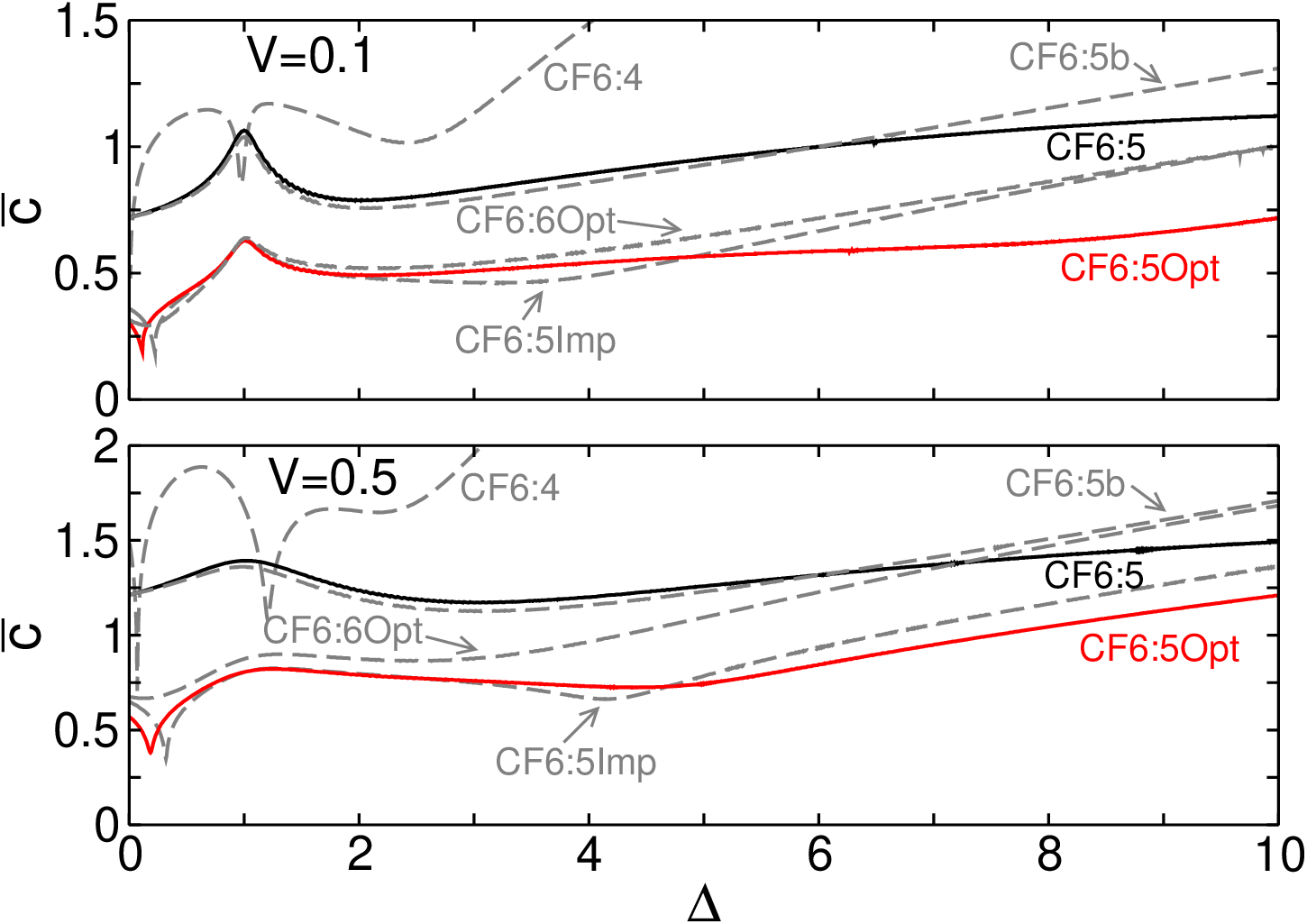} \hfill
\includegraphics[width=0.45\textwidth]{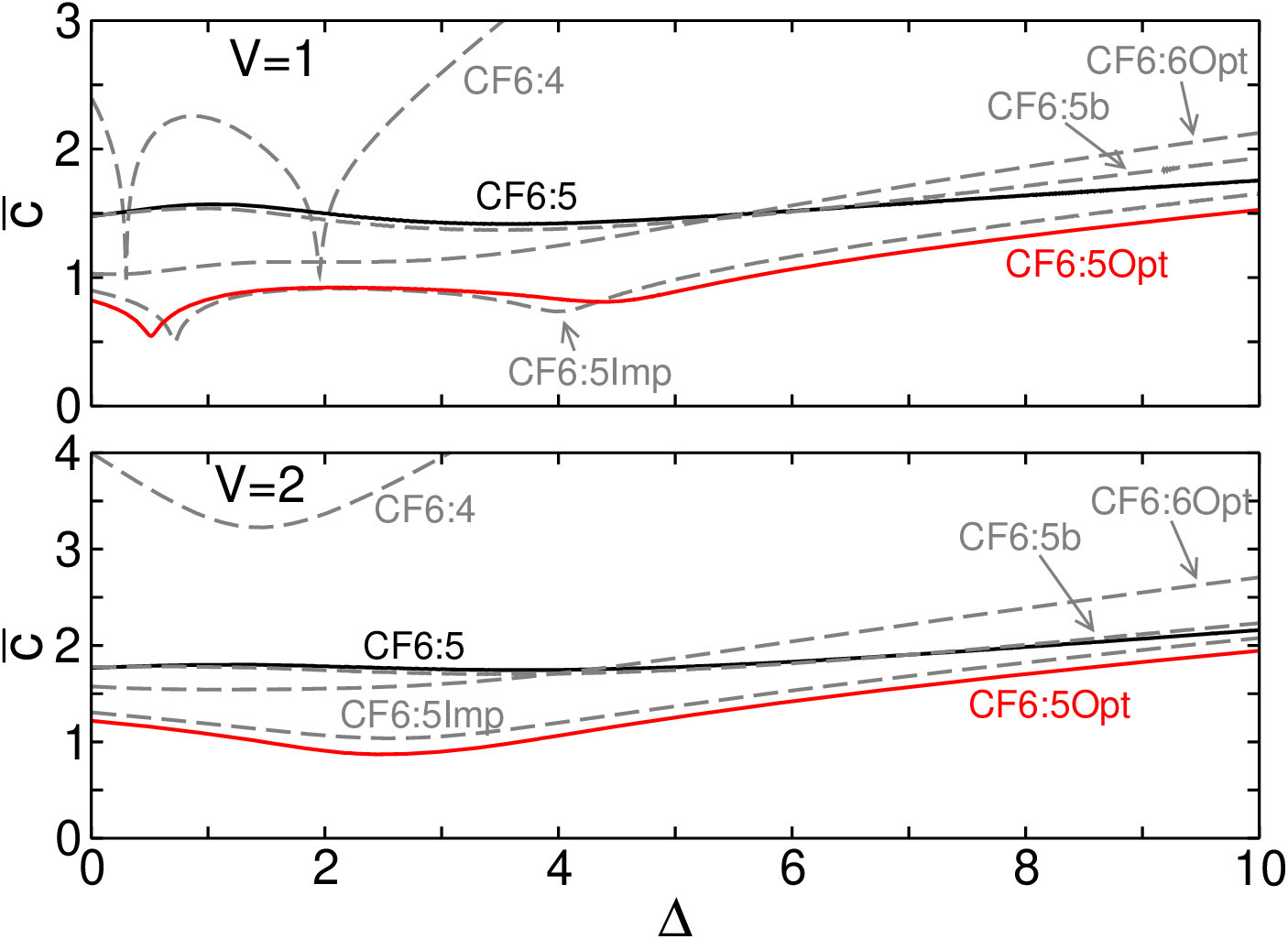}
\end{center}
\caption{Effective error constant $\bar{c}$ for different unoptimized (CF6:4, CF6:5, CF6:5b from Table~\ref{6thCFET}) and optimized (CF6:5Imp, CF6:5Opt, CF6:6Opt  from Table~\ref{6thCFETOpt}) 6th-order CFETs, as indicated.
The solid black (red) curve corresponds to CF6:5 (CF6:5Opt).
The propagation parameters are identical to Fig.~\ref{fig:Err4th}.}
\label{fig:Err6th}
\end{figure}

In Fig.~\ref{fig:Err6th} we show the effective error constant $\bar{c}$ for different 6th-order CFETs.
We can draw similar conclusions as for the 4th-order CFETs.
Since the parameter $f_{1,1}=0.2$ of the CFET CF6:5b from Ref.~\cite{BM06} is close to the optimal choice $f_{1,1}=0.16$ of CF6:5, both CFETs are comparable, with a slight advantage for CF6:5. 
The CFET CF6:4 is much less efficient, although it requires only $4$ exponentials.
The optimized CFET CF6:5Opt is generally the most efficient,
while dropping the $A_4$ term (as in CF6:5, CF6:5b) reduces the efficiency.
Notice that CF6:5Imp, including the $A_4$ term into CF6:5, is not as efficient as the fully optimized CF6:5Opt, but still significantly better than the other CFETs.
The additional freedom of choice of parameters for $6$ exponentials (CF6:6Opt) does not lead to further reduction of $\bar{c}$. 

\subsection{Comparison of CFETs of different order}\label{sec:MethOrd}

According to Eq.~\eqref{thumb}, time-propagation with smaller error, i.e. higher accuracy demands, is more efficient using higher-order CFETs.
A given $N$th-order CFET is most efficient in a certain `accuracy window', whose size depends on the respective error constant $\bar{c}$ and propagation time $T$.
The intended accuracy goal thus suggests a preferential choice of $N$ and the corresponding optimized CFET.

Consider two CFETs of order $N_1 < N_2$, with effective error constants $\bar{c}_1$, $\bar{c}_2$. 
Inverting the effort-error relation from Sec.~\ref{sec:EffErr}, we find that the $N_1$-order CFET is more efficient than the $N_2$-order CFET if 
\begin{equation}\label{ErrOrder}
  \Big( \frac{\bar{c}_1}{\bar{c}_2} \Big)^{\tfrac{N_1 N_2}{N_2-N_1}} \le \frac{\epsilon}{T} \;.
\end{equation}
The decisive quantity is the ratio $\epsilon/T$ of the maximal acceptable error $\epsilon$ and the propagation time $T$.

For a rough estimate, let us assume that the effective error constants $\bar{c}_1$, $\bar{c}_2$ are given by the number $s_1$, $s_2$ of exponentials.
With $s=1,2,5,11$ for $N=2,4,6,8$ we find the following values, which provide some orientation:
\begin{center}
\begin{tabular}{lccccccc}
\hline\noalign{\smallskip}
error $\epsilon/T$: & $\dots$ & $6 \times 10^{-2}$ & $\dots$ & $2 \times 10^{-5}$ & $\dots$ & $6 \times 10^{-9}$ & $\dots$ \\ \hline\noalign{\smallskip}
favourable $N$:  & 2 & $|$ & 4 & $|$ & 6 & $|$ & 8 \\\hline
\end{tabular}
\end{center}
As a rule of thumb, the accuracy window spans three orders of magnitude: 
4th-order CFETs are good for low (error $10^{-3}$),
6th-order for moderate (error $10^{-6}$), and 8th-order for high (error $10^{-9}$) accuracy demands. The use of 2nd-order CFETs such as the midpoint rule should be avoided.
Long propagation times shift the advantage towards higher-order CFETs.

\begin{figure}
  \begin{center}
    \includegraphics[width=0.45\textwidth]{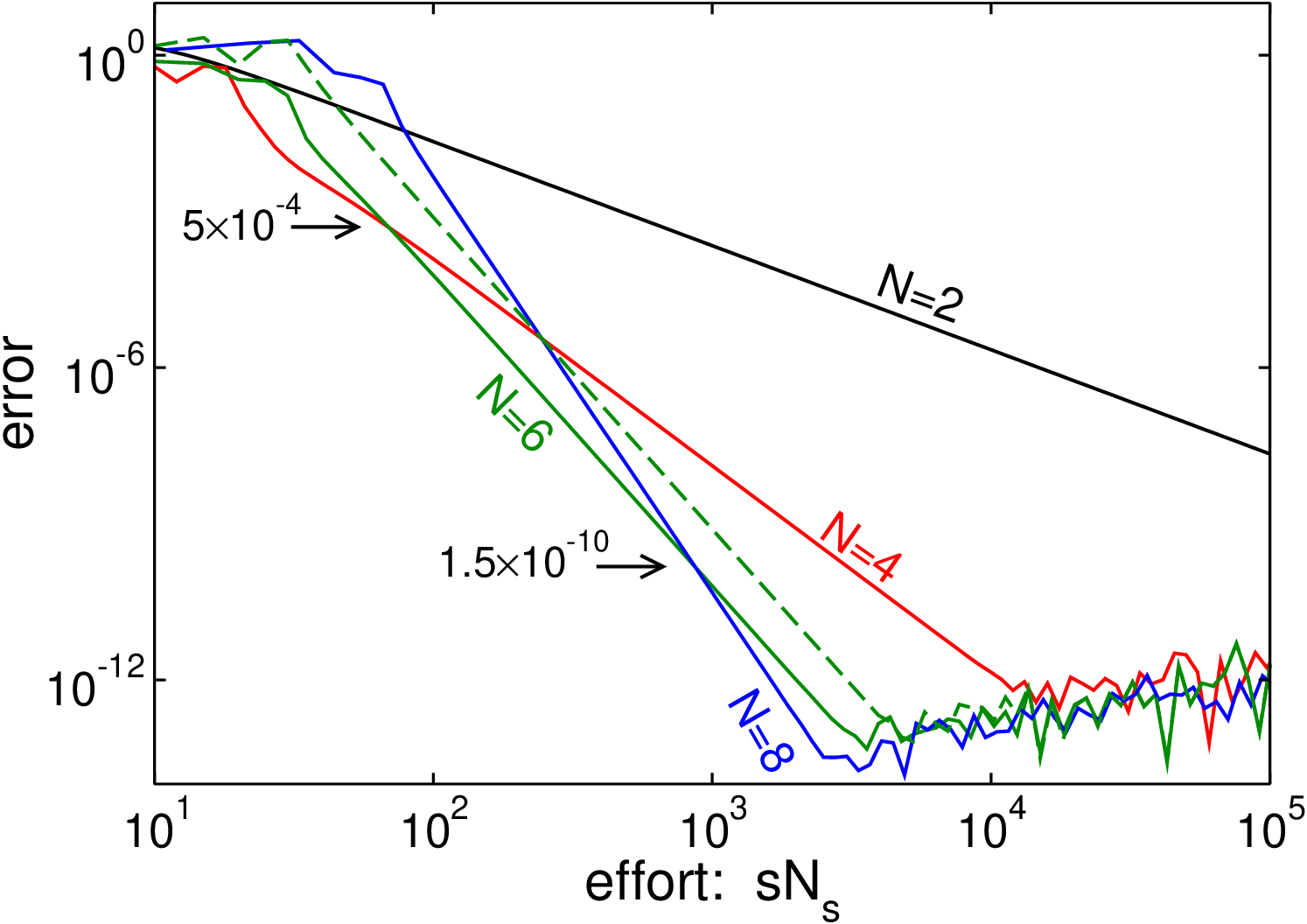}
    \hfill
    \includegraphics[width=0.45\textwidth]{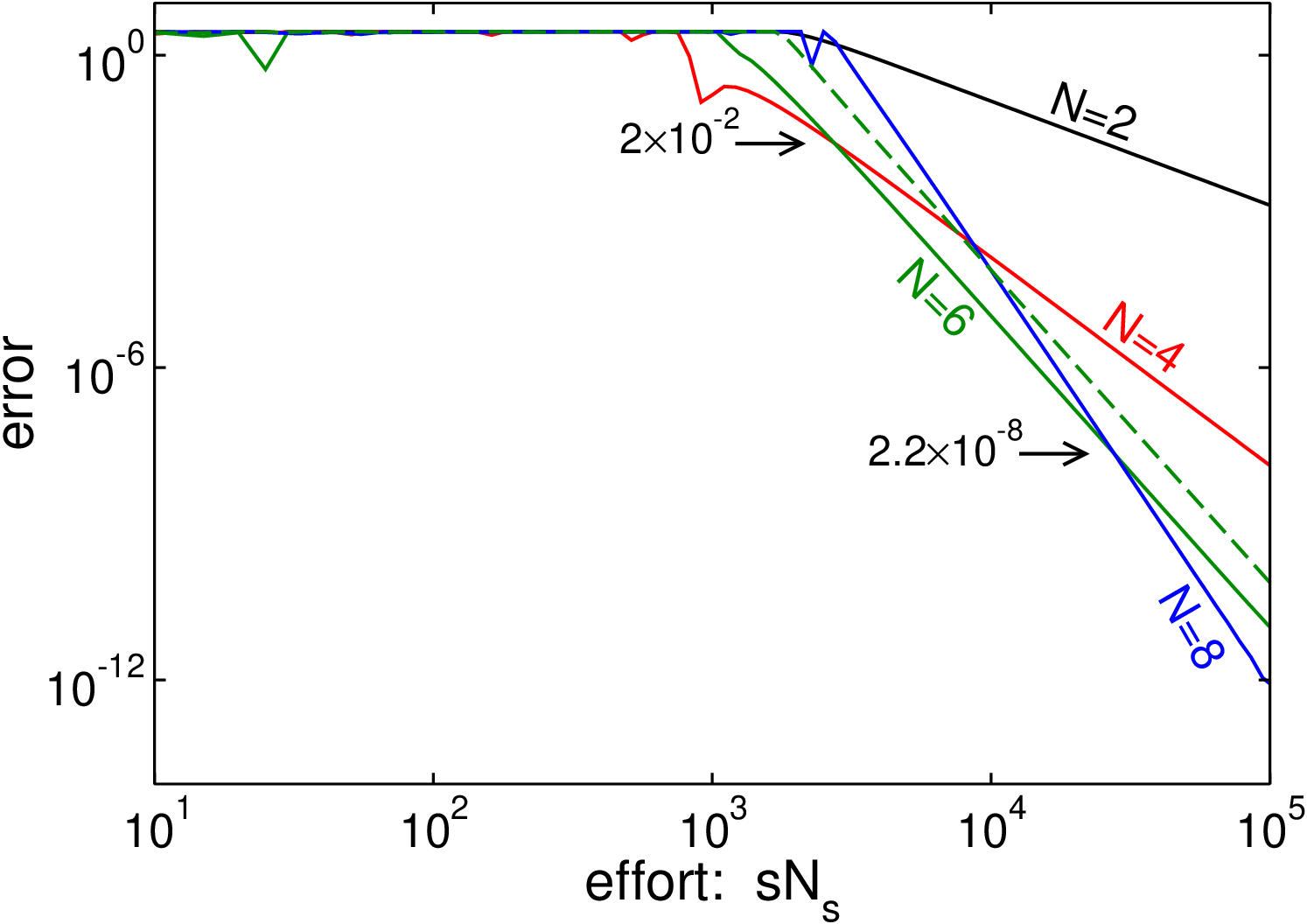}
  \end{center}
\caption{Error-effort plot for the $N$th-order CFETs CF2:1, CF4:3Opt, CF6:5Opt and CF8:11 as indicated,
applied to the driven two-level system Eq.~\eqref{TwoLevelH} with $\omega=1$, $\Delta=2$, $V=0.5$. The system is propagated over $T=5\pi/\omega$ (left panel) or $T=200\pi/\omega$ (right panel).
Also included is the error of the unoptimized CFET CF6:5 (dashed green line).
} 
\label{fig:ErrorScal}
\end{figure}

For a case study we show in Fig.~\ref{fig:ErrorScal} the error-effort plot for 2nd- to 8th-order CFETs,  applied to the two-level system from Sec.~\ref{sec:2level}, for short (left panel) and long (right panel) propagation time.
Notice that a very small error can be achieved before it saturates at machine precision.
The accuracy window of the $N$th-order CFET is bounded by the crossing with the $N\pm 2$ curves. 
For a moderate error $10^{-7}$ (the square root of machine precision for FORTRAN double precision numbers), switching from the 4th- to the 6th-order CFET reduces the effort by a factor of 2--3.
For longer propagation times (right panel) the accuracy window shifts to larger errors,
and the 8th-order CFET becomes more efficient.
The performance of the 2nd-order midpoint rule is several orders of magnitude worse.
To illustrate the benefit of optimization, we include results for the unoptimized CFET CF6:5. As can be seen, it is never competitive in comparison to the (optimized) 4th- or (unoptimized) 8th-order CFET.

\subsection{Time-step selection}\label{sec:MethStep}

In practice a prescribed accuracy goal has to be achieved without knowledge of the exact solution of the problem.
A simple, conservative approach is to perform calculations with an ever decreasing time step $\dt$ until convergence, i.e. two subsequent calculations agree within numerical round-off errors.
This approach wastes much computational time if we seek less accurate results, as it tries to construct the (numerically) exact solution.

For a better time-step selection we can
use the known scaling of the error as $\epsilon = c T \dt^N$. 
An estimate of the error constant $c$ is obtained from two calculations with different time steps $\dt_1$, $\dt_2$
according to the relation
\begin{equation}
 \max{ \|x_1(t) - x_2(t)\| } = c T  |\dt_1^N-\dt_2^N| \;.
\end{equation}
It involves only the difference between the two approximate solutions $x_1(t)$, $x_2(t)$, but not the unknown exact solution.
A reasonable choice is $\dt_1/\dt_2 = (2 \dots 3)^{1/N}$, such that the error decreases by a small but significant amount.
From the estimate of $c$ we can extrapolate to the required time step $\dt$ for the given accuracy goal.
For a reliable estimate of the final error it is recommended to perform an additional calculation with smaller $\dt$.
An alternative is to compare numerical solutions obtained with two CFETs of different order.

For applications where the time-dependence of $A(t)$ does not change significantly with $t$,
the required time-step can be determined for some finite period $0 \le t \ll T$ that is characteristic for the dynamical evolution of the system.
The solution over the entire propagation time $[0,T]$ is then computed with the fixed, predetermined value of $\dt$.
In other situations, we can proceed similar to heuristic strategies for general differential equation solvers~\cite{PFTV86,Iser09}, which achieve the global accuracy goal through control of the local time-stepping error. If the above extrapolation for $\dt$ is performed at every step, it allows for propagation with adaptive time-step selection.

\section{Gauss-Legendre quadrature}\label{sec:GaussLeg}

In numerical applications the terms $A_n$ from the Legendre expansion Eq.~\eqref{An} can be calculated with a numerical quadrature formula. 
For an optimized $N$th-order CFET the quadrature formula must be of order $N+1$. A convenient choice is Gauss-Legendre quadrature~\cite{PFTV86} with $N/2+1$ quadrature points.

\begin{table}
\begin{tabular}{cccccc}
   & \multicolumn{5}{c}{m}  \\
 \raisebox{1.5ex}[0cm][0cm]{M}  & & 1 & 2 & 3 & 4 \\\hline\hline
  & $x_m$ & 1/2 \\
  \raisebox{1.5ex}[0cm][0cm]{1} & $w_m$ & 1 \\\hline\noalign{\smallskip}
 & $x_m$ & $ 1/2 - \sqrt{3}/6$ & $ 1/2 + \sqrt{3}/6$ \\
 \raisebox{1.5ex}[0cm][0cm]{2} & $w_m$ & 1/2 & 1/2 \\\hline\noalign{\smallskip}
 & $x_m$ & $ 1/2 -\sqrt{3/20}$   & 1/2 & $ 1/2 +\sqrt{3/20}$ \\
 \raisebox{1.5ex}[0cm][0cm]{3} & $w_m$ & 5/18 & 4/9 & 5/18 \\\hline\noalign{\smallskip}
 & $x_m$ & $1/2 -\! \sqrt{\smash[b]{(3\!+\!2\sqrt{6/5})/28}}$ & $1/2 -\! \sqrt{\smash[b]{(3\!-\!2\sqrt{6/5})/28}}$ & $1/2 +\! \sqrt{\smash[b]{(3\!-\!2\sqrt{6/5})/28}}$ &  $1/2 +\! \sqrt{\smash[b]{(3\!+\!2\sqrt{6/5})/28}}$ \\
 \raisebox{1.5ex}[0cm][0cm]{4} & $w_m$ & $(18-\sqrt{30})/72$ & $(18+\sqrt{30})/72$ & $(18+\sqrt{30})/72$ & $(18-\sqrt{30})/72$ \\\hline
\end{tabular}
\caption{Points and weights for Gauss-Legendre quadrature over $[0,1]$ up to order $8$, see Eq.~\eqref{GLeg}.
Generally, $x_{M+1-m}=1-x_m$ and $w_{M+1-m}=w_m$.}
\label{tab:GaussLeg}
\end{table}

Gauss-Legendre quadrature is specified through $M$ points $x_1,\dots,x_M$,  which are the zeros of the Legendre polynomial $P_M(x)$,
and weights $w_1,\dots,w_M$ (see Table~\ref{tab:GaussLeg}).
The integral of a function $f(x)$ is approximated as
\begin{equation}\label{GLeg}
\int_0^1 f(x) dx \approx \sum_{m=1}^M w_m f(x_m) \;.
\end{equation}
Using the orthogonality of Legendre polynomials it can be shown that 
Gauss-Legendre quadrature is of order $2M$,
in the sense that this expression is exact for polynomials with maximal degree $2M-1$.
Equivalently, the error of the approximation $\int_0^\dt f(t) dt \approx \dt \sum_{m=1}^M w_m f(x_m \dt)$ scales as $\dt^{2M+1}$. 

For the integrals in Eq.~\eqref{An} Gauss-Legendre quadrature with $M=N/2+1$ points gives
\begin{equation}\label{AnGauss}
  A_n \simeq (2n-1) \, \dt \sum_{m=1}^{N/2+1} w_m P_{n-1}(x_m) A(x_m \dt)  
\end{equation}
for the terms $A_1, \dots, A_{N/2+1}$ of an optimized $N$th-order CFET.
This expression can be inserted into Eq.~\eqref{Omi} to obtain
\begin{equation}\label{OmGLeg}
 \Omega_i = \dt \sum_{m=1}^{N/2+1} g_{i,m} A(x_m \dt)
\end{equation}
as a linear combination of $A(t)$ at different times $x_m \dt$ in $[0,\dt]$,
where the new coefficients are 
\begin{equation}
  g_{i,m} = w_m \sum_{n=1}^{N/2+1} (2n-1) P_{n-1}(x_m) f_{i,n} \;.
\end{equation}
We note that, using Legendre polynomials, the calculation of the $g_{i,m}$ from the tabulated $f_{i,n}$ is much simpler than for an expansion in powers of $\dt$ (cf.~Refs.~\cite{BM06,BCOR09}).
Specifically for the CFET CF4:2 from Eq.~\eqref{UCF4:2} we obtain
\begin{equation}\label{UCF4Gauss}
   \tilde{U}^{(4)}_\mathrm{CF4:2} = \exp \Big[ \dt \left( \frac{3-2\sqrt{3}}{12} A^{(1)} + \frac{3+2\sqrt{3}}{12}  A^{(2)} \right) \Big] \,
                                    \exp \Big[ \dt \left( \frac{3+2\sqrt{3}}{12} A^{(1)} + \frac{3-2\sqrt{3}}{12}  A^{(2)} \right) \Big] \;, 
\end{equation}
where $A^{(1)} = A[ (1/2-\sqrt{3}/6) \dt]$, $A^{(2)} = A[ (1/2+\sqrt{3}/6) \dt]$.

It remains to show that Gauss-Legendre quadrature with $N/2+1$ points correctly reproduces the $\Omega_i$.
If we insert the expansion Eq.~\eqref{Aexpansion} into Eq.~\eqref{AnGauss},
we find that the $A_n$ are approximated as
\begin{equation}
  A_n \simeq (2n-1) \sum_{l \ge1} A_l \sum_{m=1}^{N/2+1} w_m P_{n-1}(x_m) P_{l-1}(x_m) \;.
\end{equation} 
The summands on the right hand side are the $N+2$-order Gauss-Legendre approximations 
\begin{equation}
\int_0^1 P_{n-1}(x) P_{l-1}(x) dx \simeq \sum_{m=1}^{N/2+1} w_m P_{n-1}(x_m) P_{l-1}(x_m)
\end{equation}
of the scalar product of Legendre polynomials.
As long as $(n-1)+(l-1) \le N+1$, i.e. $n+l \le N+3$, the approximation is exact and gives the correct value $\delta_{nl}/(2n-1)$.
In particular for $n=1,2$, all integrals for $1 \le l \le N+1$ are evaluated correctly, and Gauss-Legendre quadrature constructs the terms $A_1, A_2$ with an error of order $\dt^{N+2}$, as required for an optimized $N$th-order CFET. 
For $n \ge 3$, the integrals with $l > N-(n-3)$ are not evaluated correctly, and introduce an error of order $\dt^{N+4-n}$ into the term $A_n$.

To understand why the CFET order is nevertheless preserved we must revisit the property (P2) discussed in Secs.~\ref{sec:PropExp},~\ref{sec:OrdCond}. 
It states that every nested commutator $[A_{n_1}, \dots, A_{n_m}]$ contributing to $\Omega(t)$ fulfills the condition $n_k \le 1+\sum_{i\ne k} n_i$ for all $k=1,\dots,m$.
By rule (R2) for the CFET construction this property carries over to the approximate $\tilde{\Omega}(t)$ from Eq.~\eqref{tildeOm}.
Since the error of the term $A_{n_k}$ incurred from numerical quadrature is of the order $\dt^{N+4-n_k}$, the error of the nested commutator is of order $\dt^{N+4-n_k + \sum_{i \ne k} n_i}$ due to the multiplication with the remaining terms. By the above condition this is at least of order $\dt^{N+3}$, as required.

We note that the above argumentation shows the intrinsic connection between Gauss-Legendre quadrature and the property (P2) about the absence of certain nested commutators from $\Omega(t)$. The connection is established through expansions in orthogonal Legendre polynomials.
Of practical interest is that Gauss-Legendre quadrature with $N/2+1$ points suffices for (optimized) $N$th-order CFETs, although in principle most terms $A_n$ are reproduced with an error of lower order.


\section{Implementation}\label{sec:Implement}

Due to the simple product form of Eq.~\eqref{Ucf} the application of CFETs is straightforward.
The single difficult numerical part is the evaluation of the matrix exponentials $e^{\Omega_i}$,
which is possible with the Krylov technique discussed below.
Using Gauss-Legendre quadrature each $\Omega_i$ is obtained from $A(t)$ as a weighted sum (Eq.~\eqref{OmGLeg}).
For large-scale problems, where $A(t)$ is a sparse matrix, it implies that also the $\Omega_i$ are sparse. Moreover the sparsity pattern of $A(t)$, i.e. the distribution of nonzero entries, is preserved: Zeros add up to zeros.
This allows for seamless integration of CFETs into existing programs, 
which implement specific data storage formats or matrix-vector multiplication routines~\cite{HW10}. The extension to time-dependent Hamilton operators requires only minor modifications.
The feature of easy implementation gives CFETs additional advantage over the original Magnus expansion.

\subsection{Calculation of exponentials}\label{sec:Krylov}

Two powerful approaches for the computation of matrix exponentials,
particularly of $e^{-\ii M}$ with sparse hermitian matrices $M$,
are the Krylov~\cite{Sid98,HL97} and the Chebyshev technique~\cite{TK84}.
Both techniques calculate $e^{-\ii M} \psi$, the exponential applied to a vector, iteratively. They avoid diagonalization of the matrix $M$, which enters only through matrix-vector multiplication
as required for sparse matrices.
If $M$ is a sum of simple terms, split-operator techniques~\cite{MQ02} can reduce the computational effort considerably. 
Other methods, such as the 2nd-order Crank-Nicholson approximation $e^{-\ii M} = (1- \ii M)/(1+\ii M)$, are neither very accurate, nor suitable for large-scale problems~\cite{ML03}.

The Chebyshev technique is based on the expansion of the exponential function $e^{\ii x t}$ in a series of Chebyshev polynomials.
Similar to the calculation of spectral functions~\cite{WWAF06}, it has the advantage of low memory demands, simple implementation, and unconditional stability and concomitant accuracy for arbitrary large propagation times.
The main disadvantage, especially for time-dependent Hamilton operators, is the need to determine a-priori bounds on the eigenvalues of the matrix $M$.

The Krylov technique is based on the Lanczos iteration. Starting with the initial vector $\psi$, each multiplication with $M$ gives a new vector $M^k\psi$, which is orthogonalized to the previous vectors from the iteration.
A few $K$ iterations generate an orthogonal basis of the Krylov  subspace spanned by the vectors $\psi, M\psi, M^2 \psi, \dots , M^{K-1} \psi$.
The exponential $e^{-\ii M} \psi$ is approximately evaluated within the low-dimensional  Krylov subspace, which effectively reduces the problem to the calculation of an exponential of a dense $m\times m$ matrix~\cite{ML03}.
The success of this procedure depends on the quality of the Krylov approximation of $M$.
For the calculation of the exponential $e^{-\ii M}$, the error bound
\begin{equation}\label{KryErr}
	\mathsf{error} \le const. \times  e^{-\rho^2/K} \left(\frac{e \rho}{K}\right)^K \qquad (2\rho \le K),
\end{equation}
where the constant is independent of $K$ and $\rho$, can be established~\cite{HL97}.
Here, 
$4\rho=\lambda_\mathrm{max}-\lambda_\mathrm{min}$  is the spread of the maximal and minimal eigenvalue $\lambda_\mathrm{max}$, $\lambda_\mathrm{min}$ of $M$.
Increasing $K$ leads to fast reduction of the error. However, the finite main storage restricts the size of $K$. Therefore, the Krylov technique requires time-stepping, based on the equality $e^{-\ii H n \dt} = (e^{-\ii H \dt})^n$, if the error is not sufficiently small after a single Lanczos iteration.  
For fixed $K$, the Krylov approximation error for $e^{-\ii H \dt}$ is of order $\dt^K$.

\begin{figure}
\begin{center}
\includegraphics[width=0.45\textwidth]{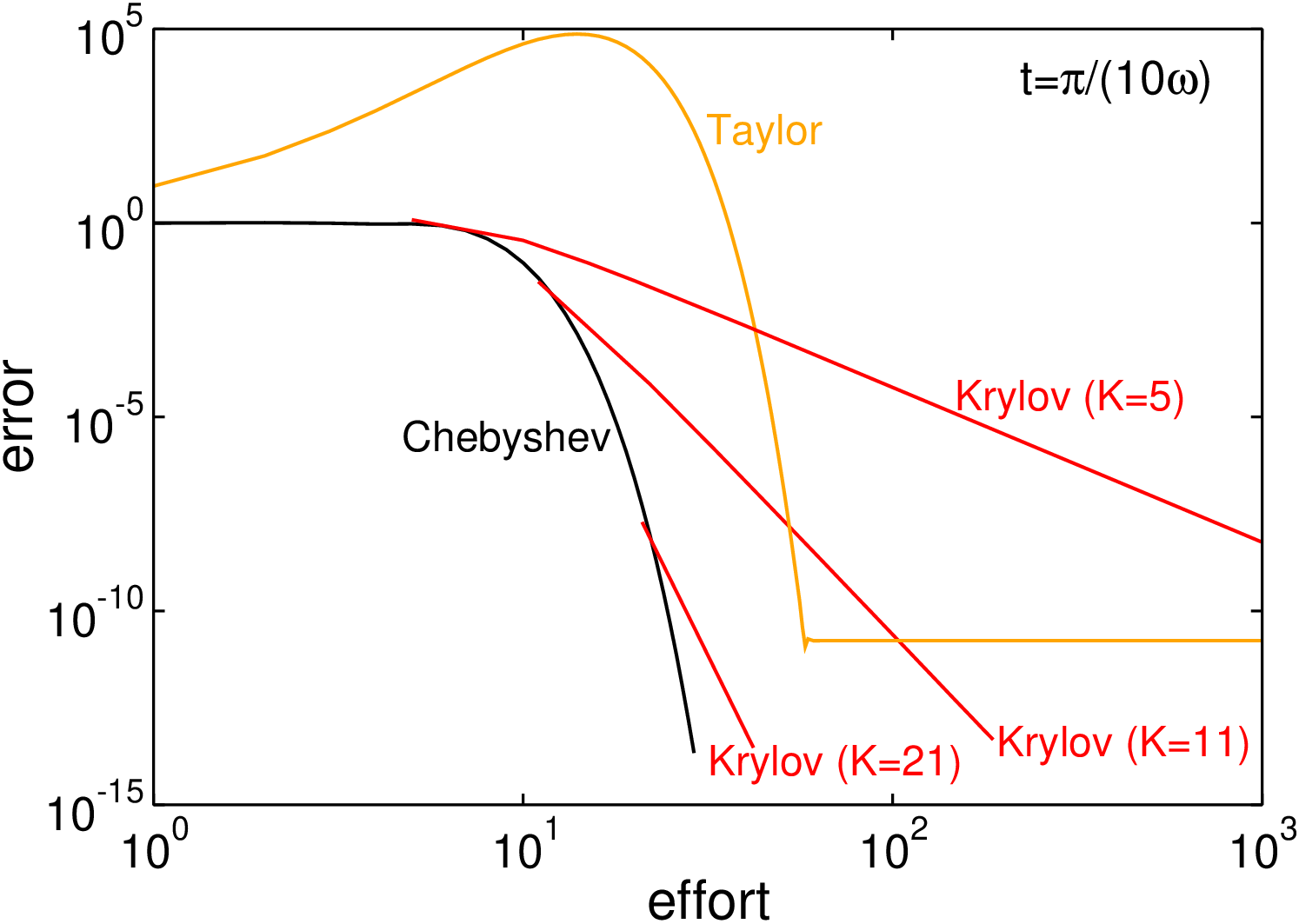}
\hfill 
\includegraphics[width=0.45\textwidth]{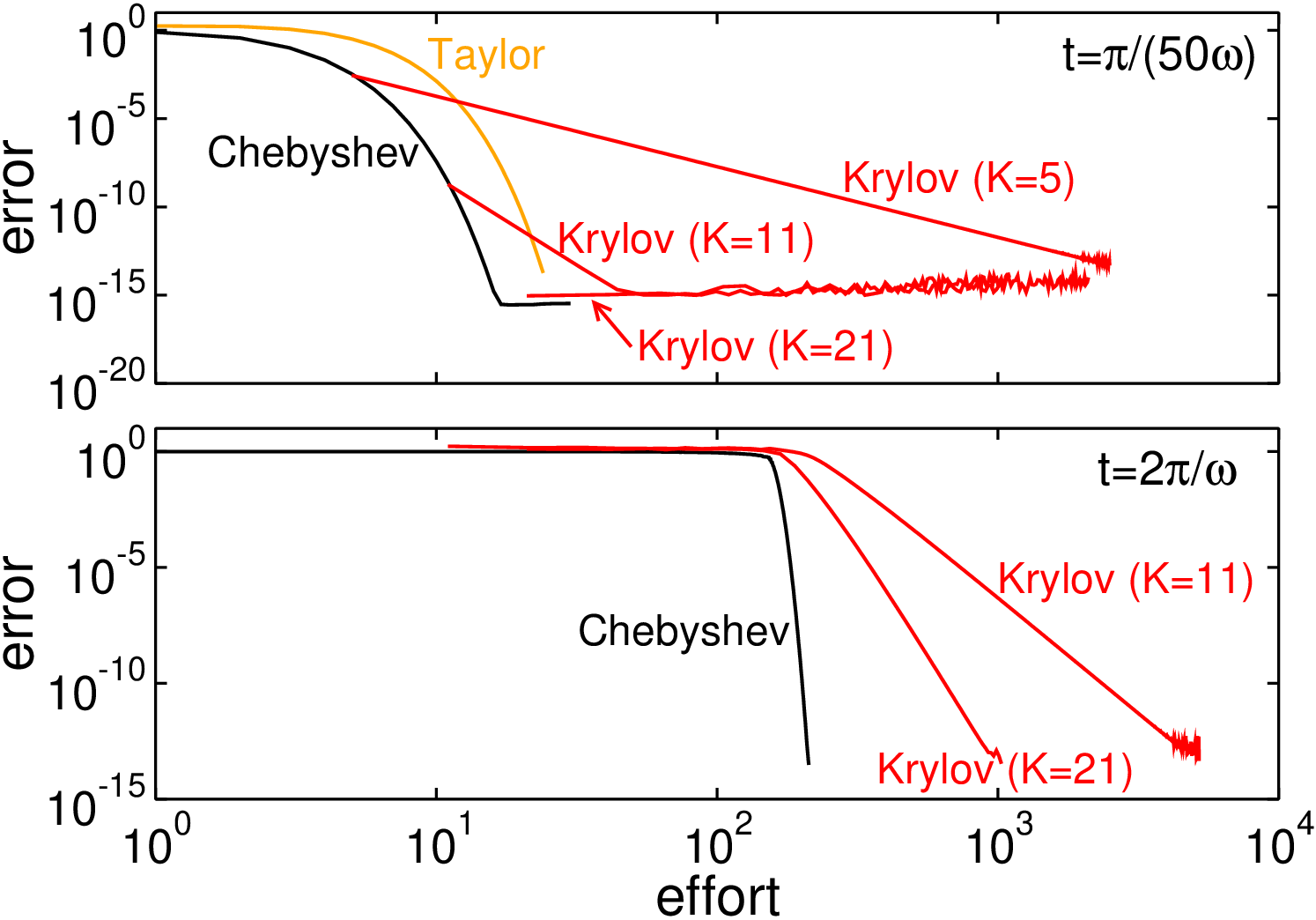} 
\end{center}
\caption{Comparison of the Krylov and Chebyshev technique with the Taylor expansion for the calculation of the exponential $e^{- \ii t  H} \psi$, for the Hamilton operator of the harmonic oscillator ($H_{nn}=n\omega$) as explained in the text.
Left panel: Propagation over a 20th oscillator period ($t=\pi/(10\omega)$). 
Right panel: Propagation over a 100th ($t=\pi/(50\omega)$) and a full oscillator period ($t=2\pi/\omega$).
Notice the asymptotic decay of the Krylov error $\sim \mathsf{effort}^{K-1}$.}
\label{fig:NumExp}
\end{figure}

\subsection{Comparison of the Krylov and Chebyshev technique}\label{sec:Krylov2}

In Fig.~\ref{fig:NumExp} we compare the Krylov and Chebyshev technique with an $m$th-order Taylor expansion of the exponential $e^{- \ii t H} \psi$,
where $H$ is the diagonal $50 \times 50$ matrix with entries $H_{nn} = n\omega$ and the vector elements $\psi_n$ are chosen at random (prior to normalization of $\psi$). 
 This corresponds to time-propagation for the quantum harmonic oscillator (cf. Sec.~\ref{sec:ParaOsci}). 
The error is given by the $\ell^2$-norm $\epsilon=|\psi(t) - \tilde{\psi}(t)|$
between the numerical result $\tilde{\psi}(t)$ and the exact $\psi(t)$, here with $\psi_n(t) = e^{-\ii t \omega n} \psi_n$. The effort is equal to the number of evaluations of $H \psi$ in the computation.
For the Chebyshev technique and Taylor expansion, which evaluate the exponential at once,
this corresponds to the number of terms kept in the respective series.
These definitions of error and effort are also used in the following examples. 
 
For the left panel in Fig.~\ref{fig:NumExp}, the system is propagated for $t=\pi/(10\omega)$, i.e. a 20th of the oscillator period.
The plot shows the typical problems of the Taylor expansion, whose error grows initially before it saturates far above machine precision.
Since the Taylor expansion violates unitarity, the large errors of the exponential spoil the stability of time-propagation.
Notice that an $N$th-order Runge-Kutta method applied to $e^{-\ii H t}$ is equivalent to using the Taylor expansion, which explains their diminished usefulness for quantum systems.
For the Chebyshev technique, the error decays fast after the first $10-20$ terms. Unitarity is 	again achieved only at the level of machine precision, which however now can be reached easily.
The Krylov technique is competitive for sufficiently many Krylov vectors in the iteration ($K \gtrsim 10$) and moderate accuracy demands.
Notice that the eigenvalues of the quantum harmonic oscillator occur as multiples of $\omega$,
which leads to a large eigenvalue spread $\rho$ in Eq.~\eqref{KryErr})
and increases the computational effort more than the `classical' time-scale $1/\omega$ may suggest. 
 
The Krylov technique becomes more efficient for small time-steps \dt.
It complements the Chebyshev technique which excels for longer propagation times.
Both scenarios are depicted in the right panel in Fig.~\ref{fig:NumExp}.
This makes the Krylov technique the more suitable choice for combination with CFETs, where the length of the time-step is restricted by the time-dependence of $A(t)$ (or $H(t)$).
Its central advantage, however, is that it strictly preserves unitarity even for finite error.
This allows us to dispense with the evaluation of the exponentials $e^{\Omega_i}$ to very high accuracy when the overall error is dominated by the CFET error.
Instead, we can use the Lanczos iteration with small $K$ (it must $K>N$ for an $N$th-order CFET).
The reduction of the time-step \dt, in order to decrease the CFET error, reduces the Krylov error at the same time.
While we recommend the use of the Krylov technique for CFETs we must also note
that the present example shows that the Chebyshev technique should not be finally dismissed
even for short-time propagation.

\section{Example: CFETs applied to the parametric harmonic oscillator}\label{sec:ParaOsci}

A genuine example for driven system is the quantum parametric harmonic oscillator
\begin{equation}\label{Hpara}
 H(t) = \frac{1}{2} \hat{p}^2 + \frac{\omega(t)^2}{2} \hat{q}^2 \;,
\end{equation}
where we allow for a time-dependent oscillator frequency $\omega(t)=\omega_0^2 + \xi \cos \Omega t$. Position $\hat{q}$ and momentum operator $\hat{p}$ obey the canonical commutation relation $[\hat{q},\hat{p}]=\ii$. 
The oscillator position $q(t)=\langle \hat{q} \rangle (t) = \langle\psi(t)|\hat{q}|\psi(t)\rangle$,
given as the expectation value of $\hat{q}$, follows the classical equation of motion -- the Mathieu equation --
\begin{equation}\label{Mathieu}
\ddot{q} + ( \omega_0^2+\xi \cos \Omega t  ) q = 0 \;.
\end{equation}

\subsection{Classical oscillator}

\begin{figure}
\begin{center}
\includegraphics[width=0.45\textwidth]{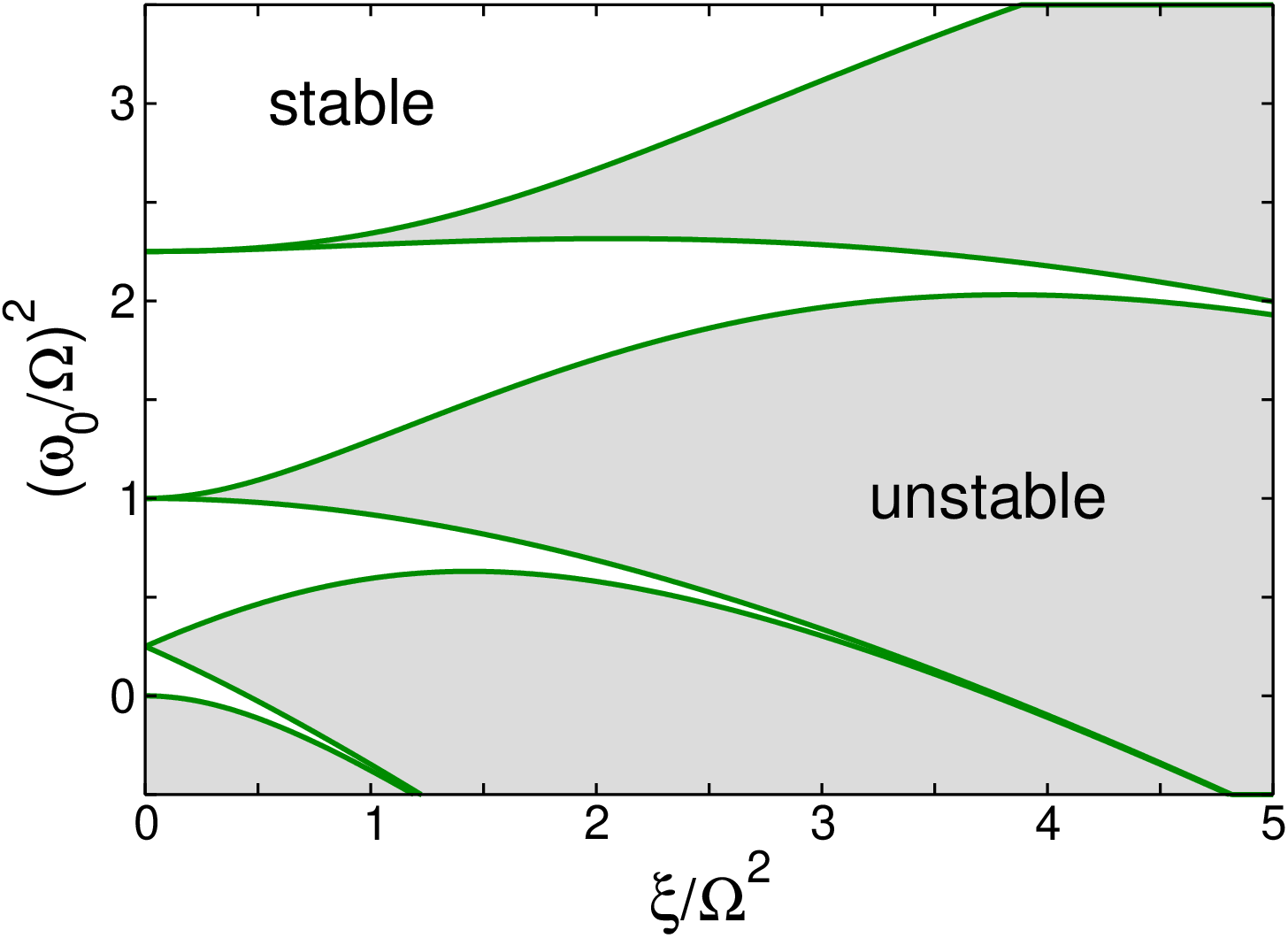} 
\hfill
\includegraphics[width=0.45\textwidth]{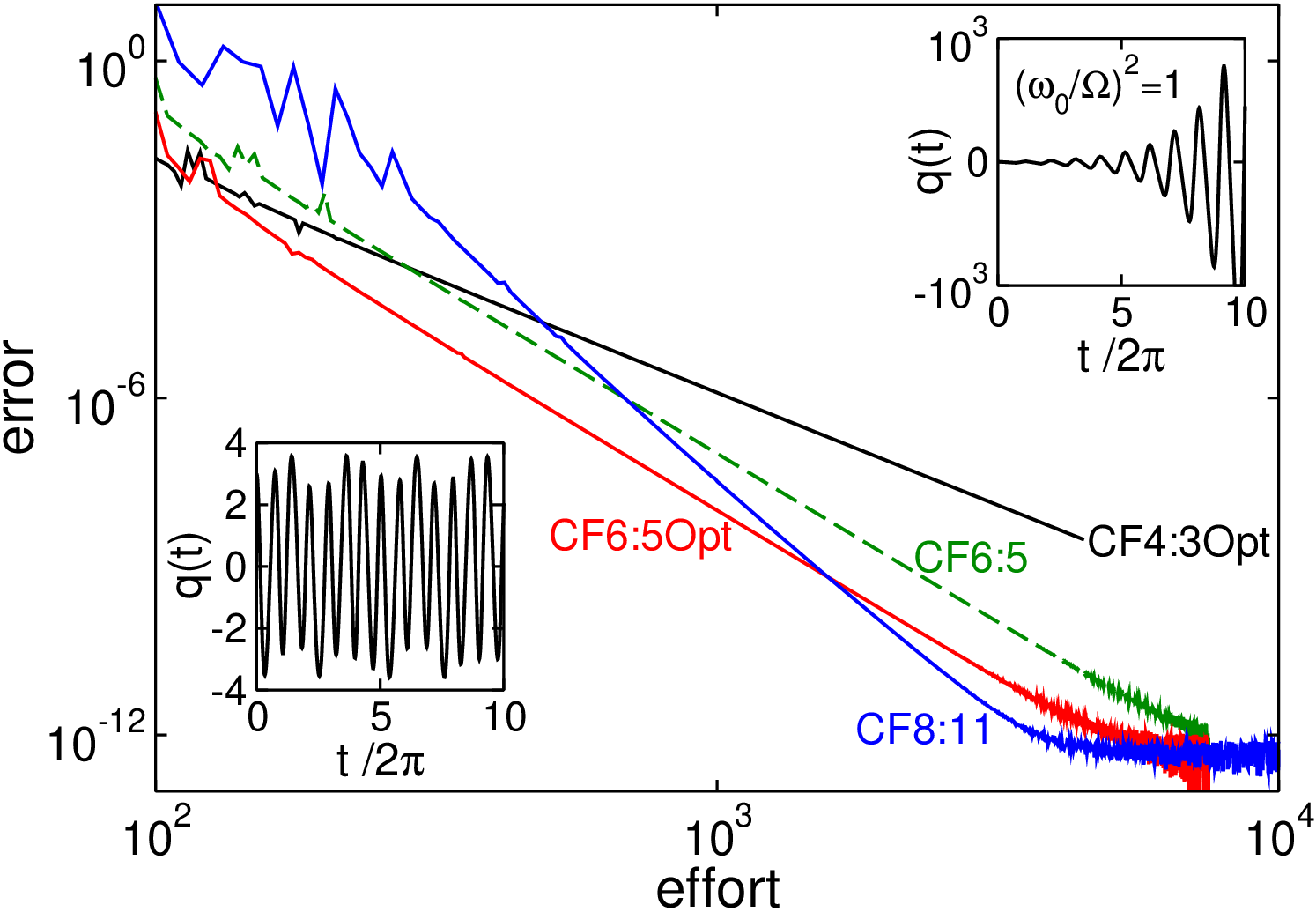}
\end{center}
\caption{
Left panel: Stability chart of the parametric harmonic oscillator, and the associated Mathieu equation Eq.~\eqref{Mathieu}. In the shaded regions of parameter space exponentially growing solutions exist, i.e. the system is unstable. For $\xi\to 0$, isolated instabilities occur at integer values of $2 \omega_0/\Omega$.
Right panel: Error-effort plot for the classical parametric harmonic oscillator,  
propagating a solution $q(t)$ of Eq.~\eqref{Mathieu} for $(\omega_0/\Omega)^2=2$, $\xi/\Omega^2=1$ over $0 \le t \le T=20\pi/\omega$, i.e over $10$ periods of the frequency modulation.
The insets display $q(t)$, with initial condition $q(0)=3$, $\dot{q}=0$,
for the same $\xi/\Omega^2=1$ and $(\omega_0/\Omega)^2=2$ (stable regime, left inset) or
$(\omega_0/\Omega)^2=1$ (unstable regime, right inset).
}
\label{fig:StabilityChart}
\end{figure}

The solution of the Mathieu equation provides us with the classical propagator $U(t,0)$, which is a $2\times 2$ matrix. According to Floquet theory, the eigenvalues $\lambda$ of $U(2\pi/\Omega)$, the propagator over one period, determine the stability of the classical system: It is stable, i.e the solutions of Eq.~\eqref{Mathieu} are bounded, if all $|\lambda| \le 1$, and unstable otherwise. 
The left panel of Fig.~\ref{fig:StabilityChart} shows the stability chart of the parametric oscillator, which we obtained with the CFET CF6:5Opt.

In the right panel of Fig.~\ref{fig:StabilityChart} we show the corresponding error-effort plot for CFETs of different order, where the error $\epsilon=\max |q(t)-\tilde{q}(t)|$ is measured as the difference between the exact  and numerical position $q(t)$ and $\tilde{q}(t)$ over $10$ periods $0 \le t \le T=20\pi/\Omega$.
The optimized 6th-order CFET  CF6:5Opt is advantageous 
for practical accuracy demands.
In the left panel of Fig.~\ref{fig:KrylovPara} we compare different CFETs 
over a range of $\xi, \omega_0$ values.  Shown is the effort needed to achieve an error of $10^{-6}$ or better. As expected, the CFET CF6:5Opt is the most efficient.

\subsection{Quantum oscillator}

For the quantum oscillator, 
we represent position $\hat{q} = (b+b^\dagger)/\!\sqrt{2\omega_0}$ and momentum operator  $\hat{p}= \ii \sqrt{\omega_0/2} (b^\dagger-b)$ through bosonic ladder operators $[b,b^\dagger]=1$.
The Hamilton operator is given by
\begin{equation}
  H(t) = \frac{\omega_0}{4} \left[ \Big( \frac{\omega(t)^2}{\omega_0^2} -1  \Big) ({b^\dagger}^2 + b^2)
    +  \Big( \frac{\omega(t)^2}{\omega_0^2} +1 \Big)  (2 b^\dagger b + 1) \right] \;.
\end{equation}
For $\xi = 0$, with $\omega(t) \equiv \omega_0$, we recover the standard Hamilton operator $H = \omega_0 (b^\dagger b + 1/2)$.
Truncation of the infinite-dimensional bosonic Hilbert space, excluding high energy states, is required to obtain the Hamilton operator as a matrix. 
For the examples we keep the lowest $N_b=50$ Fock states $|n\rangle$,
with $b^\dagger b |n\rangle = n|n\rangle$.

\begin{figure}
\begin{center}
\includegraphics[width=0.45\textwidth]{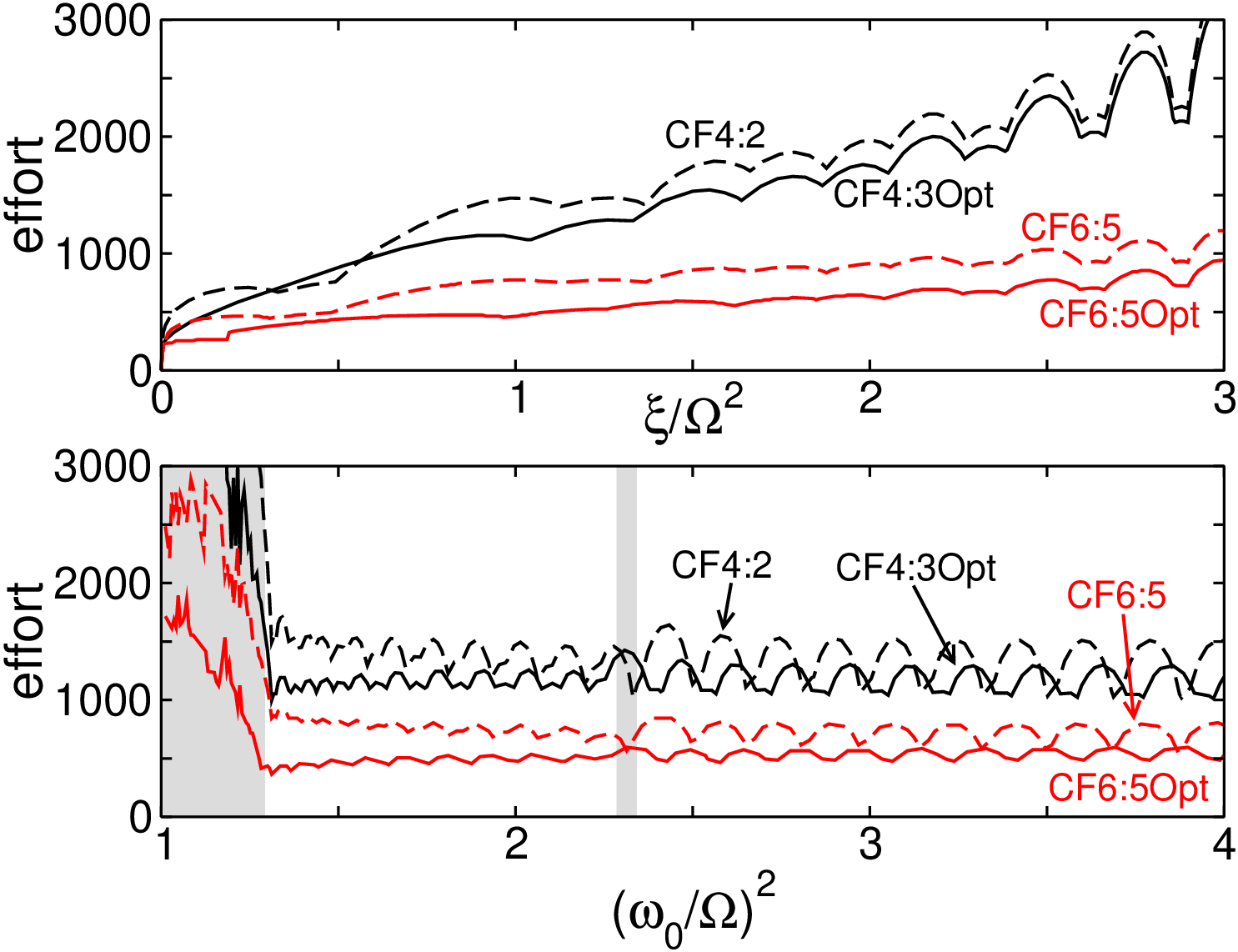} 
\hfill
\includegraphics[width=0.45\textwidth]{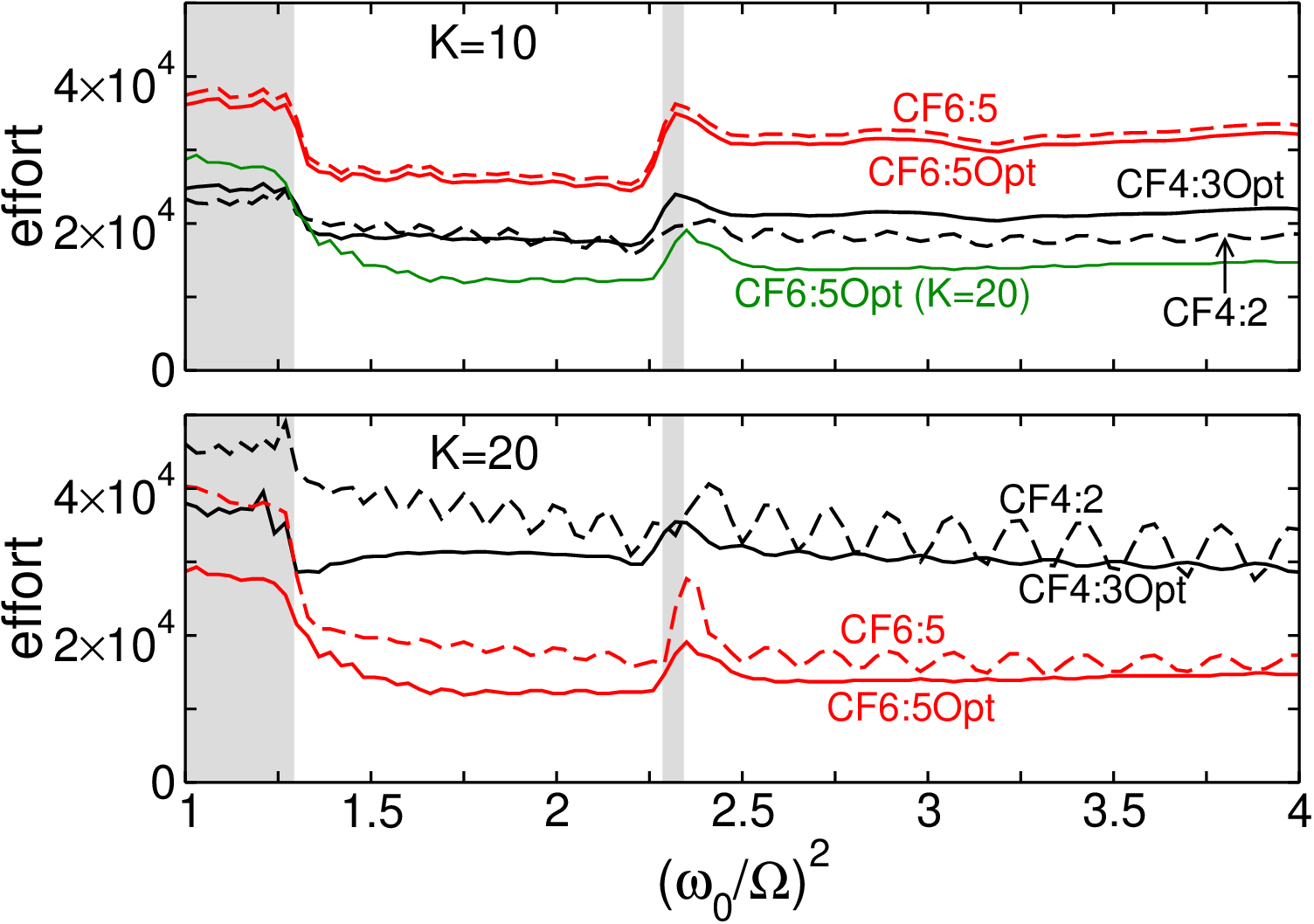} 
\end{center}
\caption{Effort required to achieve an error $10^{-6}$ for time-propagation of the classical (left panels) and quantum (right panels) parametric harmonic oscillator Eq.~\eqref{Hpara} over
$ 0 \le t \le T=20\pi/\omega$, i.e. $10$ periods of the frequency modulation.
Solid (dashed) curves corresponds to optimized (unoptimized) CFETs as indicated.
For values in the gray shaded regions, unstable solutions of Eq.~\eqref{Mathieu} exist according to Fig.~\ref{fig:StabilityChart}. 
Left panels: For the classical oscillator, the effort is given as the number of evaluated exponentials. 
Shown are results for $(\omega_0/\Omega)^2=2$ as a function of  $\xi/\Omega^2$ (upper panel), and for $\xi/\Omega^2=1$ as a function of $(\omega_0/\Omega)^2$ (lower panel). 
Right panels: For the quantum oscillator with $\xi/\Omega^2=1$, the initial coherent state $|\psi_c\rangle$ is propagated using $N_b=50$ Fock states. The effort counts the number of applications of the Hamilton operator.
The exponentials are evaluated with the Krylov technique, for $K=10$ (upper panel) and $K=20$ (lower panel) Krylov vectors. The green curve in the upper panel reproduces the result for CF6:5Opt, $K=20$ from the lower panel, as indicated.
}
\label{fig:KrylovPara}
\end{figure}

In the right panel of Fig.~\ref{fig:KrylovPara} we compare different CFETs for the quantum oscillator.
As the initial wave function we choose a coherent state
$|\psi_c\rangle$
with 
$\langle\psi_c|\hat{q}|\psi_c\rangle = 3$, $\langle\psi_c|\hat{p}|\psi_c\rangle = 0$.
The error is measured by the deviation $\epsilon=\max |\psi(t)-\tilde{\psi}(t)|$ between the exact and numerical wave function $\psi(t)$ and $\tilde{\psi}(t)$, over $10$ periods $0 \le t  \le 20\pi/\Omega$.
Shown is the effort needed to achieve an error of $10^{-6}$ or better, as for the classical case.

A new aspect in comparison with the CFET error analysis for the classical oscillator
is the numerical evaluation of the exponentials $e^{\Omega_i}$ with the Krylov technique. 
For few Krylov vectors ($K=10$, upper panel) the Krylov error from the approximate exponentials dominates. In this case, short time-steps are preferential to reduce the Krylov error sufficiently, with the consequence that the unoptimized 4th-order CFET CF4:2 is most efficient since it uses the smallest number of exponentials.
With more Krylov vectors ($K=20$, lower panel) the exponentials are evaluated to much higher accuracy also for longer time steps, and the expected advantage of optimized higher-order CFETs is recovered.
The overall most efficient propagation is achieved with the CFET CF6:5Opt for $K=20$.

Notice that the necessary Hilbert space truncation limits calculations in the unstable regimes shown in Fig.~\ref{fig:StabilityChart}, as the classical instability manifest itself for the quantum system in the excitation of high energy Fock states. Although the truncated Hamilton operator remains hermitian and can be used for time-propagation, the position expectation value $\langle \hat{q} \rangle (t)$ cannot be expected to obey the classical equation Eq.~\eqref{Mathieu}.

\subsection{The interaction picture for numerical time-propagation}\label{sec:IntPic}

Standard time-dependent perturbation theory is based on the interaction picture. 
Consider a decomposition $A(t) = D + B(t)$, where $D$ is a constant diagonal matrix.
The propagator for $D$ is the exponential $e^{t D}$.
The interaction picture is defined by $x^I(t) = e^{-t D} x(t)$.
If $B(t) \equiv 0$, $x^I(t)$ is constant.
Otherwise, it obeys the equation of motion 
\begin{equation}
\partial_t x^I(t) = e^{-tD} B(t) e^{tD} x^I(t) = B^I(t) x^I(t) \;,
\end{equation}
where $B^I(t) = e^{-tD} B(t) e^{tD}$.
Since $D$ is diagonal, the matrix elements of $B^I(t)$ are easily calculated, 
with $B^I_{mn} = e^{(D_{nn}-D_{mm})t} B^I_{mn}$.
Notice that the diagonal elements of $B(t)$ do not change, and a sparsity pattern is preserved.

The interaction picture is useful if it simplifies the equation of motion when $B(t)$ is a small perturbation. 
That is is generally not the case can be understood for the driven two-level system Eq.~\eqref{TwoLevelH} from Sec.~\ref{sec:2level}, where $D$ is identified with the term $ \Delta \sigma_z$.
The equation of motion in the interaction picture is identical to the original equation of motion with new parameters $\Delta^I=0$, $\omega^I = \omega-\Delta$.
As can be seen from Eq.~\eqref{TwoLevelU}, the propagator in the interaction picture is identical to the original propagator apart from an additional rotating phase $e^{\pm \ii \Delta t}$.
This implies that the calculation in the interaction picture has not simplified.
From the perspective of numerical time propagation the difficulty even increases
 since $B^I(t)$ varies faster than $B(t)$ due to the additional time-dependence acquired in the transformation with $e^{tD}$.
This is particularly true if $B(t)$ is a small perturbation, since $|\omega - \Delta| > \omega$ for large $\Delta$.

Notice that for the present problem the choice $D= \pm \omega \sigma_z$ leads to a constant Hamilton operator in the interaction picture, which allows for the construction of the exact propagator Eq.~\eqref{TwoLevelU}.
Indeed, the celebrated rotating wave approximation is exact for this particular case.
Despite its persistence
in quantum optics it does not easily generalize to other situations.

\subsection{The interaction picture for the harmonic oscillator}

\begin{figure}
\begin{center}
\includegraphics[width=0.45\textwidth]{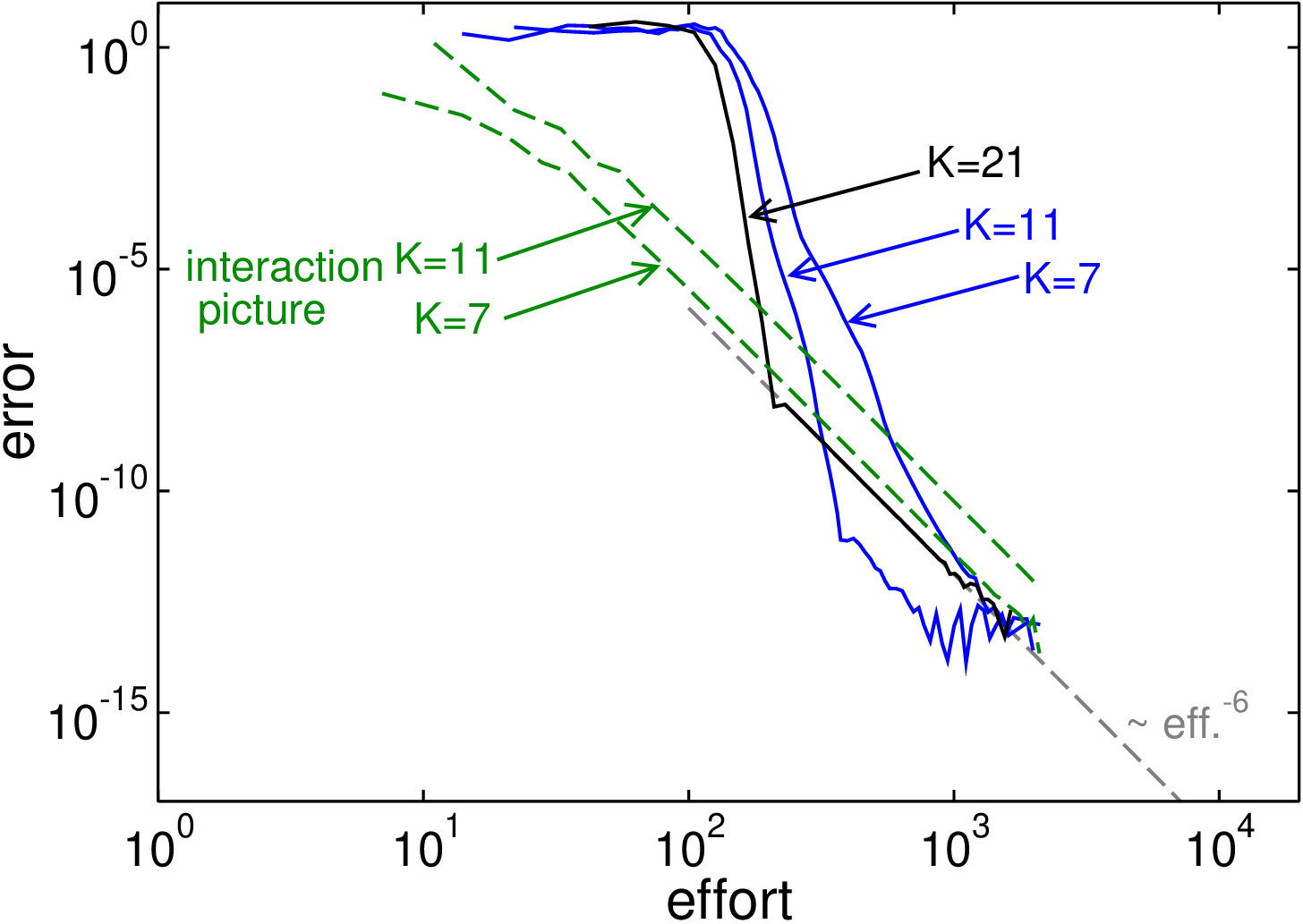} 
\hfill
\includegraphics[width=0.45\textwidth]{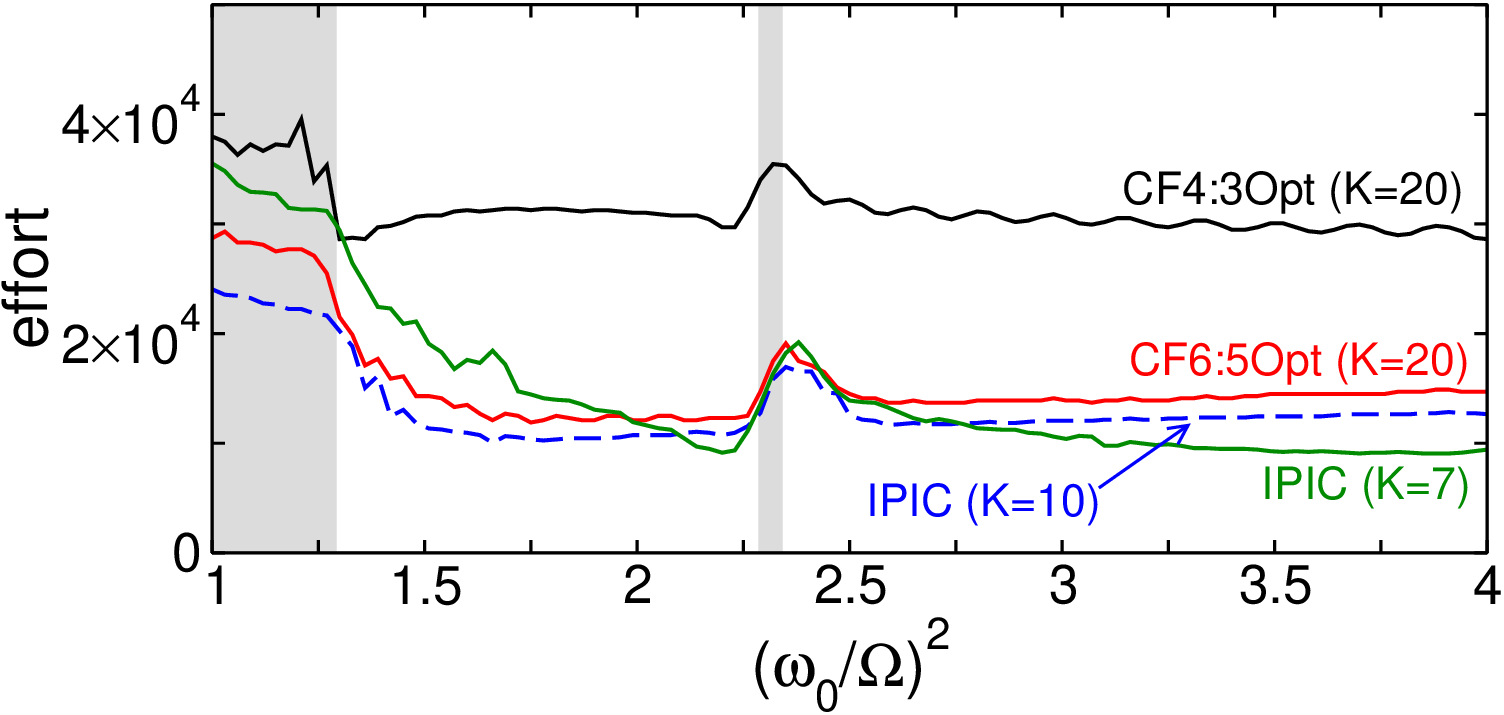} 
\end{center}
\caption{Time-propagation in the interaction picture (IPIC) of the quantum parametric harmonic oscillator Eq.~\eqref{Hpara} with $\xi/\Omega^2=1$ over $0 \le t \le T=20\pi/\omega$.
Left panel: Error-effort plot for $(\omega/\Omega)^2=4$,  using the CFET CF6:5Opt and the Krylov technique with $K$ Krylov vectors as indicated. The solid curves give the error from standard propagation, the dashed curves from working in the interaction picture.
The gray dashed line shows the asymptotic scaling of the CFET error $\sim \mathsf{effort}^{-6}$.
Right panel: Effort required to achieve an error $10^{-6}$ with the CFET CF6:5Opt,
as a function of $(\omega_0/\Omega)^2$.
For comparison, results for standard propagation with $K=20$ of Krylov vectors are reproduced from the lower right panel in Fig.~\ref{fig:KrylovPara}.
}
\label{fig:KrylovInter}
\end{figure}

While the interaction picture per se does not simplify time-propagation, it can be useful to reduce the computational effort associated with the numerical evaluation of exponentials.
As discussed in Sec.~\ref{sec:Krylov2}, the quantum harmonic oscillator is an example where large eigenvalues $n \omega$ increase the effort.
Switching to the interaction picture, with $D=\omega b^\dagger b$, increases the CFET error 
because of the additional time-dependence on the scale of  $\omega$, but simplifies the evaluation of exponentials since the large diagonal entries $n\omega$ are removed from the matrix.
We illustrate this possibility with the error-effort plot for the CFET CF6:5Opt in Fig.~\ref{fig:KrylovInter} (left panel), where results for the interaction picture are compared to those from standard propagation for a different number $K$ of Krylov vectors.
We see that in the interaction picture the Krylov error is much reduced so that the CFET error, with scaling $\sim \textsf{effort}^6$, dominates over the entire range. 
For moderate accuracy demands, with errors down to $10^{-7}$, the interaction picture with only $K=7$ Krylov vectors is most efficient.
For smaller error, the interaction picture is again less favourable, since the CFET error has increased in comparison to standard propagation.
The right panel of Fig.~\ref{fig:KrylovInter} shows the effort to achieve an error $10^{-6}$,
supporting the expectation that the interaction picture becomes rather efficient at larger $(\omega_0/\Omega)^2$.
Whether there is a benefit of using the interaction picture also for non-bosonic systems remains to be studied.

\section{Comparison of CFETs to Floquet approaches}\label{sec:Floquet}

For problems with a periodic time-dependence Floquet theory suggests exploitation
of the periodicity of the propagator. A notable implementation of this idea is the $(t,t')$-method~\cite{PKM94}. Introducing time as an additional variable $t'$, 
the wave function $\psi(t)$ is recovered from the solution $\Psi(t,t') = \exp(-\ii \mathcal{H} t) \Psi(0,t')$
of the Schr\"odinger equation with a time-independent Hamilton operator $\mathcal{H} = H(t')-\ii \partial_{t'}$ as $\psi(t) = \Psi(t,t)$. The validity of this procedure can be checked by evaluation of $\partial_t \Psi(t,t)$, with initial condition $\Psi(0,t') = \psi(0)$. 
In computations, the auxiliary degree of freedom $t'$ is represented with a Fourier basis of periodic functions $\phi_n(t')=e^{2 \pi \ii n t'/T}$.
The calculation of the matrix exponential $ \exp(-\ii \mathcal{H} t)$ in the enlarged Hilbert space is ideally suited for the Chebyshev technique providing solutions for one or more periods at once. The accuracy is determined by the number $N_F$ of Fourier modes kept in the calculation.

\begin{figure}
\begin{center}
\includegraphics[width=0.45\textwidth]{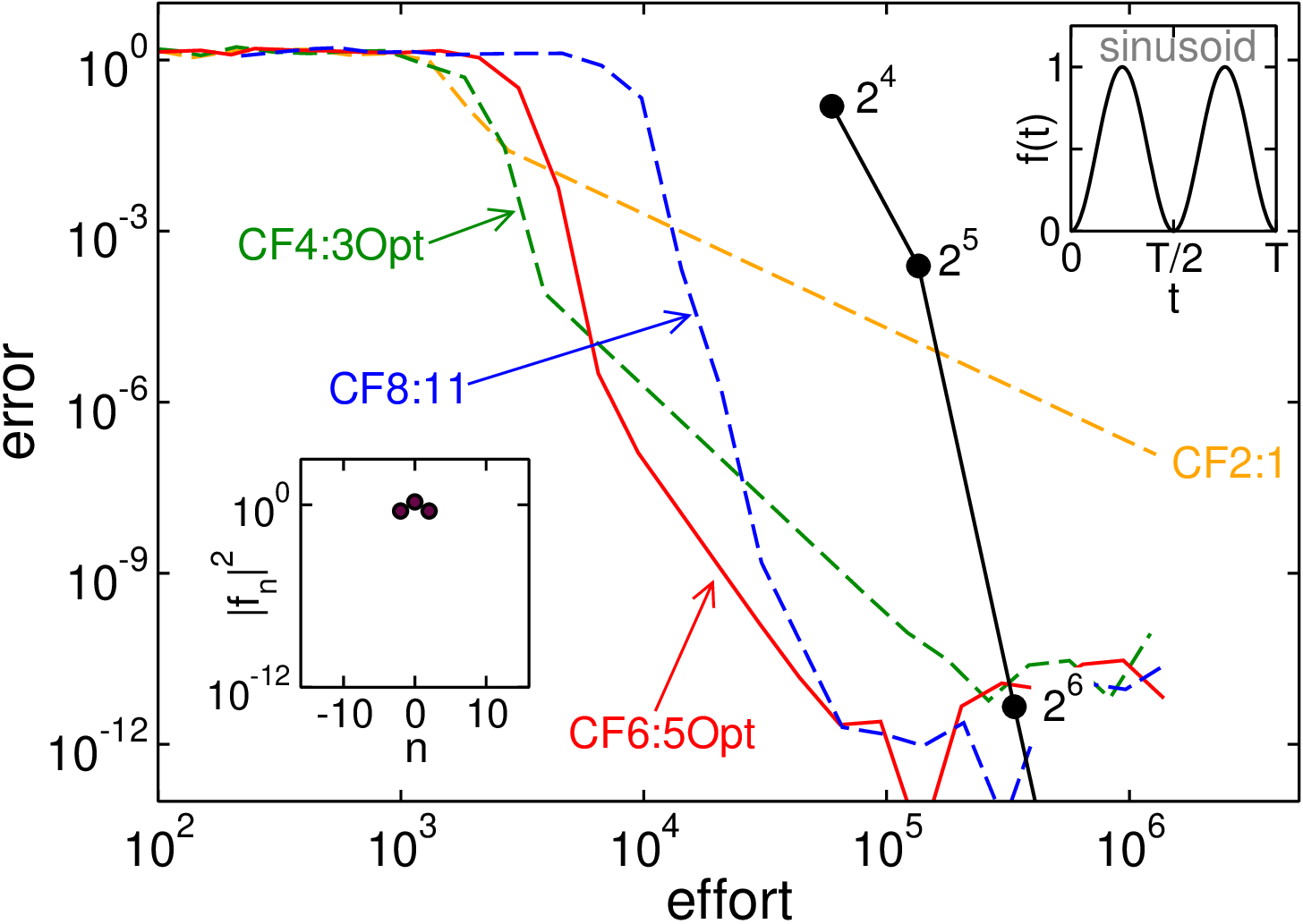} 
\hfill
\includegraphics[width=0.45\textwidth]{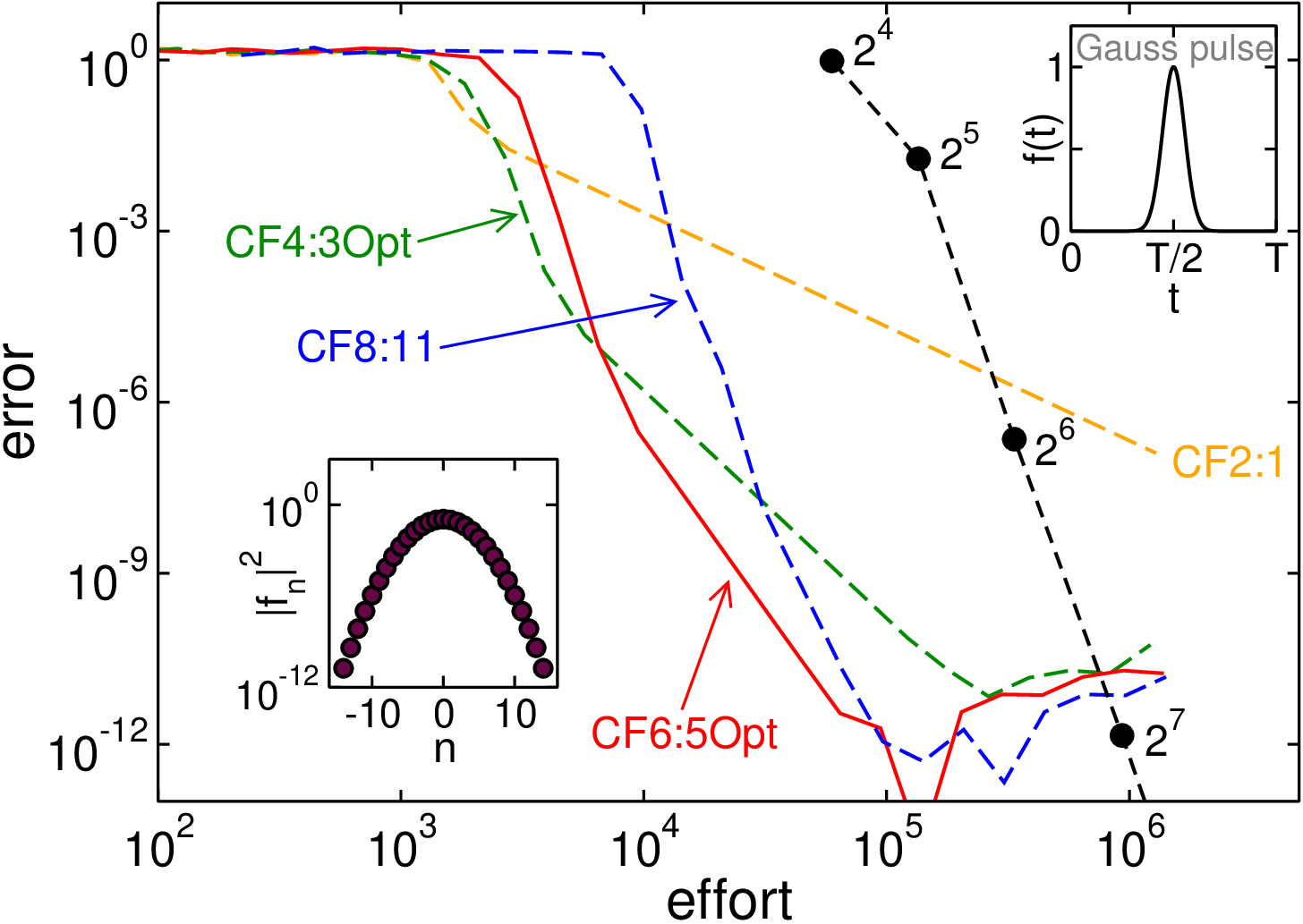} 
\end{center}
\caption{Comparison of the $(t,t')$-method with CFETs for the forced quantum harmonic oscillator.
Similar to Figs.~\ref{fig:KrylovPara},~\ref{fig:KrylovInter}, the initial coherent state $|\psi_c\rangle$ 
is propagated in a truncated Hilbert space with $N_b=50$ Fock states over $10$ periods 
of  the sinusoidal force (left panel) and the Gaussian pulse (right panel) specified in the text.
The black circles correspond to the $(t,t')$-method with $N_F$ Fourier modes as indicated,
the other curves to CFET/Krylov propagation with $K=10$ Krylov vectors. 
The insets display the Fourier coefficients $|f_n|^2$ and $f(t)$ itself.
}
\label{fig:ttprime}
\end{figure}

In Ref.~\cite{PKM94}, the $(t,t')$-method was compared to a 2nd-order Magnus propagator.
It was found that the $(t,t')$-method is far more efficient and allows for reduction of the error down to machine precision with moderate effort.
Following these examinations, we consider the quantum harmonic oscillator $H = \hat{p}^2/2 + \hat{q}^2/2 + f(t) \hat{q}$ with a time-dependent periodic force $f(t)=f(t+T)$.
We propagate the initial coherent state $|\psi_c\rangle$ over $10$ periods ($T=5/3 \pi$) with (i) a sinusoidal force $f(t) = \sin^2  2 \pi t /T$, (ii) a Gaussian pulse $f(t) = \exp(-((t-T/2)/0.4)^2)$,
The error-effort plot in Fig.~\ref{fig:ttprime} compares the $(t,t')$-method with higher-order CFETs.

We see that in both examples the $(t,t')$-method is significantly less efficient than any but the 2nd-order CFET. 
Although the $t$-$t'$ error drops rapidly once $N_F$ is sufficiently large to represent the Fourier components of the auxiliary wave function $\Psi(t,t')$,
even moderate accuracy requires  $N_F \ge 2^6$ and proportionately large effort.
For the Gaussian pulse more Fourier modes must be kept, since weight is distributed to higher Fourier coefficients $f_n =(1/T) \int_0^T f(t) e^{2\pi\ii n t/T} dt$ of the driving force $f(t)$. This restricts the use of the $(t,t')$-method if memory limitations are a concern. Notice that splitting the periodic problem into several time-steps increases the effort further, in particular since the Fourier coefficients of the then discontinuous force decay more slowly.

The poor efficiency of the $(t,t')$-method in comparison to the higher-order CFETs is not a failure of the Floquet approach.
If we associate a fictitious time-step $T/N_F$ with the representation of the wave function $\Psi(t,t')$ through $N_F$ Fourier nodes per period,
it is much larger than the time-step in the CFET time-stepping.
This is in accordance with the expectation that for periodic problems Fourier decomposition provides a better representation of the propagator than the concatenation of step-wise constant propagators.
A related observation is the increased accuracy of the Fourier transform for integration of periodic functions over the combination of finite order polynomial integration formulae.
In total, the $(t,t')$-method requires less application of the Hamilton operator $\mathcal{H}$ for propagation over the entire $10$ periods than the CFET/Krylov technique with short time-steps.
However, the practically relevant effort of computations in the Fourier space is just larger by $N_F$, which, effectively, renders the $(t,t')$-method less efficient than higher-order CFETs.

\section{Further applications}
\label{sec:App}

We complete our study of the practical applicability of CFETs with calculations for two complex quantum systems, for which neither exact solutions nor classical analogues are known: 
A chain of interacting spins -- or two-level atoms -- in pulsed magnetic fields (Sec.~\ref{sec:MultSpin}),
and the hydrogen atom in an electric field (Sec.~\ref{sec:Hydro}).
Both systems feature non-trivial physical effects,
and require computation of exponentials for moderate-to-large sparse matrices.

\subsection{Driven spin chain}\label{sec:MultSpin}

In first approximation atoms in a strong light field can be described by interacting spins $1/2$ in a magnetic field. We consider the Hamilton operator
\begin{equation}\label{SpinChain1}
  H = \sum_{s=1}^S H^{(s)} + J \sum_{s=1}^{S-1}  (\sigma^{(s)}_x \sigma^{(s+1)}_x + \sigma^{(s)}_y \sigma^{(s+1)}_y )
\end{equation}
of a spin chain with $S$ spins, where
\begin{equation}\label{SpinChain2}
 H^{(s)} = \Delta \sigma_z^{(s)} + \Re V(t) \, \sigma^{(s)}_x - \Im V(t) \, \sigma^{(s)}_y 
 = \begin{pmatrix} \Delta & V(t) \\ V^*(t) & -\Delta \end{pmatrix}
\end{equation}
is the Hamilton operator of a single spin, subjected to a magnetic field similar to Eq.~\eqref{TwoLevelH}.
If the system is initially prepared in the ground state, the magnetic field induces transitions to excited states.
For the choice 
\begin{equation}\label{BPulse}
  V(t) = \frac{V e^{-2 \ii \omega t}}{\cosh t/\tau} \;,
\end{equation}
a magnetic pulse of half-width $\approx 1.32 \tau$ and frequency $\omega$,
the transition probability for a single spin ($S=1$) can be deduced from the result for the Rosen-Zener model~\cite{RZ32}.
Specifically, the transition probability $P_\infty=|\langle{\uparrow}|U(\infty,-\infty)|{\downarrow}\rangle|^2$, i.e. the probability that the spin is flipped through the pulse, is 
\begin{equation}
  P_\infty = \frac{\sin^2 \pi V \tau}{\cosh^2 \pi (\Delta-\omega) \tau} \;.
\end{equation}

\begin{figure}
\begin{center}
\includegraphics[width=0.45\textwidth]{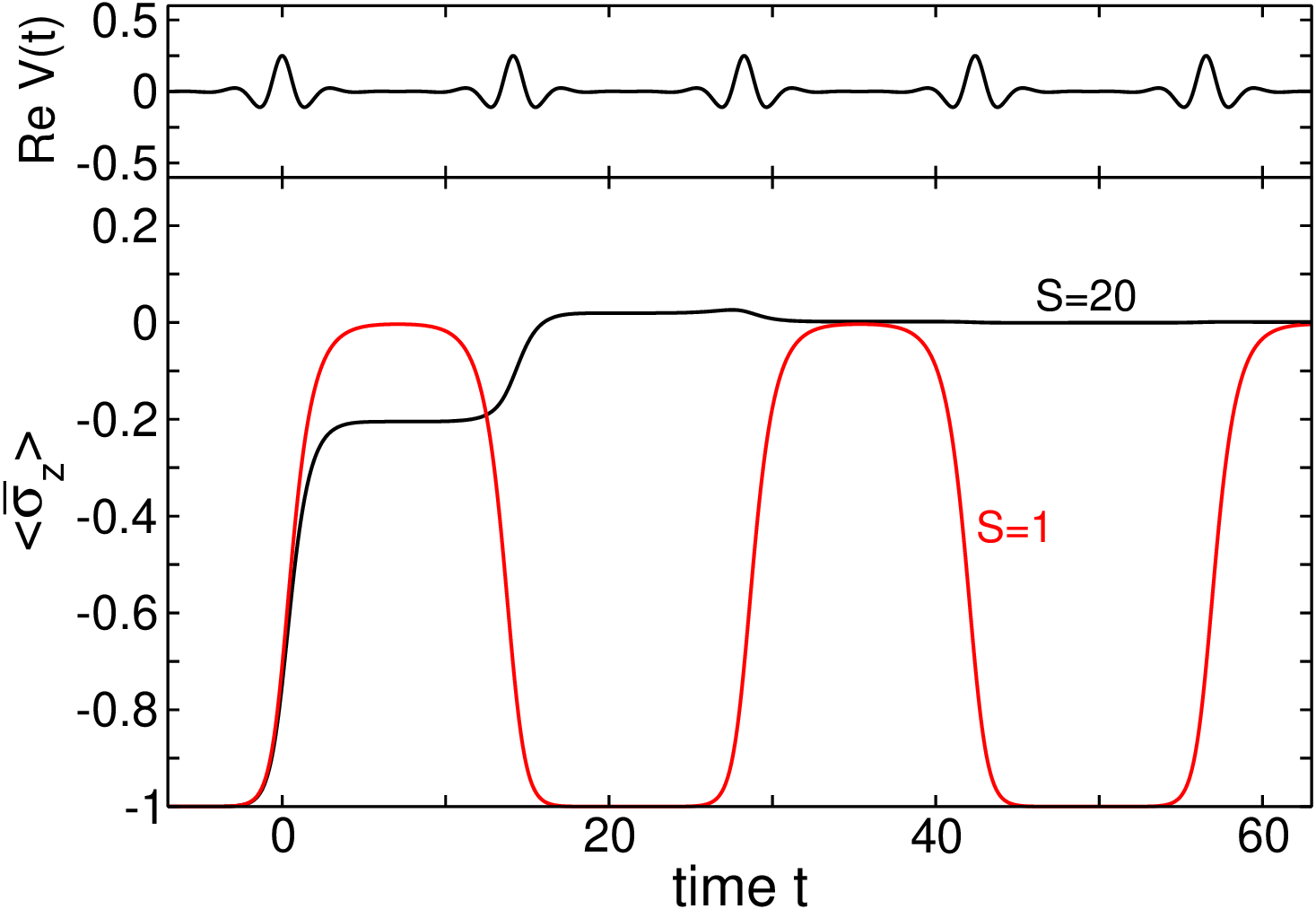} 
\hfill
\includegraphics[width=0.45\textwidth]{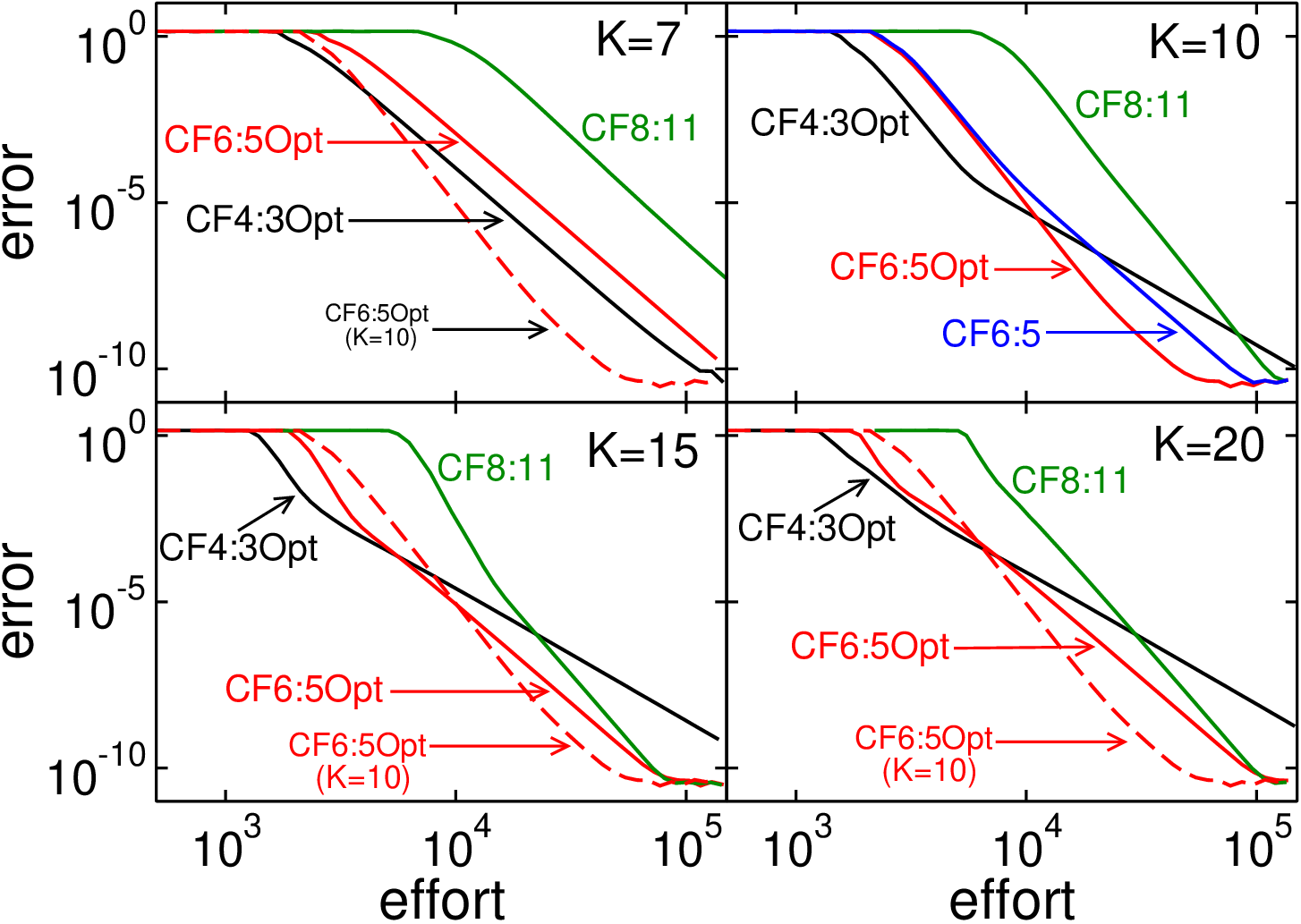} 
\end{center}
\caption{Time-propagation of the driven spin chain Eqs.~\eqref{SpinChain1},~\eqref{SpinChain2},
with $\Delta=1$, $J=0.1$.
The system is pumped by a sequence of resonant pulses (Eq.~\eqref{BPulse} with $\omega=\Delta$, $\tau=1$, $V=1/(4\tau)$), centered at multiples of $t_0=9\pi/2$.
The initial state at $-t_0/2$ is the product state $|{\downarrow \dots \downarrow}\rangle$.
Left panels: Total spin $z$-component $\bar{\sigma}_z(t)$ as a function of time, for $S=20$ spins. For comparison, the curve for a single spin ($S=1$) is included.
The upper panel shows the pulsed field $\Re V(t)$.
Right panel: Error-effort plot for the time-propagation shown in the left panel, for different CFETs and number of Krylov vectors $K$ as indicated.
For comparison, the panels include either the curve for CF6:5Opt, $K=10$, or for CF6:5 in the upper right panel.}
\label{fig:ManySpins}
\end{figure}

In Fig.~\ref{fig:ManySpins} (left panel) we show the expectation value $\bar{\sigma}_z(t) \equiv (1/S) \sum_{s=1}^S \langle \psi(t)|\sigma_z^{(i)}|\psi(t) \rangle$ for a sequence of magnetic field pulses. The pulse sequence brings a single spin (curve for $S=1$) back to its initial state after two subsequent pulses.
For several interacting spins ($S=20$), dephasing leads to a state with $\bar{\sigma}_z(t) \equiv 0$ after the first few pulses. 

The right panel in Fig.~\ref{fig:ManySpins} compares the efficiency of different CFETs with a different number $K$ of Krylov vectors.
This example shows, similar as for the harmonic oscillator, the importance of balancing the Krylov and CFET error.
For small $K=7$ the Krylov error dominates, which gives the 4th-order CFET CF4:3Opt an advantage over higher-order CFETs because it requires less exponentials per time-step.
The Krylov error is however less dominant than for the harmonic oscillator, and the 6th-order CFET CF6:5Opt with $K=10$ results in the most efficient time-propagation.
Notice that the unoptimized CFET CF6:5 (upper right panel) is about $50\%$ less efficient.
As an interesting feature we note that the slope of the curves for $K=10$ resembles that of a 9th-order relation ($\mathsf{error} \sim 1/\mathsf{effort}^9$), which is the expected scaling of the Krylov error for $K=10$ (cf. Eq.~\eqref{KryErr}). 
The `bend' from the 9th-order scaling to a 4th-order scaling is clearly seen in the curve for CF4:3Opt. The error of higher-order CFETs remains smaller than the Krylov error, and 9th-order scaling persists down to machine precision. 

\subsection{The hydrogen atom in an electric field}\label{sec:Hydro}

Our last example is that of a hydrogen-like atom in a classical monochromatic electric field along the $z$-axis. The Hamilton operator in dipole approximation is $H= -\nabla^2 -\frac{2}{r} + E_z(t) d_z$, where $E_z(t)$ denotes the field strength and $d_z \equiv z$ is the $z$-component of the dipole operator.
Working in the basis of hydrogen eigenstates $|nlm\rangle$,
with energy $\omega_n = - 1/n^2$ for $E_z(t) \equiv 0$,
the quantum number $m$ is conserved for the above Hamiltonian. We consider only the $m=0$ sector.
 The required matrix elements of the dipole operator $d_z$ can be calculated analytically or with a one-dimensional numerical integration.
They are non-zero only between states for which the respective $l$ differs by $\pm 1$.

\begin{figure}
\begin{center}
\includegraphics[width=0.45\textwidth]{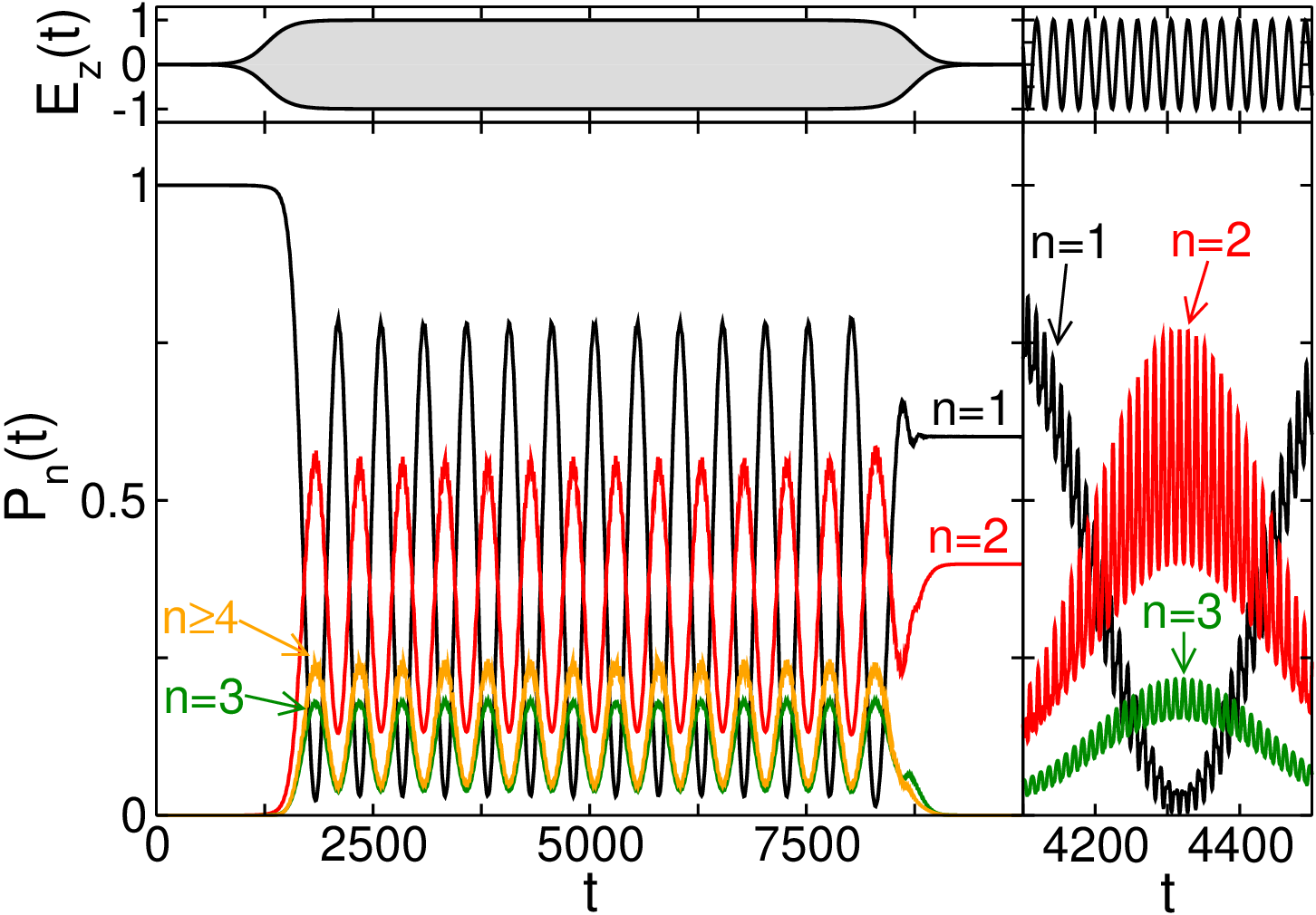} 
\hfill
\includegraphics[width=0.45\textwidth]{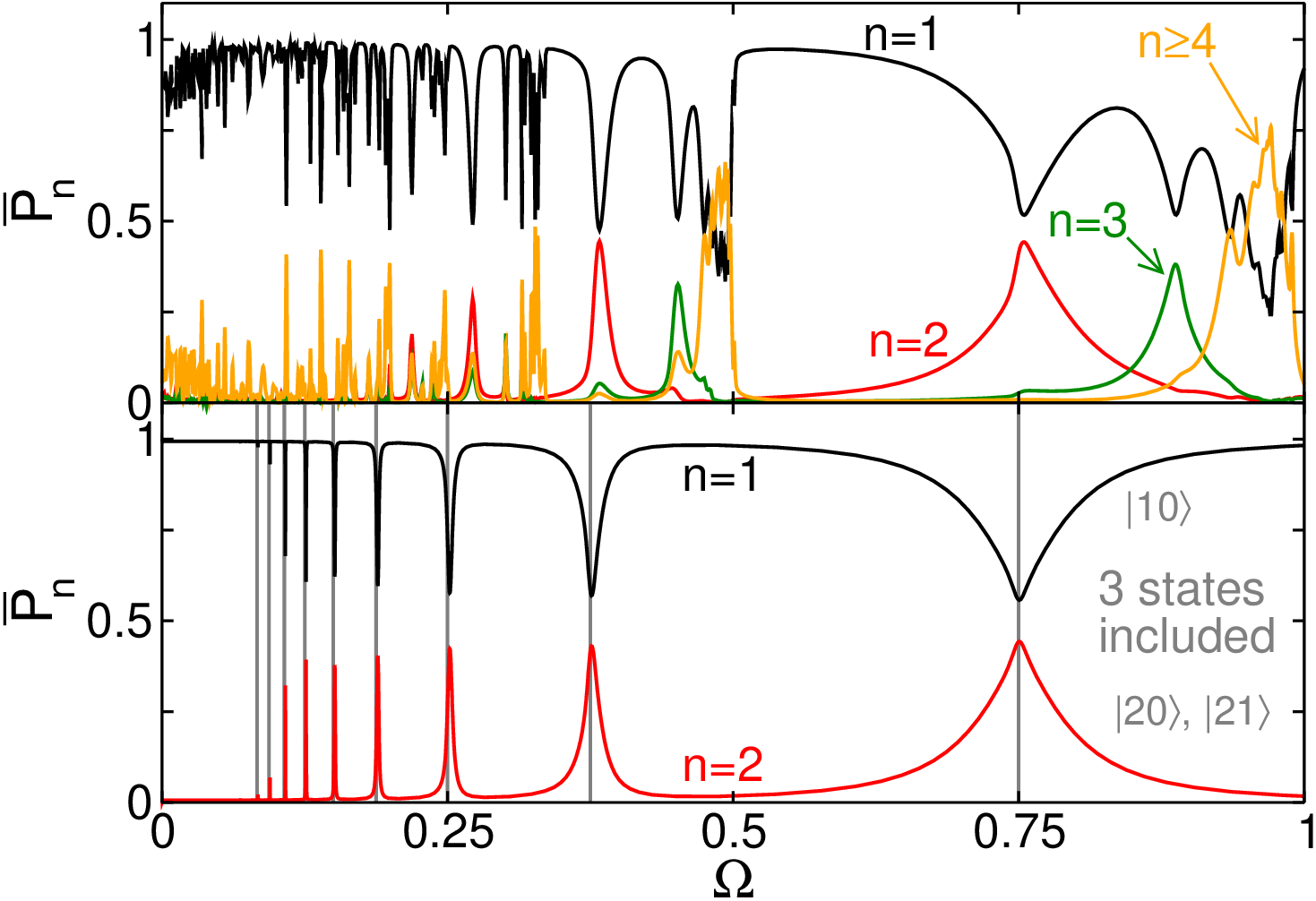} 
\end{center}
\caption{Time evolution of a hydrogen atom in an alternating electric field.
Left panel: Occupation probability $P_n(t)$ as a function of time,
for $\Omega=0.27$, $E_z^0=0.1$. The upper row gives $E_z(t)$. 
On the left we show the envelopes, while the magnification on the right resolves the fast oscillations (the curve for $n \ge 4$ is omitted). 
The shown propagation time covers $430$ field oscillations.
All $n \le 50$ states are kept in the calculation.
Right panel: Time-averaged occupation probability $\bar{P}_n$ as a function of $\Omega$,
for $E_z^0=0.1$.
The lower panel shows the result for the three level system $|10\rangle$, $|20\rangle$, $|21\rangle$.
The vertical gray lines indicate the resonance frequencies $3/(4k)$, $k=1,2,\dots$,
for transitions $|10\rangle \leftrightarrow |21\rangle$.
}
\label{fig:hydro}
\end{figure}

The system is initially prepared in the $E_z(t) \equiv 0$ ground state $\psi(t=0) = |10\rangle$,
and $\psi(t)$ is calculated for $0 \le t \le T=10^4$ using the CFET CF6:5Opt in combination with the Krylov technique ($K=10$). The electric field 
is given by $E_z(t) = E_z^0 h(t) \cos \Omega t$,
where $h(t)=(1+a)/(1+a \exp(-b (x-t_0)^2))$ is an envelope function with $a=b=10^{-6}$, $t_0=5000$.
In Fig.~\ref{fig:hydro} we show the summed occupation probability $P_n(t) = \sum_{l=0}^{n-1} |\langle n l|\psi(t)\rangle|^2$ (left panel) and its time average $\bar{P}_n =  (1/T) \int_0^T P_n(t) dt$ (right panel).
In the weak coupling limit $E_z \to 0$, 
resonances occur if the transition frequency $|\omega_{n_1}-\omega_{n_2}|$
between states $|n_1,l\rangle$ and $|n_2,l \pm 1\rangle$ is a multiple of the field frequency $\Omega$.
This behaviour is clearly seen if only the three $n=1,2$ states $|10\rangle$, $|20\rangle$, $|21\rangle$ are included in the calculation (lower right panel in Fig.~\ref{fig:hydro}).
The broad resonance at $\Omega = \omega_1-\omega_2 = 3/4$ is most pronounced, 
while the resonances at $\Omega = 3/8, 3/12, \dots$ become increasingly sharp (for a non-classical field, these would correspond to multi-photon absorption).
Inclusion of states with larger $n$ (upper right panel, with $n \le 50$ in the numerical calculation) shifts the frequencies of the $n=1 \leftrightarrow n=2$ transition, and leads to the numerous sharp resonances of transitions to higher excited states.

\section{Conclusions}
\label{sec:Conc}

The development of practicable techniques for the propagation of driven quantum systems 
requires realization of high theoretical efficiency gains under the restrictions of actual applications. 
In the present paper we studied a particular class of numerical techniques,
the commutator-free exponential time-propagators, which combine favourable theoretical properties, such as preservation of unitarity and high approximation order, with the virtue of simple implementation.

Conceptually, CFETs are related to the more traditional Magnus expansion.
From the practical point of view, they are in fact the better alternative, at least for the problems studied here. Avoiding commutators makes them easier to implement and also more efficient, since the complicated structure of the original Magnus expansion and all the bookkeeping it requires is replaced by their simple exponential product form.

We dealt with the derivation, optimization, and application of CFETs from the common point of view of the practitioner who wants to solve the Schr\"odinger equation.
For every issue the present work extends the existing literature.
Our construction and analysis of CFETs relies essentially on 
the use of Legendre polynomials and their orthogonality properties.
In this way we can provide a comprehensive and self-contained presentation. It also simplifies the error analysis and allows us to identify the importance of including higher-order terms for the CFET optimization.
We provide coefficients of fully optimized 4th- and 6th-order CFETs, as well as of a good albeit unoptimized 8th-order CFET.
As both the theoretical and practical error analysis show full optimization is successful in further reducing the error, leading to about 50\% higher efficiency in comparison to the partly optimized counterparts.
While the potential of 6th-order CFETs is probably largely exhausted,
optimization of 8th-order CFETs remains promising.

We have discussed the practical application of CFETs at great length, paying particular attention to 
realistic situations where exponentials can not be calculated in closed form. 
Based on our findings, we generally recommend the use of the CFET CF6:5Opt together with a Krylov calculation of the exponential using about 10-15 Krylov vectors.
The results for the examples presented show that very accurate results can be obtained with moderate effort. 
They provide evidence that for the Schr\"odinger equation optimized higher-order CFETs are substantially more efficient than alternative techniques such as general purpose Runge-Kutta methods or numerical Floquet approaches. 
Most importantly, CFETs are robust: They are unconditionally stable, and their quality does not substantially decline 
at points of resonance.
CFETs are thus a good choice for library routines for time-propagation. 
We believe that the implementation and optimization of a general purpose time-propagation routine provides most potential for further significant efficiency gains.
Irrespective of machine dependent implementation details,
this has to include refined strategies for the automated choice of the step-size and the number of Krylov vectors, as well as tracking of the accumulated error.
Even now CFETs are a viable and convenient technique for the time-propagation of driven quantum systems.

\section{Acknowledgments}
This work was financed by Deutsche Forschungsgemeinschaft via Sonderforschungsbereich 652 and AL 1317/1-1.

\appendix

\section{Recursion for the Magnus expansion}
\label{app:RecMag}

It is possible to write every $\Omega_n(t)$ from the Magnus expansion as an $n$-fold time-ordered integral 
\begin{equation}
\Omega_n(t) = \int\limits_0^t dt_1 \int\limits_0^{t_1} dt_2 \dots \int\limits_0^{t_{n-1}} dt_n \, Z_n(t_1,\dots,t_n) = \int\limits_{\Delta_n[t|1,\dots,n]} Z_n(1,\dots,n)
\end{equation}
of a multivariate function $Z_n(1,\dots,n) \equiv Z_n(t_1,\dots,t_n)$.
With regard to Eq.~\eqref{MagnusOm}, we have
\begin{equation}\label{AppZ}
Z_1(t_1) = A(t_1) \;, \;
Z_2(t_1,t_2) = \frac{1}{2} [A(t_1),A(t_2)] \;, \; 
Z_3(t_1,t_2,t_3) = \frac{1}{6} [A(t_1),[A(t_2),A(t_3)]] + \frac{1}{6} [[A(t_1),A(t_2)],A(t_3)] \;.
\end{equation}
The integration domain of the time-ordered integral is the set of decreasing $n$-tuples
\begin{equation}
  \Delta_n[t| i_1, \dots,  i_n] = \{(t_1, \dots, t_n) \in \mathbb{R}^n |
 t > t_{i_1} > t_{i_2} > \dots > t_{i_n} > 0 \} \;,
\end{equation}
where ${i_1,\dots,i_n}$ is a permutation of ${1,\dots,n}$ denoting the arrangement of the tuple elements. For example, $\Delta_2[t|1,2] = \{(t_1,t_2) \in \mathbb{R}^2 | t > t_1 > t_2 > 0\}$ and
$\Delta_2[t|2,1] = \{(t_1,t_2) \in \mathbb{R}^2 | t > t_2 > t_1 > 0\}$.
Every permutation selects one of the $n!$ wedge-shaped subsets of the $n$-dimensional hypercube $[0,t]^n$, which is the disjoint union of all these sets
(up to points from an $n-1$-dimensional subset, which as a set of measure zero is irrelevant for integration).  
The time derivative $\dot{\Omega}_n(t)$ is given by the $n-1$-fold integral
\begin{equation}
\dot{\Omega}_n(t) = \int\limits_{\Delta_{n-1}[t|1,\dots,n-1]} Z_n(t,1,\dots,n-1) \;.
\end{equation}
We also note that
\begin{equation}\label{AppA1}
 \int\limits_0^t dt' \int\limits_{\Delta_n[t'|1,\dots,n]} f(t',1,\dots,n) = \int\limits_{\Delta_{n+1}[t|1,\dots,n+1]} f(1,\dots,n+1) \;.
\end{equation}

According to Eq.~\eqref{OmRec}, $\Omega_{n+1}$ is given as 
\begin{equation}\label{AppA2}
\begin{split}
\Omega_{n+1} (t)  &=  \sum_{m=1}^{n} \frac{(-1)^{m+1}}{(m+1)!} \sum_{\substack{n_1, \dots, n_{m+1} \ge 1 \\ n_1 + \dots + n_{m+1} = n+1}} \int\limits_0^t dt'
 [ \dots [\dot{\Omega}_{n_1}(t'),\Omega_{n_2}(t')],\dots, \Omega_{n_{m}}(t')],\Omega_{n_{m+1}}(t')] \\
  &=  \sum_{m=1}^{n} \frac{(-1)^{m+1}}{(m+1)!} \sum_{\substack{n_1, \dots, n_{m+1} \ge 1 \\ n_1 + \dots + n_{m+1} = n+1}} \int\limits_0^t dt' \int\limits_{\mathcal{I}_{n}[t'|n_1-1,n_2,\dots,n_{m+1}]} 
  [Z_{n_1},Z_{n_2},\dots,Z_{n_{m+1}}](t',1,\dots,n) \;,
  \end{split}
\end{equation}
where we use the notation ($n=n_1+\dots+n_k$)
\begin{equation}
[Z_{n_1},Z_{n_2},\dots,Z_{n_k}](1,\dots,n) = [\dots [Z_{n_1}(1,\dots,n_1),Z_{n_2}(n_1+1,\dots,n_1+n_2) ],\dots,Z_{n_k}(n-n_k+1,\dots,n)  ]  
\end{equation}
for the nested commutator in the integrand. The integration domain is a product set
\begin{equation}
\mathcal{I}_{n}[t | n_1,\dots,n_k] = \Delta_{n_1}[t|1,\dots,n_1] \times \Delta_{n_2}[t|1,\dots,n_2] \times \dots \times  \Delta_{n_k}[t| 1,\dots,n_k]    \qquad (n=n_1+\dots+n_k) \;.
\end{equation}

To bring the integrals in Eq.~\eqref{AppA2} into time-ordered form, the integration domain 
is split into disjoint pieces that are mapped onto the `wedge' sets $\Delta_n[t|\dots]$ through a permutation of the integration variables. 
For every $(t_1,\dots,t_n) \in \mathcal{I}_{n}[t | n_1,\dots,n_k]$ a unique permutation $\pi$ exists that orders the $n$-tuple such that $(t_{\pi^{-1}(1)}, \dots, t_{\pi^{-1}(n)}) \in \Delta_n[t|1,\dots,n]$.
The admissible permutations are those that respect the order of elements corresponding to each of the $\Delta_{n_i}[t|\dots]$ factors in $\mathcal{I}_{n}[t|n_1,\dots,n_k]$. 
These form the set
\begin{equation}
\begin{split}
 \mathcal{P}_n[n_1,\dots,n_k] = \{\pi \text{ is permutation of } \{1,\dots,n\} | & \, \pi(1) < \dots < \pi(n_1) \text{ and } \pi(n_1+1) < \dots < \pi (n_1+n_2) \\
 &\dots \text{ and } \pi (n-n_k+1) < \dots < \pi(n) \} \;,
\end{split}
\end{equation}
where still $n=n_1+\dots+n_k$. It has $n!/(n_1! \cdots n_k!)$ elements.
In particular, $\mathcal{P}_n[1,1,\dots,1]$ is the set of all permutations,
while $\mathcal{P}_n[n]$ contains only the identity.

The decomposition of $\mathcal{I}_{n}[t|n_1,\dots,n_k]$ into disjoint subsets congruent with $\Delta_n[t|1,\dots,n]$ is given by 
\begin{equation}
\mathcal{I}_{n}[t|n_1,\dots,n_k] = \biguplus\limits_{\pi \in \mathcal{P}_n[n_1,\dots,n_k]} \Delta_n[t|\pi^{-1}(1),\dots,\pi^{-1}(n)] \;.
\end{equation}
Permutation of the integration variables then gives the identity
\begin{equation}
 \int\limits_{\mathcal{I}_{n}[t|n_1,\dots,n_k]} f(1,\dots,n) =
\sum_{\pi \in \mathcal{P}_n[n_1,\dots,n_k]} \;\;  \int\limits_{\Delta_n[t|1,\dots,n]}  f(\pi(1),\dots,\pi(n)) \;.
\end{equation}
This identity allows us to express the integrals in Eq.~\eqref{AppA2} as time-ordered integrals. 
The final integration over $t'$ preserves time-ordering according to Eq.~\eqref{AppA1}.

After these preparations we can finally state the recursion 
\begin{equation}\label{ZRec}
\begin{split}
&Z_1(t)  = A(t) \,, \\
 &Z_{n+1}(0,1,\dots,n) =  \sum_{m=1}^n \frac{(-1)^{m+1}}{(m+1)!} \sum_{\substack{n_1, \dots, n_{m+1} \ge 1 \\ n_1 + \dots + n_{m+1} = n+1}}  \sum_{\pi \in \mathcal{P}_n[n_1-1,\dots,n_{m+1}]} 
 [Z_{n_1},\dots,Z_{n_{m+1}}](0,\pi(1),\dots,\pi(n)) \;.
\end{split}
\end{equation}
While the first terms $Z_n$ can be obtained by hand, 
the calculation of higher terms 
is better left to the computer.
Consider exemplarily the calculation of $Z_3$. 
The sum over $n_1,\dots,n_{m+1}$ contains $2+1$ terms for $m=1,2$.
Thus,
\begin{equation}
\begin{split}
 Z_3(0,1,2) &= \frac{1}{2} \sum_{\pi \in \mathcal{P}_2[0,2]} [Z_1,Z_2](0,\pi(1),\pi(2))
  + \frac{1}{2} \sum_{\pi \in \mathcal{P}_2[1,1]} [Z_2,Z_1](0,\pi(1),\pi(2))
  - \frac{1}{6} \sum_{\pi \in \mathcal{P}_2[0,1,1]} [Z_1,Z_1,Z_1](0,\pi(1),\pi(2)) \\
  &=\frac{1}{2} [Z_1,Z_2](0,1,2) +\frac{1}{2} ([Z_2,Z_1](0,1,2) + [Z_2,Z_1](0,2,1) )
  -\frac{1}{6} ([Z_1,Z_1,Z_1](0,1,2)+ [Z_1,Z_1,Z_1](0,2,1)) \\
  &= \frac{1}{4} [0,[1,2]] + \frac{1}{4} [[0,1],2] +\frac{1}{4} [[0,2],1] - \frac{1}{6} [[0,1],2] - \frac{1}{6} [[0,2],1] = \frac{1}{6} [0,[1,2]] + \frac{1}{6} [[0,1],2] \;,
 \end{split}
 \end{equation}
 writing $[0,[1,2]]=[A(0),[A(1),A(2)]]=[A(t_0),[A(t_1),A(t_2)]]$ etc. as a short-hand notation.
 This reproduces the term from Eqs.~\eqref{MagnusOm},~\eqref{AppZ}.

\section{Free Lie algebras and Hall bases}
\label{app:Free}

Avoiding formal definitions, the basic concept of a free Lie algebra can be understood in simple terms. For more thorough accounts, see Refs.~\cite{MO99,DeG00}.

A free Lie algebra is a vector space equipped with a function in two arguments $[\cdot,\cdot]$, the commutator. It consists of all nested commutators of the generators $A_1, A_2, \dots$ and all linear combinations thereof. 
In addition to the standard vector space properties, one demands 
bilinearity $[X+Y,Z]=[X,Z]+[Y,Z]$, $[cX,Y]=c[X,Y]$ and anti-symmetry $[X,Y]=-[Y,X]$ of the commutator, 
together with the Jacobi identity $[X,[Y,Z]]+[Y,[Z,X]]+[Z,[X,Y]]=0$.
No further relations hold: Two elements of the free Lie algebra are different if they cannot be transformed into each other with these identities.
In other words, only the minimal relations characteristic for a commutator hold.

Anti-symmetry and the Jacobi identity imply linear dependencies between nested commutators of the generators. In particular, they do not form a vector space basis of the free Lie algebra. 
For three elements $X,Y,Z$ exist 12 commutator combinations 
\begin{equation}\label{app:CommComb}
\begin{split}
&[X,[Y,Z]], \; [X,[Z,Y]], \; [Y,[X,Z]], \; [Y,[Z,X]], \; [Z,[X,Y]], \; [Z,[Y,X]], \\
&[[X,Y],Z], \; [[X,Z],Y], \; [[Y,X],Z], \; [[Y,Z],X], \; [[Z,X],Y], \; [[Z,Y],X] \;.
\end{split}
\end{equation}
Any three of them are linearly dependent,  such that we must select two for a basis.
With this in mind, the Hall basis construction defines a systematic selection rule.
First, define an order ``$<$'' on the generators and nested commutators.
For the generators, set $A_i < A_j$ if $i < j$.
For the commutators, set $[X,Y] < [V,W]$ if $X<V$ or $X=V, Y<W$.
Set generally $X < Y$ if $Y$ is composed out of more commutators than $X$. 
The Hall basis is now defined recursively:
(H1) All generators $A_i$ are in the Hall basis,
(H2) a commutator $[A_i,A_j]$ is in the Hall basis if $A_i < A_j$ (i.e. $i<j$),
(H3) if X,Y,Z are in the Hall basis, so is $[X,[Y,Z]]$ provided that
 $[Y,Z]$ is in the Hall basis and $Y \le X < [Y,Z]$.
 
 To understand rule (H3), observe first that it removes the ambiguity due to anti-symmetry,
 since it enforces $X<Y$ for Hall basis elements $[X,Y]$. 
 Now consider a nested commutator $[X,[Y,Z]]$ from the Hall basis.
 It is $Y \le X < [Y,Z]$ by (H3), and also $Y < Z$. Consequently, $Y<[X,Z]$.
 Both properties rule out most commutators from Eq.~\eqref{app:CommComb}
 apart from $[X,[Y,Z]]$ itself and 
 $[Y,[X,Z]]$, $[Z,[Y,X]]$, $[[Y,X],Z]$.
 If $X=Y$, only $[Y,[X,Z]]$ is non-zero.
 Otherwise, for $Y < X$,  $[Y,[X,Z]]$ violates (H3).
Then, depending on whether $Z \lessgtr [Y,X]$, 
either $[Z,[Y,X]]$ or $[[Y,X],Z]$ fulfills (H3).
If $Z=[X,Y]$, both commutators vanish.
In any case, at most one commutator from Eq.~\eqref{app:CommComb} is a Hall basis element in addition to $[X,[Y,Z]]$. 
This argument implies linear independence of the basis elements,
and can be turned into an inductive proof.
Moreover, rule (H3) amounts to a recursive algorithm to check for membership of the Hall basis.
 
Completeness of the Hall basis can be shown with a similar argumentation.
Based on this, a recursive algorithm can be devised to express commutators $[X,Y]$ as linear combinations of the Hall basis elements.
In Table~\ref{app:HallTable} { we} show the first $23$ Hall basis elements involving the generators $A_1$, $A_2$.
For example, the last commutator $[[A_1,[A_1,A_2]],[A_2,[A_1,A_2]]]$ fulfills (H3) with $X=[A_1,[A_1,A_2]]$, $Y=A_2$, $Z=[A_1,A_2]$.
Another example is to write the four-fold nested commutator $[A_1,[A_2,[A_1,[A_2,A_1]]]] = [[A_1,A_2],[A_1,[A_2,A_1]]]-[A_2,[A_1,[A_1,[A_1,A_2]]]]$ as the unique sum of two elements from the table.
As discussed in Secs.~\ref{sec:Prop},~\ref{sec:CF},
only a small subset of all Hall basis elements needs to be considered for the Magnus expansion or CFET construction.

\begin{table}
  \begin{equation*}
  \begin{gathered}
 A_1, \;
A_2, \;
[A_1,A_2], \;
[A_1,[A_1,A_2]], \;
[A_2,[A_1,A_2]], \\ 
[A_1,[A_1,[A_1,A_2]]], \;
[A_2,[A_1,[A_1,A_2]]], \;
[A_2,[A_2,[A_1,A_2]]], \\ 
[A_1,[A_1,[A_1,[A_1,A_2]]]], \;
[A_2,[A_1,[A_1,[A_1,A_2]]]], \;
[A_2,[A_2,[A_1,[A_1,A_2]]]], \\ 
[A_2,[A_2,[A_2,[A_1,A_2]]]], \;
[[A_1,A_2],[A_1,[A_1,A_2]]], \;
[[A_1,A_2],[A_2,[A_1,A_2]]], \\ 
[A_1,[A_1,[A_1,[A_1,[A_1,A_2]]]]], \;
[A_2,[A_1,[A_1,[A_1,[A_1,A_2]]]]], \;
[A_2,[A_2,[A_1,[A_1,[A_1,A_2]]]]], \\ 
[A_2,[A_2,[A_2,[A_1,[A_1,A_2]]]]], \;
[A_2,[A_2,[A_2,[A_2,[A_1,A_2]]]]], \;
[[A_1,A_2],[A_1,[A_1,[A_1,A_2]]]], \\ 
[[A_1,A_2],[A_2,[A_1,[A_1,A_2]]]], \;
[[A_1,A_2],[A_2,[A_2,[A_1,A_2]]]], \;
[[A_1,[A_1,A_2]],[A_2,[A_1,A_2]]]
  \end{gathered}
\end{equation*}
  \caption{The 23 Hall basis elements with generators $A_1$, $A_2$ and up to $5$ commutators.}
  \label{app:HallTable}
\end{table}

\section{Order conditions for 6th-order CFETs}
\label{app:6th}

The order conditions for 6th-order CFETs can be largely solved by algebraic manipulations.
For 6th-order CFETs with 5 exponentials, one has 7 equations for the 8 coefficients
$f_{1,1}, f_{1,2}, f_{1,3}, f_{2,1}, f_{2,2}, f_{2,3}, f_{3,1}, f_{3,3}$, as follows:
\begin{equation}\label{app:6thOrderConds}
\begin{split}
A_1: \; 1 &=   2 f_{1,1} +2 f_{2,1} + f_{3,1} \\[0.5ex]
A_3: \; 0 &=  2 f_{1,3} +2 f_{2,3} + f_{3,3} \\[0.5ex]
[A_1,A_2]: \; - \frac{1}{6} &= - f_{1,1}f_{1,2} - 2 f_{2,1}f_{1,2} - f_{2,1}f_{2,2} - f_{3,1}f_{1,2} - f_{3,1}f_{2,2}\\[0.5ex]
[A_2,A_3]: \; - \frac{1}{30} &=  \phantom{+} f_{1,2}f_{1,3} +2 f_{1,2}f_{2,3} + f_{1,2}f_{3,3} + f_{2,2}f_{2,3} + f_{2,2}f_{3,3} \\[0.5ex]   
[A_1,[A_1,A_3]]: \; \frac{1}{60} &= 
+\frac{1}{3} f_{1,1}f_{2,1}f_{1,3} -\frac{2}{3} f_{1,1}f_{2,1}f_{2,3} -\frac{1}{3} f_{1,1}f_{2,1}f_{3,3} +\frac{1}{6} f_{1,1}f_{3,1}f_{1,3} -\frac{1}{3} f_{1,1}f_{3,1}f_{2,3} -\frac{1}{6} f_{1,1}f_{3,1}f_{3,3} \\
&\phantom{=} -\frac{1}{3} f_{1,1}^{2}f_{2,3} -\frac{1}{6} f_{1,1}^{2}f_{3,3} +\frac{2}{3} f_{2,1}f_{3,1}f_{1,3} +\frac{1}{6} f_{2,1}f_{3,1}f_{2,3} -\frac{1}{6} f_{2,1}f_{3,1}f_{3,3} +\frac{2}{3} f_{2,1}^{2}f_{1,3} -\frac{1}{6} f_{2,1}^{2}f_{3,3} \\ 
&\phantom{=}+\frac{1}{6} f_{3,1}^{2}f_{1,3} +\frac{1}{6} f_{3,1}^{2}f_{2,3}\\[0.5ex]
[A_2,[A_1,A_2]] : \;  - \frac{1}{60}  &= 
-\frac{1}{3} f_{1,1}f_{1,2}^{2} -1 f_{2,1}f_{1,2}f_{2,2} -1 f_{2,1}f_{1,2}^{2} -\frac{1}{3} f_{2,1}f_{2,2}^{2} -1 f_{3,1}f_{1,2}f_{2,2} -\frac{1}{2} f_{3,1}f_{1,2}^{2} -\frac{1}{2} f_{3,1}f_{2,2}^{2}
\\[0.5ex]
 [A_1,[A_1,[A_1,A_2]]]: \; \frac{1}{360} &= 
 \phantom{+}\frac{1}{3} f_{1,1}f_{2,1}f_{3,1}f_{1,2} +\frac{1}{2} f_{1,1}f_{2,1}f_{3,1}f_{2,2} +\frac{1}{3} f_{1,1}f_{2,1}^{2}f_{1,2} +\frac{1}{3} f_{1,1}f_{2,1}^{2}f_{2,2} +\frac{1}{12} f_{1,1}f_{3,1}^{2}f_{1,2} \\
 &\phantom{=} +\frac{1}{6} f_{1,1}f_{3,1}^{2}f_{2,2} +\frac{1}{3} f_{1,1}^{2}f_{2,1}f_{1,2} +\frac{1}{6} f_{1,1}^{2}f_{2,1}f_{2,2} +\frac{1}{6} f_{1,1}^{2}f_{3,1}f_{1,2} +\frac{1}{6} f_{1,1}^{2}f_{3,1}f_{2,2} +\frac{1}{12} f_{1,1}^{3}f_{1,2} \\
 &\phantom{=} +\frac{1}{12} f_{2,1}f_{3,1}^{2}f_{2,2} +\frac{1}{6} f_{2,1}^{2}f_{3,1}f_{2,2} +\frac{1}{12} f_{2,1}^{3}f_{2,2}
\end{split}
\end{equation}
The order conditions for 6 exponentials have a similar structure, but are too long to be shown here.

Apart from degenerate cases, the order conditions can be reduced to a single polynomial equation.
We consider $f_{1,1}$ as a free parameters.
Then, if $f_{2,1}$ is the solution of $p(f_{1,1},f_{2,1}) =0$ with the polynomial
\begin{equation}
\begin{split}
 p(x,y) = & \phantom{+(} -2 + 30 x - 192 x^2 + 680 x^3 - 1440 x^4 + 1815 x^5 - 1250 x^6 + 
 360 x^7 \\
 & + (18 - 232 x + 1230 x^2 - 3440 x^3 + 5345 x^4 - 4350 x^5 + 1440 x^6) y  \\
& + (- 60 + 650 x - 2740 x^2 + 5655 x^3 - 5710 x^4 + 2250 x^5 ) y^2  
 + (90 - 800 x + 2535 x^2 - 3450 x^3 + 1710 x^4) y^3 \\
 & + ( - 60 + 425 x - 920 x^2 + 630 x^3 ) y^4 + (15 - 80 x + 90 x^2) y^5
\end{split}
\end{equation}
of degree 5 in $y$, the remaining coefficients are given by
\begin{equation}\label{app:RemCoeff}
\begin{split}
 f_{2,2} &= \frac{1 + 5 f_{1,1} (f_{1,1}-1) }{30 (f_{1,1} + f_{2,1}-1) (f_{1,1} + f_{2,1}) ( 2 f_{1,1} + f_{2,1} -1 )}  \;, \;
  f_{1,2} =  \frac{1 - 6 f_{2,2} + 12 f_{1,1} f_{2,2} + 6 f_{2,1} f_{2,2}}{6 (1 - f_{1,1})} \;, \\[1ex]
  f_{1,3} &= \frac{ ( 2 f_{1,1} -1 ) ( 2 f_{1,1} + f_{2,1} -1 ) - 3 f_{2,2} }{ 30 (f_{1,2} (2 f_{1,1} -1 ) ( 2 f_{1,1} + f_{2,1} -1 ) +  f_{2,2} (1 + 8 f_{1,1}^2 + 2 (f_{2,1}-2) f_{2,1} + 
       (8 f_{2,1}-7)f_{1,1} ))} \;, \\[1ex] 
 f_{2,3} &= \frac{ f_{1,1} + 3 f_{1,2} + 4 f_{1,1} f_{2,1} + 2 ( f_{2,1}-1) f_{2,1} + 
 6 f_{2,2} -1}{30 (f_{1,2} (2 f_{1,1}-1) ( 2 f_{1,1} + f_{2,1} -1 ) +  f_{2,2} (1 + 8 f_{1,1}^2 + 2 ( f_{2,1} -2 ) f_{2,1} + 
      (8 f_{2,1}-7)f_{1,1} ) )} \;, \\[1ex]
   f_{3,1} &= 1 - 2 f_{1,1} - 2 f_{2,1} \;, \quad f_{3,3} = -2 f_{1,3} - 2 f_{2,3} \;.
\end{split}
\end{equation}

Several solutions exist with simple explicit expressions for the coefficients, such as the ten solutions shown
in Table~\ref{app:6thCoeff}.
Unfortunately, none of these is competitive with the CFETs from Tables~\ref{6thCFET},~\ref{6thCFETOpt}.
For the CFET CF6:5 with $f_{1,1} = 0.16 = 4/25$, the coefficient $f_{2,1}=0.387524052\dots$ is the single real root of the polynomial $p(x)= -126131602 + 1646347450 x - 7919062500 x^2 + 16950031250 x^3 - 
 15834375000 x^4 + 5498046875 x^5$.

\begin{table}
\begin{center}
  \begin{tabular}{ll}
  \hline
    \multicolumn{2}{l}{6th-order, 5 exponentials} \\\hline\hline\noalign{\smallskip}
    $f_{1,1}= (5-\sqrt{5})/10$ & $f_{2,1}=(23-4\sqrt{5})/60 $ \\\hline\noalign{\smallskip}
    $f_{1,1}= (5+\sqrt{5})/10 $ & $f_{2,1}=(23+4\sqrt{5})/60 $ \\\hline\noalign{\smallskip}
    $f_{1,1}= (65 \pm \sqrt{1005})/90 $ & $f_{2,1}=3/10 $ \\\hline\noalign{\smallskip}
    $f_{1,1}=3/10$ & $f_{2,1}=(553 \pm 3\sqrt{201})/2400$ \\\hline\noalign{\smallskip}
    $f_{1,1}=1$ & $f_{2,1}=(30 \pm \sqrt{\smash[b]{290 \pm 50\sqrt{5}}})/60 $ \\\hline
\end{tabular} 
\end{center}
  \caption{Explicit simple solutions of the order conditions Eq.~\ref{app:6thOrderConds} for 6th-order CFETs with $5$ exponentials.
  The remaining coefficients can be found with Eq.~\eqref{app:RemCoeff}.
  In the last row, all four combinations of the signs are allowed.}
  \label{app:6thCoeff}
\end{table}


\end{document}